\documentclass[twoside]{amsart}
\usepackage{amsmath,amsthm,amsfonts,amscd,amssymb,mathrsfs,latexsym}
\usepackage[all]{xy}

\newcommand{\cD}{{\mathcal D}}
\newcommand{\cE}{{\mathcal E}}

\newcommand{\Ce}{{\mathbb C}}
\renewcommand{\Re}{{\mathbb R}}

\newcommand{\Ze}{{\mathbb Z}}
\newcommand{\Ne}{{\mathbb N}}

\newcommand{\spn}{{\rm span}}
\newcommand{\LCTVS}{{\rm LCTVS}}
\newcommand{\dist}{{\rm dist}}
\newcommand{\mul}{{\rm mul}}
\newcommand{\oq}{{\rm oq}}
\newcommand{\GH}{{\rm GH}}
\newcommand{\rH}{{\rm H}}
\newcommand{\q}{{\rm q}}
\newcommand{\CM}{{\rm CM}}
\newcommand{\CQM}{{\rm CQM}}
\newcommand{\sa}{{\rm sa}}
\newcommand{\Cov}{{\rm Cov}}
\newcommand{\diam}{{\rm diam}}
\newcommand{\Af}{{\rm Af}}
\newcommand{\SUB}{{\rm SUB}}
\newcommand{\Ind}{{\rm Ind}}
\newcommand{\fin}{{\rm fin}}
\newcommand{\tr}{{\rm tr}}

\newcommand{\pa}{\|}

\newcommand{\ie}{{\it i.e. }}

\theoremstyle{definition}
\newtheorem{theorem}{Theorem}[section]
\newtheorem{corollary}[theorem]{Corollary}
\newtheorem{lemma}[theorem]{Lemma}
\newtheorem{proposition}[theorem]{Proposition}
\newtheorem{remark}[theorem]{Remark}
\newtheorem{notation}[theorem]{Notation}
\newtheorem{example}[theorem]{Example}
\newtheorem{question}[theorem]{Question}
\newtheorem{definition}[theorem]{Definition}

\newtheorem*{acknowledgments}{Acknowledgments}

\begin{document}

\title{Order-unit quantum Gromov-Hausdorff distance}
\author{Hanfeng Li}

\address{Department of Mathematics \\
University of Toronto \\
Toronto, ON M5S 3G3, CANADA} \email{hli@fields.toronto.edu}
\date{March 31, 2005}


\begin{abstract}
We introduce a new distance $\dist_{\oq}$ between compact quantum
metric spaces. We show that $\dist_{\oq}$ is Lipschitz equivalent to
Rieffel's distance $\dist_{\q}$, and give criteria for when a
parameterized family of compact quantum metric spaces is
continuous with respect to $\dist_{\oq}$. As applications, we show
that the continuity of a parameterized family of quantum metric
spaces induced by ergodic actions of a fixed compact group is
determined by the multiplicities of the actions, generalizing
Rieffel's work on noncommutative tori and integral coadjoint
orbits of semisimple compact connected Lie groups; we also show
that the $\theta$-deformations of Connes and Landi are continuous
in the parameter $\theta$.
\end{abstract}

\maketitle

\section{Introduction}
\label{intro:sec}

In  \cite{Connes89} Connes initiated
the study of metric spaces in noncommutative setting
in the framework of his spectral triple \cite{Connes94}.
The main ingredient of a spectral triple is a Dirac operator $D$.  On the
one hand, it captures the differential structure by setting $df=[D,
f]$. On the other hand, it enables us to recover the Lipschitz
seminorm $L$, which is usually defined as
\begin{eqnarray} \label{dist to Lip:eq}
L(f):=\sup \{\frac{|f(x)-f(y)|}{\rho(x, y)}:x \neq y\},
\end{eqnarray}
where $\rho$ is the geodesic metric on the Riemannian manifold,
instead by means of $L(f)=\pa [D, f]\pa$, and then one recovers the metric $\rho$  by
\begin{eqnarray} \label{Lip to dist:eq}
\rho(x, y)=\sup_{L(f)\le 1} |f(x)-f(y)|.
\end{eqnarray}
In Section 2 of \cite{Connes89} Connes went further by considering
the (possibly $+\infty$-valued) metric on the state space of the
algebra defined by (\ref{Lip to dist:eq}). Motivated by what
happens to ordinary compact metric spaces, in \cite{Rieffel98b,
Rieffel99b, Rieffel00} Rieffel introduced ``compact quantum metric
spaces'' which requires the metric on the state space to induce
the weak-$*$ topology.
Many interesting examples of compact quantum metric spaces have
been constructed
 \cite{Rieffel98b, Rieffel02, Rieffel03O, Li9}.
Rieffel's theory of compact quantum metric space does not require
$C^*$-algebras, and is set up on  more general spaces, namely
order-unit spaces. Also, one does not need Dirac operators, but
only the seminorm $L$.

Motivated by questions in string theory, in \cite{Rieffel00}
Rieffel also introduced a notion of quantum Gromov-Hausdorff
distance for compact quantum metric spaces, as an analogue of the
Gromov-Hausdorff distance $\dist_{\GH}$ \cite{Gromov81} for ordinary
compact metric spaces. This is defined as a modified ordinary
Gromov-Hausdorff distance for the state-spaces. This distance
$\dist_{\q}$ is a metric on the set $\CQM$ of all isometry
classes of compact quantum metric spaces, and has many nice
properties. Two nontrivial examples of convergence with respect to
$\dist_{\q}$ have been established by Rieffel. One is that the
$n$-dimensional noncommutative tori $T_{\theta}$'s equipped with
the quantum metrics induced from the canonical action of
$\mathbb{T}^n$ are continuous, with the parameter $\theta$ as
$n\times n$ real skew-symmetric matrices
\cite[Theorem 9.2]{Rieffel00}. The other one is that some natural
matrices related to representations of a semisimple compact
connected Lie group converge to integral coadjoint orbits of this
group \cite[Theorem 3.2]{Rieffel01}. In general, it is not easy to
show the continuity of a parameterized family of compact quantum
metric spaces. In particular, the methods used in these two
examples are quite different.

In view of the principle of noncommutative geometry, it may be
more natural to define the quantum distance as a modified
 Gromov-Hausdorff distance for the order-unit spaces (or
 $C^*$-algebras) directly. Under this guidance, we define
an \emph{order-unit quantum Gromov-Hausdorff distance},
$\dist_{\oq}$, as a modified ordinary Gromov-Hausdorff distance for
certain balls in the order-unit spaces (Definition~\ref{dist_oq:def}). We also introduce a variant
$\dist^R_{\oq}$ for the compact quantum metric spaces with radii
bounded above by $R$. Denote by $\CQM^R$  the set of all
isometry classes of these compact quantum metric spaces.
It turns out that these order-unit quantum distances are Lipschitz equivalent to
Rieffel's quantum distance.

\begin{theorem} \label{dist_q=dist_oq:thm}
$\dist_{\q}$ and $\dist_{\oq}$ are Lipschitz equivalent metrics on
$\CQM$, that is
\begin{eqnarray*}
\frac{1}{3}\dist_{\oq}\le \dist_{\q}\le
5\dist_{\oq};
\end{eqnarray*}
while $\dist_{\q}$ and $\dist^R_{\oq}$ are Lipschitz equivalent metrics
on  $\CQM^R$, that is
\begin{eqnarray*}
\frac{1}{2}\dist^R_{\oq}\le \dist_{\q}\le
\frac{5}{2}\dist^R_{\oq}.
\end{eqnarray*}
\end{theorem}

As an advantage of our approach, we can give criteria for when a
parameterized family of compact quantum metric spaces is
continuous with respect to the order-unit quantum distance. We
introduce a notion of continuous fields of compact quantum metric
spaces (Definition~\ref{cont field of CQM:def}), as a concrete way
of saying ``a parameterized family''. This is an analogue of
continuous fields of Banach spaces \cite[Section 10.1]{Dixmier77}.
Roughly speaking, these criteria say that the family is continuous
under quantum distances if and only if continuous sections are
uniformly dense in the balls (the set $\cD(A_t)$
in below)
we use to define the order-unit quantum
distance.

\begin{theorem}  \label{criterion of conv:thm}
Let $(\{(A_t, L_t)\}, \Gamma)$ be a continuous field of compact
quantum metric spaces over a locally compact Hausdorff space $T$.
Let $t_0\in T$, and let $\{f_n\}_{n\in \Ne}$ be a sequence in
$\Gamma$, the space of continuous sections, such that
$(f_n)_{t_0}\in \cD(A_{t_0})$ for each $n\in
\Ne$ and the set $\{(f_n)_{t_0}:n\in \Ne\}$ is dense in
$\cD(A_{t_0})$. Then the following are
equivalent:

(i) $\dist_{\oq}(A_t, A_{t_0})\to 0$ as $t\to t_0$;

(ii) $\dist_{\GH}(\cD(A_t), \cD(A_{t_0}))\to 0$ as $t\to t_0$;

(iii) for any $\varepsilon>0$, there is an $N$ such that the open
$\varepsilon$-balls in $A_t$ centered at $(f_1)_t, {\cdots}, (f_N)_t$
cover $\cD(A_t)$ for all $t$ in some
neighborhood $\mathcal{U}$ of $t_0$.
\end{theorem}

Similar criteria are also given for convergence with respect to
$\dist^R_{\oq}$ (Theorem~\ref{criterion 2 of conv:thm}), which is
useful when the radii of the compact quantum metric spaces are
known to be bounded above by $R$.

An important class of compact quantum metric spaces come from
ergodic actions of compact groups \cite{Rieffel98b}. Let $G$ be a
compact group with a fixed length function $\mathnormal{l}$ given
by $\mathnormal{l}(x)=d(x, e_{G})$, where $x\in G$ and $d$ is a
left-invariant metric on $G$ and $e_G$ is the identity. For an
ergodic action $\alpha$ of $G$ on a unital $C^*$-algebra
$\mathcal{A}$ (\ie the only $\alpha$-invariant elements are the
scalar multiples of the identity of $\mathcal{A}$),
Rieffel proved that the seminorm $L(a)=\sup \{\frac{\pa
\alpha_x(a)-a\pa}{\mathnormal{l}(x)}: x\in G, x\neq e_G\}$ makes
$\mathcal{A}$ into  a compact quantum metric space \cite[Theorem
2.3]{Rieffel98b}. This includes the examples of noncommutative
tori and coadjoint integral orbits mentioned above. In general,
one can talk about ergodic actions of $G$ on complete order-unit
spaces $\overline{A}$. When the action $\alpha$ is finite in the
sense that the multiplicity $\mul(\overline{A}^{\Ce}, \gamma)$ of
every equivalence class of
irreducible representations $\gamma\in \hat{G}$ in the
induced action $\alpha\otimes I$ on
$\overline{A}^{\Ce}=\overline{A}\otimes \Ce$ is finite (which is
always true in $C^*$-algebra case \cite[Proposition 2.1]{HLS81}),
the same construction also makes $\overline{A}$ into a compact
quantum metric spaces.
 Using our criteria for quantum distance convergence
(Theorem~\ref{criterion 2 of conv:thm}) we give a unified proof
for the two examples above about continuity of noncommutative tori
and convergence of matrix algebras to integral coadjoint orbits,
and show in general that a parameterized family of compact quantum
metric spaces induced by ergodic finite actions of $G$ is
continuous with respect to $\dist_{\oq}$ if and only if the
multiplicities of the actions are locally constant:

\begin{theorem} \label{criterion of cont field of action:thm}
Let $\{\alpha_t\}$ be a continuous field of strongly continuous finite ergodic
actions of $G$ on a continuous field of order-unit spaces
$(\{\overline{A_t}\}, \Gamma)$
over a locally compact Hausdorff space $T$.
Then the induced field $(\{(A_t, L_t)\}, \Gamma)$ (for a fixed $\mathnormal{l}$)
 is a continuous field
of compact quantum metric spaces. For any  $t_0\in T$ the
following are equivalent:

(i) $\lim_{t\to t_0}\mul(\overline{A_t}^{\Ce}, \gamma)=
\mul(\overline{A_{t_0}}^{\Ce}, \gamma)$ for all $\gamma \in \hat{G}$;

(ii) $\limsup_{t\to t_0}\mul(\overline{A_t}^{\Ce}, \gamma)\le
\mul(\overline{A_{t_0}}^{\Ce}, \gamma)$ for all $\gamma \in \hat{G}$;

(iii) $\dist_{\oq}(A_t, A_{t_0})\to 0$ as $t\to
    t_0$.
\end{theorem}

In \cite{CL01} Connes and Landi introduced a one-parameter
deformation $S^4_{\theta}$ of the $4$-sphere with the property
that the Hochschild dimension of $S^4_{\theta}$ equals that of
$S^4$. They also considered general $\theta$-deformations, which
was studied further by Connes and Dubois-Violette in \cite{CD01}
(see also \cite{Sitarz01}). In general, the $\theta$-deformation
$M_{\theta}$  of a manifold $M$ equipped with a smooth action of
the $n$-torus $T^n$ is determined by defining the algebra of
smooth functions  $C^{\infty}(M_{\theta})$ as the invariant
subalgebra (under the diagonal action of $T^n$) of the algebra
$C^{\infty}(M\times
T_{\theta}):=C^{\infty}(M)\hat{\otimes}C^{\infty}(T_{\theta})$ of
smooth functions on $M\times T_{\theta}$; here   $\theta$ is a
real skew-symmetric $n\times n$ matrix and $T_{\theta}$ is the
corresponding noncommutative $n$-torus. When $M$ is a compact spin
manifold, Connes and Landi showed that the canonical Dirac
operator $(D, \mathcal{H})$ on $M$ and a deformed anti-unitary
operator $J_{\theta}$ together give a spectral triple for
$C^{\infty}(M_{\theta})$, fitting it into Connes' noncommutative
Riemannian geometry framework \cite{Connes94, Connes96}.

Intuitively, the $\theta$-deformations are continuous in the parameter
$\theta$. Quantum distances provide a concrete way for us to
express the continuity. In \cite{Li9} we showed that when $M$ is connected,
$(C^{\infty}(M_{\theta}))_{\sa}$ equipped with  the seminorm
$L_{\theta}$ determined by the Dirac operator $D$ is a compact
quantum metric space. Denote by $\Theta$ the space of all $n\times
n$ real skew-symmetric matrices. In Section~\ref{Contdeform:sec}
we shall see that there is a natural continuous field of
$C^*$-algebras over $\Theta$
with fibres $C(M_{\theta})$. Denote by $\Gamma^M$ the space of
continuous sections of this field. As another application of our
criteria for quantum distance convergence, we show that
$\theta$-deformations are continuous with respect to $\dist_{\oq}$:

\begin{theorem} \label{theta-deform cont:thm}
Let $M$ be a connected compact spin manifold with a smooth action
of $\mathbb{T}^n$. Then the field $(\{
((C^{\infty}(M_{\theta}))_{\sa}, L_{\theta})\}, (\Gamma^M)_{\sa})$ is a
continuous field of compact quantum metric spaces over $\Theta$.
And
$\dist_{\oq}((C^{\infty}(M_{\theta}))_{\sa}, (C^{\infty}(M_{\theta_0}))_{\sa})\to 0$ as $\theta\to \theta_0$.
\end{theorem}

This paper is organized as follows. In Section~\ref{Prelim:sec} we
review briefly the Gromov-Hausdorff distance for compact metric
spaces and Rieffel's quantum distance for compact quantum metric
spaces. Via  a characterization of state-spaces of compact quantum
metric spaces, a formula for Rieffel's distance $\dist_{\q}$ is given
in Section~\ref{CharSSCQM:sec}.

In Section~\ref{DefOQGH:sec} we define the order-unit
Gromov-Hausdorff distance $\dist_{\oq}$ and prove
Theorem~\ref{dist_q=dist_oq:thm}. One important aspect of the
theory of (quantum) Gromov-Hausdorff distance is the (quantum)
compactness theorem. In Section~\ref{compact:sec} we give a
reformulation of Rieffel's quantum compactness  theorem in terms
of the balls we use to define the order-unit distance. The notion
of continuous fields of compact quantum metric spaces is
introduced in Section~\ref{ConField:sec}. In
Section~\ref{CritConv:sec} we prove our criteria for quantum
distance convergence.

The sections~\ref{Lip&FinMul:sec}-\ref{ContFieldCQMbyG:sec} are
devoted to an extensive study of compact quantum metric spaces
induced by ergodic compact group actions, where we show how
multiplicities of the actions dominate the metric aspect of such
spaces. In Section~\ref{Lip&FinMul:sec} we show that an ergodic
action induces a compact quantum metric space only when the action
is finite. In Section~\ref{comp&bm:sec} we investigate when a
family of compact quantum metric spaces induced from ergodic
actions of a fixed compact group is totally bounded.
Theorem~\ref{criterion of cont field of action:thm} is proved in
Section~\ref{ContFieldCQMbyG:sec}.

Finally, we prove Theorem~\ref{theta-deform cont:thm} in Section~\ref{Contdeform:sec}.

\begin{acknowledgments}
  This is part of my Ph.D. dissertation submitted to UC Berkeley in 2002.
  I am indebted to my advisor, Professor Marc Rieffel, for many helpful discussions, suggestions,
  and for his support throughout my time at Berkeley.
I also thank Thomas Hadfield and Fr\'ed\'eric Latr\'emoli\`ere
  for valuable conversations.
\end{acknowledgments}

\section{Preliminaries}
\label{Prelim:sec}

In this section we review briefly the Gromov-Hausdorff distance
for compact metric spaces \cite{Gromov99, Sakai96, BBI01} and
Rieffel's quantum distance for compact quantum metric spaces
\cite{Rieffel98b, Rieffel99b, Rieffel00, Rieffel01, Rieffel03}.

Let $(X, \rho)$ be a metric space, \ie $\rho$ is a metric on the space
$X$. For any subset $Y\subseteq X$ and $r>0$, let
\begin{eqnarray*}
\mathcal{B}_r(Y)=\{x\in X:\rho(x,y)<r \mbox{ for some } y\in Y\}
\end{eqnarray*}
be the set of points with distance less than $r$ from $Y$.
When $Y=\{x\}$, we also write it as
$\mathcal{B}_r(x)$ and call it the {\it open ball} of radius $r$
centered at $x$.

For nonempty subsets $Y, Z\subseteq X$, we can measure the
distance between $Y$ and $Z$ inside of $X$ by the \emph{Hausdorff distance}
$\dist^{\rho}_{\rH}(Y, Z)$ defined by
\begin{eqnarray*}
 \dist^{\rho}_{\rH}(Y, Z):=\inf \{r>0:Y\subseteq \mathcal{B}_r(Z),\,
 Z\subseteq \mathcal{B}_r(Y)\}.
\end{eqnarray*}
We will also use the notation $\dist^X_{\rH}(Y, Z)$ when there is no
confusion about the metric on $X$.

For any compact metric spaces $X$ and $Y$, Gromov \cite{Gromov81}
introduced the \emph{Gromov-Hausdorff
 distance}, $\dist_{\GH}(X, Y)$, which is defined as
\begin{eqnarray*}
\dist_{\GH}(X, Y)&:=&\inf \{\dist^{Z}_{\rH}(h_X(X), h_Y(Y))|\,
h_X:X\rightarrow Z,\, h_Y:Y\rightarrow Z \mbox{ are }\\ & & \mbox{
isometric embeddings into some metric space } Z\}.
\end{eqnarray*}
It is possible to reduce the space $Z$ in above to be the disjoint union $X\coprod Y$. A
distance $\rho$ on $X\coprod Y$ is said to be {\it admissible} if the
inclusions $X, Y\hookrightarrow X\coprod Y$ are isometric
embeddings. Then it is not difficult to check that
\begin{eqnarray*}
\dist_{\GH}(X, Y)=\inf \{ \dist^{\rho}_{\rH}(X, Y):\rho
\mbox{ is an admissible distance on } X\coprod Y\}.
\end{eqnarray*}

For a compact metric space $(X, \rho)$, we shall denote by
$\diam(X):=\max \{\rho(x, y)|\\x, y\in X\}$
the \emph{diameter} of $X$.
Also let $r_X=\frac{\diam(X)}{2}$ be the \emph{radius} of $X$.
For any $\varepsilon>0$, the \emph{covering number}
$\Cov_{\rho}(X, \varepsilon)$ is defined as the smallest number of open balls of
radius $\varepsilon$ whose union covers $X$.

Denote by $\CM$ the set of isometry classes of compact metric
spaces. One important property of Gromov-Hausdorff distance is the
completeness and compactness theorems by Gromov \cite{Gromov81}:

\begin{theorem}[Gromov's Completeness and Compactness Theorems] \label{GH:thm}
The space $(\CM, \dist_{\GH})$ is a complete metric space. A subset
$\mathcal{S}\subseteq \CM$ is totally bounded
(\ie has compact closure) if and only if

(1) there is a constant $D$ such that $\diam(X, \rho)\le D$ for all
$(X, \rho)\in \mathcal{S}$;

(2) for any $\varepsilon>0$, there exists a constant $K_{\varepsilon}>0$ such that
$\Cov_{\rho}(X, \varepsilon)\le K_{\varepsilon}$ for all $(X, \rho)\in \mathcal{S}$.

\end{theorem}

Next we recall Rieffel's quantum Gromov-Hausdorff distance
$\dist_{\q}$ for compact quantum metric spaces.

Rieffel has found that the right framework for compact quantum
metric spaces is that of order-unit spaces. There is an abstract
characterization of order-unit spaces due to Kadison
\cite{Kadison51, Alfsen71}. An \emph{order-unit space} is a real
partially ordered vector space, $A$, with a distinguished element
e (the order unit) satisfying:

\noindent
1) (Order unit property) For each $a\in A$ there is an $r\in \Re$ such
that $a\le re$;

\noindent
2) (Archimedean property) For $a\in A$, if $a\le re$ for all $r\in \Re$
with $r>0$, then $a\le 0$.

On an order-unit space $(A, e)$, we can define a norm as
\begin{eqnarray*}
\pa a\pa=\inf \{r\in \Re: -re\le a\le re\}.
\end{eqnarray*}
Then $A$ becomes a normed vector space and we can consider its
dual, $A'$, consisting of the bounded linear functionals, equipped
with the dual norm $\pa \cdot \pa'$.

By a {\it state} of an order-unit space $(A, e)$, we mean a $\mu
\in A'$ such that $\mu(e)={\pa \mu\pa'} =1$. States are
automatically positive. Denote the set of all states of $A$ by
$S(A)$. It is a compact convex subset of $A'$ under the
weak-$*$ topology. Kadison's basic representation theorem
\cite{Alfsen71} says that the natural pairing between $A$ and
$S(A)$ induces an isometric order isomorphism of $A$ onto a dense
subspace of the space $\Af_{\Re}(S(A))$ of all affine $\Re$-valued
continuous functions on $S(A)$, equipped with the supremum norm
and the usual order on functions.

For an order-unit space $(A, e)$ and a seminorm $L$ on $A$,
we can define an ordinary
metric, $\rho_L$, on $S(A)$ (which may take value $+\infty$) by
(\ref{Lip to dist:eq}).
We say that $L$ is a \emph{Lipschitz seminorm}
on $A$ if it satisfies:

1) For $a\in A$, we have $L(a)=0$ if and only if $a\in \Re e$.

\noindent We call $L$ a \emph{Lip-norm}, and call the pair $(A,
L)$ a \emph{compact quantum metric space} \cite[Definitions 2.1,
2.2]{Rieffel00} if $L$ satisfies further:

2) The topology on $S(A)$ induced by the metric $\rho_L$ is the
    weak-$*$ topology.

\noindent
The \emph{diameter} $\diam(A)$, the \emph{radius} $r_A$, and the
\emph{covering number} $\Cov(A, \varepsilon)$ of $(A, L)$ are defined to
be those of $(S(A), \rho_L)$.

Let $(A, e)$ be an order-unit space with a Lipschitz seminorm $L$. Then
$L$ and $\pa \cdot \pa$ induce norms ${\tilde L}$ and $\pa
\cdot\pa^{\sim}$ respectively
on the quotient space ${\tilde A}=A/\Re e$.
The dual of $({\tilde A},\pa \cdot \pa^{\sim})$ is
exactly $A'^0=\{ \lambda \in A':\lambda(e)=0\}$.
Now ${\tilde L}$ induces a dual seminorm $L'$ on $A'^0$, which
may take value $+\infty$. The metric on $S(A)$ induced by
(\ref{Lip to dist:eq}) is related to $L'$ by:
\begin{eqnarray} \label{metric and dual:eq}
\rho_L(\mu, \nu)=L'(\mu-\nu)
\end{eqnarray}
for all $\mu, \nu \in S(A)$.

\begin{notation} \label{ball:notation}
For any $r\ge 0$, let
\begin{eqnarray*}
\cD_{r}(A):=\{a\in A: L(a)\le 1, \pa a\pa \le
r\}.
\end{eqnarray*}
When $L$ is a Lip-norm on $A$, set
\begin{eqnarray*}
\cD(A):=\cD_{r_A}(A).
\end{eqnarray*}
\end{notation}

\begin{proposition}\cite[Proposition 1.6, Theorem 1.9]{Rieffel98b}
\label{criterion of Lip:prop}
Let $(A, e)$ be an order-unit space with a Lipschitz seminorm $L$.
Then $L$ is a Lip-norm if and only if

\noindent \, \, \, \, \, \, \,\,
(1) there is a constant $K\ge 0$ such that $L'\le K\pa \cdot \pa'$ on $A'^0$;

\noindent \, \, \, \, \,
or
(1') there is a constant $K\ge 0$ such that $\pa \cdot\pa^{\sim}\le K
\tilde{L}$ on $\tilde{A}$;

\noindent
and \, \, \, \,
(2) for any $r\ge 0$, the  ball $\cD_{r}(A)$
is totally bounded in $A$ for $\pa \cdot \pa$;

\noindent \, \, \, \, \,
or
(2') for some $r> 0$, the  ball $\cD_{r}(A)$
is totally bounded in $A$ for $\pa \cdot \pa$.

In this event, the minimal $K$ is exactly $r_A$.
\end{proposition}

Let $A$ be an order-unit space. By a \emph{quotient} $(\pi, B)$ of $A$, we mean
an order-unit space $B$ and a surjective linear positive map $\pi:A\rightarrow B$
preserving the order-unit. Via the dual map $\pi':B'\rightarrow A'$, one may identify
$S(B)$ with a closed convex subset of $S(A)$.
This gives a bijection between isomorphism classes
of quotients of $A$ and closed convex subsets of $S(A)$
\cite[Proposition 3.6]{Rieffel00}.
If $L$ is a Lip-norm on $A$, then the quotient seminorm $L_B$ on $B$, defined by
\begin{eqnarray*}
L_B(b):=\inf \{ L(a): \pi (a)=b\}
\end{eqnarray*}
is a Lip-norm on $B$, and $\pi'|_{S(B)}:S(B)\rightarrow S(A)$
is an isometry for the corresponding metrics
$\rho_L$ and $\rho_{L_B}$ \cite[Proposition 3.1]{Rieffel00}.

Let $(A, L_A)$ and $(B, L_B)$ be compact quantum metric spaces.
The direct sum $A\oplus B$, of vector spaces, with $(e_A, e_B)$ as
order-unit, and with the natural order structure is also an
order-unit space. We call a Lip-norm $L$ on $A\oplus B$ {\it
admissible} if it induces $L_A$ and $L_B$ under the natural
quotient maps $A\oplus B\rightarrow A$ and $A\oplus B\rightarrow
B$. Rieffel's \emph{quantum Gromov-Hausdorff distance} $\dist_{\q}(A,
B)$ \cite[Definition 4.2]{Rieffel00} is defined by
\begin{eqnarray*}
\dist_{\q}(A,B)= \inf \{\dist^{\rho_L}_{\rH}(S(A), S(B)): L \mbox{ is an
admissible Lip-norm on } A\oplus B\}.
\end{eqnarray*}

Let $(A, L)$ be a compact quantum metric space. Let $\bar{A}$ be
the completion of $A$ for $\pa \cdot \pa$. Define a seminorm,
$\bar{L}$, on $\bar{A}$ (which may take value $+\infty$) by
\begin{eqnarray*}
\bar{L}(b):=\inf \{\liminf_{n\to \infty}L(a_n): a_n\in A, \lim_{n\to
\infty}a_n=b\}.
\end{eqnarray*}
The \emph{closure} of $L$, denoted by $L^c_A$, is defined as the
restriction of $\bar{L}$ to the subspace
\begin{eqnarray*}
A^c:=\{b\in \bar{A}:\bar{L}(b)<\infty\}.
\end{eqnarray*}
Then $L^c$ is a Lip-norm on $A^c$, and $\rho_{L}=\rho_{L^c}$ on $S(A)=S(A^c)$ \cite[Theorem 4.2, Proposition 4.4]{Rieffel99b}.
Identify $\bar{A}$ with $\Af_{\Re}(S(A))$. Then
$A^c$ is exactly the space of Lipschitz functions in $\Af_{\Re}(S(A))$, and $L^c$ is just the Lipschitz seminorm defined
by (\ref{dist to Lip:eq}) \cite[Proposition 6.1]{Rieffel00}.
We say that $L$ is \emph{closed} if $L$ equals its closure.

Let $(A, L_A)$ and $(B, L_B)$ be compact quantum metric spaces. By
an \emph{isometry} from $(A, L_A)$ to $(B, L_B)$ we mean an order
isomorphism $\varphi$ from $A^c$ onto $B^c$ such that
$L^c_A=L^c_B\circ \varphi$. The isometries from $(A, L_A)$ to $(B,
L_B)$ are in natural bijective correspondence with the affine
isometries from $(S(B), \rho_{L_B})$ onto $(S(A), \rho_{L_A})$
through $\varphi\mapsto \varphi'|_{S(B)}$ \cite[Corollary
6.4]{Rieffel00}.

Denote by $\CQM$ the set of isometry classes of compact quantum
metric spaces. Rieffel also proved a quantum version of Gromov's
completeness and compactness theorems \cite[Theorems 12.11 and
13.5]{Rieffel00}:

\begin{theorem}[Rieffel's Quantum Completeness and Compactness Theorems]
 \label{QGH:thm}
The space $(\CQM, \dist_{\q})$ is a complete metric space. A
subset $\mathcal{S}\subseteq \CQM$ is totally bounded if
and only if

(1) there is a constant $D$ such that $\diam(A, L)\le D$ for all $(A,
L)\in \mathcal{S}$;

(2) for any $\varepsilon>0$, there exists a constant $K_{\varepsilon}>0$
such that $\Cov(A, \varepsilon)\le K_{\varepsilon}$ for all $(A, L)\in
\mathcal{S}$.
\end{theorem}

\section{A characterization of state-spaces of compact quantum metric spaces}
\label{CharSSCQM:sec}

In this section we give a characterization of state-spaces of
compact quantum metric spaces in Proposition~\ref{criterion of
SSQ:prop}, and use it to give a formula for Rieffel's $\dist_{\q}$ in
Proposition~\ref{new dist_q:prop}.

Proposition 5.7 and Corollary 6.4 in \cite{Rieffel00} tell us that
for compact quantum metric spaces $(B_i, L_i)$, $i=1,2$, if their
state-spaces are affinely isometrically embedded into the
state-space $S(A)$ of some other compact quantum metric space $(A,
L)$, then
\begin{eqnarray*}
 \dist_{\q}(B_1, B_2)\le \dist^{S(A)}_{\rH}(S(B_1), S(B_2)).
\end{eqnarray*}
This provides a powerful way of getting upper bounds for $\dist_{\q}(B_1, B_2)$.
In practice, it is quite easy to embed the
state-space of a quantum metric space into some other compact metric
space. So we need to find out what kind of compact metric spaces can
be the state-space of a compact quantum metric space.

Throughout the rest of this section, locally convex topological
vector spaces ({\LCTVS}) will all be Hausdorff. Let $\mathfrak{X}$
be a compact convex subset of a {\LCTVS} $V$ over $\Re$. Then
$(\Af_{\Re}(\mathfrak{X}), 1_{\mathfrak{X}})$ is an order-unit
space. For each $\mu \in \mathfrak{X}$, the evaluation at $\mu$
induces a linear function $\sigma(\mu)$ on $\Af_{\Re}(\mathfrak{X})$.
Clearly
\begin{eqnarray*}
(\sigma(\mu))(1_{\mathfrak{X}})=1=\pa \sigma(\mu)\pa.
\end{eqnarray*}
So $\sigma(\mu)$ is a state of $\Af_{\Re}(\mathfrak{X})$. This defines
an affine map $\sigma:\mathfrak{X}\rightarrow S(\Af_{\Re}(\mathfrak{X}))$.
Let $\rho$ be a metric on $\mathfrak{X}$. We say that $\rho$ is
\emph{midpoint-balanced} \cite[Definition 9.3]{Rieffel99b}
if for any $\mu, \nu, \mu', \nu' \in \mathfrak{X}$ with
$\frac{\mu+\nu'}{2}=\frac{\mu'+\nu}{2}$, we have $\rho(\mu,
\nu)=\rho(\mu', \nu')$. We say that $\rho$ is \emph{convex}
if for
any $\mu, \nu, \mu', \nu' \in \mathfrak{X}$ and $0\le t\le 1$, we
have
\begin{eqnarray*}
\rho(t\mu+(1-t)\mu',t\nu+(1-t)\nu')\le t\rho(\mu,
\nu)+(1-t)\rho(\mu', \nu').
\end{eqnarray*}

\begin{proposition} \label{criterion of SSQ:prop}
Let $\mathfrak{X}$ be a compact convex subset of a {\LCTVS} $V$,
and let $\rho$ be a metric on $\mathfrak{X}$ compatible with the
topology. Then $(\mathfrak{X}, \rho)$ is affinely isometric to
$(S(A), \rho_L)$ for some compact quantum metric space $(A, L)$ if
and only if the metric $\rho$ is convex and midpoint-balanced. In
this event, the closed compact quantum metric space is
$(\Af_{\Re}(\mathfrak{X})_L, L_{\rho})$, unique up to isometry,
where $\Af_{\Re}(\mathfrak{X})_L$ is the space of Lipschitz
functions in $\Af_{\Re}(\mathfrak{X})$ and $L_{\rho}$ is the
Lipschitz seminorm defined by (\ref{dist to Lip:eq}).
\end{proposition}
\begin{proof}
Assume that $(\mathfrak{X}, \rho)$ is affinely isometric to
$(S(A), \rho_L)$ for some compact quantum metric space $(A, L)$.
It is easy to check directly from (\ref{Lip to dist:eq}) that the
metric $\rho_L$ and hence $\rho$ are convex and midpoint-balanced.

Conversely, assume that the metric $\rho$ is convex and balanced.
Elements in the dual $V'$ separate
the points in V by the Hahn-Banach theorem.
Since the restrictions of elements in $V'$ to $\mathfrak{X}$
    are all in $\Af_{\Re}(\mathfrak{X})$, we see that functions in
    $\Af_{\Re}(\mathfrak{X})$ separate the points of $\mathfrak{X}$.
Theorem II.2.1 in \cite{Alfsen71} tells us that $\sigma$ is a
homeomorphic embedding of $\mathfrak{X}$ into
$S(\Af_{\Re}(\mathfrak{X}))$, and that $\sigma(\mathfrak{X})$
contains the set of extreme points of $S(\Af_{\Re}(\mathfrak{X}))$.
Since $\sigma(\mathfrak{X})$ is convex and closed, we see that
$\sigma$ is surjective. Hence we may identify $\mathfrak{X}$ and
$S(\Af_{\Re}(\mathfrak{X}))$.  By \cite[Lemma 2.1]{Rieffel99b} we
have
$(\Af_{\Re}(\mathfrak{X}))'^0=\Re(S(\Af_{\Re}(\mathfrak{X}))-S(\Af_{\Re}(\mathfrak{X})))=
\Re(\mathfrak{X}-\mathfrak{X})$ (see the
  discussion preceding Notation~\ref{ball:notation}). By \cite[Theorem
9.7]{Rieffel99b} there is a norm $M$ on
$(\Af_{\Re}(\mathfrak{X}))'^0=\Re(\mathfrak{X}-\mathfrak{X})$ such
that $\rho(\mu, \nu)=M(\mu-\nu)$ for all $\mu, \nu \in
\mathfrak{X}$.  Then
\cite[Theorem 9.8]{Rieffel99b} (see also the discussion right
after the proof of Proposition 1.1 in \cite{Rieffel02}) asserts
that $(\Af_{\Re}(\mathfrak{X})_L, L_{\rho})$ is a closed compact
quantum metric space and $(\mathfrak{X}, \rho)$ is its
state-space. The uniqueness of such a closed compact quantum
metric space follows from \cite[Corollary 6.4]{Rieffel00}.
\end{proof}

Consequently we have the following description of the quantum
distance $\dist_{\q}$:

\begin{proposition} \label{new dist_q:prop}
Let $(A, L_A)$ and $(B, L_B)$ be compact quantum metric spaces.
Then we have
\begin{eqnarray*}
\dist_{\q}(A, B)=\inf\{\dist^V_{\rH}(h_A(S(A)), h_B(S(B))):
h_A \mbox{ and }h_B \mbox{ are affine isometric }\\
\mbox{ embeddings of } S(A) \mbox{ and }S(B) \mbox{ into some real normed space }V\}.
\end{eqnarray*}
\end{proposition}
\begin{proof} Denote the right hand side of the above equation by
  $\dist'_{\q}(A, B)$. For any admissible Lip-norm $L$ on $A\oplus B$
  let $V=(A\oplus B)'^{0}$ equipped with the norm $L'$ (see the
  discussion preceding Notation~\ref{ball:notation}). Pick an element
  $p$ in $S(A\oplus B)$, and let $\varphi: S(A\oplus B)\to V$
  be the translation $x\mapsto x-p$. Then $\varphi$ is an affine
  isometric embedding from $(S(A\oplus B), \rho_L)$ to $V$ according
 to (\ref{metric and dual:eq}). Hence $\dist'_{\q}(A,B)\le
  \dist^V_{\rH}(\varphi(S(A)), \varphi(S(B)))=\dist^{\rho_L}_{\rH}(S(A), S(B))$.
  Thus $\dist'_{\q}(A, B)\le  \dist_{\q}(A, B)$.

Now let $V, h_A$ and $h_B$ be as in Proposition~\ref{new
  dist_q:prop}. Let $\mathfrak{X}$ be the convex hull of
  $h_A(S(A))\cup h_B(S(B))$. Clearly $\mathfrak{X}$ equipped
  with the distance induced from the norm in $V$ is
  compact, and hence is the state-space of some compact quantum metric
 space $(C, L_C)$ by Proposition~\ref{criterion of SSQ:prop}.
 Therefore $\dist_{\q}(A, B)\le \dist^V_{\rH}(h_A(S(A)),
  h_B(S(B)))$ by \cite[Proposition 5.7, Corollary 6.4]{Rieffel00}.
  Consequently $\dist_{\q}(A, B)\le \dist'_{\q}(A, B)$.
\end{proof}


\section{Definition of the order-unit quantum Gromov-Hausdorff\\ distance}
\label{DefOQGH:sec}

In this section we define the order-unit Gromov-Hausdorff distance and prove
Theorem~\ref{dist_q=dist_oq:thm}.

Rieffel's definition of quantum Gromov-Hausdorff distance is a
modified ordinary Gromov-Hausdorff distance for the state-spaces.
In the view of Noncommutative Geometry, whose principle is the
duality between ordinary spaces and appropriate vector spaces of
functions over the spaces, it may be more natural to do everything
on the vector spaces of functions directly, avoiding referring
back to the state-spaces. So it may be more natural to measure the
ordinary Gromov-Hausdorff distance for the vector spaces of
functions directly. But the order-unit spaces of functions are not
compact, so we can not apply the ordinary Gromov-Hausdorff distance
to them. One way to get around this difficulty is to consider some
core of the vector spaces of functions which captures all the
information of the order-unit spaces. One natural choice is the
unit ball. But, unless the order-unit space is finite dimensional,
the unit ball is not compact either. It also does not remember the
Lip-norm. Now comes the candidate, $\cD(A)$ (see
Notation~\ref{ball:notation}) for closed $(A, L_A)$. When $r_A>0$,
$\cD(A)$ is absorbing, \ie for every $a\in A$ there
is some $\varepsilon>0$ such that $\lambda a\in \cD(A)$ for all
$0\le \lambda <\varepsilon$. Thus $\cD(A)$ equipped with the metric induced by the norm of $A$
encodes the normed space structure of $A$. It also captures the
Lip-norm:

\begin{lemma} \label{ball to Lip:lemma}
Let $(A, L)$ be a closed compact quantum metric space. Then for
any $R\ge r_A$ we have
\begin{eqnarray} \label{ball to Lip:eq}
\{a\in A:L(a)\le 1\}=\Re e_A+\cD_{ R}(A).
\end{eqnarray}
Conversely, let $(B, e_B)$ be an order-unit space, and let $X$ be
a balanced (\ie $\lambda x\in X$ for all $x\in X$ and $\lambda \in
\Re$ with $|\lambda |\le 1$), absorbing (\ie $\{\lambda x:\lambda\in \Re_{+}, x\in X\}=B$),
compact convex subset of $B$ (under the order-unit norm topology).
Let $R$ be
the radius of $X$. If $X=\{b\in (X+\Re e_B):\pa b\pa \le R\}$,
then there is a unique closed Lip-norm $L$ on $B$ such that
$X=\cD_{ R}(B)$. In this case $L$ is also characterized as
the unique seminorm on $B$ satisfying $X+\Re e_B=\{b\in B: L(b)\le
1\}$.
\end{lemma}
\begin{proof} (\ref{ball to Lip:eq})
follows directly from Proposition~\ref{criterion of Lip:prop}. Now
let $X$ be as in Lemma~\ref{ball to Lip:lemma}. Then clearly
$X+\Re e_B$ is also a balanced absorbing convex set.
Since $X$ is compact, $X+\Re e_B$ is closed. Let $L$ be the
Minkowski functional \cite[Theorem 37.4]{Berberian74}
corresponding to $X+\Re e_B$, \ie the unique seminorm on $B$
satisfying that $X+\Re e_B=\{b\in B: L(b)\le 1\}$. Clearly
$L(e_B)$=0. Suppose that $L(b)=0$. Then for any $n\in \Ne$ we have
$nb\in X+\Re e_B$. Thus there exist $x_n\in X$ and $\lambda_n\in
\Re$ such that $nb=x_n+\lambda_ne_B$. Since $\pa x_n\pa\le R$, we
have $\pa \tilde{b}\pa^{\sim}=\pa \frac{1}{n}x_n\pa^{\sim}\le
\frac{1}{n}R$ in $\tilde{B}=B/ \Re e_B$. Thus $\pa
\tilde{b}\pa^{\sim}=0$, and hence $b\in \Re e_B$. Therefore $L$ is
a Lipschitz seminorm on $B$. Clearly the condition (1') in
Proposition~\ref{criterion of Lip:prop} is satisfied with $K=R$.
The assumption $X=\{b\in (X+\Re e_B):\pa b\pa \le R\}$ means that
$X=\cD_{ R}(B)$. Note that $R>0$ since $X$ is absorbing.
Thus the condition (2') in
Proposition~\ref{criterion of Lip:prop} is also satisfied with
$r=R$. By Proposition~\ref{criterion of Lip:prop}
$L$ is a Lip-norm on $B$, and $r_B\le R$. Since $X+\Re e_B$ is
closed, $L$ is closed. The uniqueness of such a closed Lip-norm
follows from  (\ref{ball to Lip:eq}).
\end{proof}

Most importantly, $\cD(A)$ is compact with the
distance induced from the norm on $A$ by
Proposition~\ref{criterion of Lip:prop}. So we can use it to
redefine the quantum Gromov-Hausdorff distance. There is one
subtle point: we do not know whether $\cD(A)$ remembers the order-unit $e_A$ or not (see
Remark~\ref{dist'_oq:ques}). We shall come back to this point later.

Now the question is what kind of modified Gromov-Hausdorff
distance we should put on $\cD(A)$. Certainly this
modified Gromov-Hausdorff distance should reflect the convex
structure on $\cD(A)$. If we look at the definition
of $\dist_{\GH}$ in Section~\ref{Prelim:sec}, one immediate choice
for the modified distance is $\inf\{\dist^V_{\rH}(h_A(\cD(A)), h_B(\cD(B)))\}$, where the infimum runs
over affine isometric embeddings $h_A$ and $h_B$ of $\cD(A)$ and $\cD(B)$ into some real normed space
$V$. On the other hand, notice that $\cD(A)$ is the
state-space of some compact quantum metric space $(A, L_A)'$
according to Proposition~\ref{criterion of SSQ:prop}. So we may
try to use Rieffel's quantum distance for $(A, L_A)'$ and $(B,
L_B)'$. Proposition~\ref{new dist_q:prop} tells us that these two
possible definitions agree. Notice that when $r_A>0$ we can extend
$h_A$ uniquely to an affine isometric embedding of $A$ into $V$.
When $r_A=0$, the space $A$ is one-dimensional, so we can also
extend $h_A$ to $A$ (by enlarging $V$ if $V=\{0\}$). Therefore the
infimum actually runs over affine isometric embeddings $h_A$ and
$h_B$ of $A$ and $B$ into real normed spaces $V$. These embeddings
may not be linear since $h_A(0_A)$ and $h_B(0_B)$ need not be
$0_V$. But we can always assume that $h_A$ is linear by composing
both $h_A$ and $h_B$ with the translation $x\mapsto x-h_A(0_A)$ in
$V$. To makes things easier, we choose to require both $h_A$ and
$h_B$ to be linear. Since we do not know whether $\cD(A)$ remembers the order-unit $e_A$ or not (see
Remark~\ref{dist'_oq:ques}),
 we need to consider also $\pa
h_A(r_Ae_A)-h_B(r_Be_B)\pa$. Now we get to:

\begin{definition} \label{dist_oq:def}
Let $(A, L_A)$ and $(B, L_B)$ be compact quantum metric spaces. We
define the {\it order-unit quantum Gromov-Hausdorff distance}
between them, denoted by $\dist_{\oq}(A, B)$, by
\begin{eqnarray*}
\dist_{\oq}(A, B)
 :=\inf \{\max (\dist^V_{\rH}(h_A(\cD(A)), h_B(\cD(B))),
\pa h_A(r_Ae_A)-h_B(r_Be_B)\pa)\},
\end{eqnarray*}
and, for $R\ge 0$, the
 {\it $R$-order-unit quantum
Gromov-Hausdorff distance} between them, denoted by
$\dist^R_{\oq}(A, B)$, by
\begin{eqnarray*}
\dist^R_{\oq}(A, B)
 := \inf \{\max
(\dist^V_{\rH}(h_A(\cD_{ R}(A)), h_B(\cD_{ R}(B))),
\pa h_A(Re_A)-h_B(Re_B)\pa) \},
\end{eqnarray*}
where the infima are taken over all triples $(V, h_A, h_B)$ consisting of a real normed space
$V$ and linear isometric embeddings $h_A:A\rightarrow V$ and $h_B:B\rightarrow V$.
\end{definition}

\begin{remark} \label{dist_oq:remark}
(1) To simply the notation, usually we shall identify $A$ and $B$ with
    their images $h_A(A)$ and $h_B(B)$ respectively, and just say that
    $V$ is a normed space containing both $A$ and $B$;

(2) See the discussion preceding Theorem~\ref{criterion 2 of conv:thm} for the motivation
of introducing $\dist^R_{\oq}$;

(3) We choose to use the terms $\pa
h_A(r_Ae_A)-h_B(r_Be_B)\pa$ and $\pa h_A(Re_A)-h_B(Re_B)\pa$ to take
care of the order-units. As another choice, one may also
omit these terms and require $h_A(e_A)=h_B(e_B)$
in Definition~\ref{dist_oq:def}. Denote the resulting distances by
$\dist^*_{\oq}$ and $\dist^{R*}_{\oq}$.
It is easy to see that $\dist_{\oq}\le \dist^*_{\oq}
$ and $\dist^R_{\oq}\le \dist^{R*}_{\oq}$. One may also check that the
proofs of Propositions~\ref{dist_oq<dist_q:prop},
\ref{dist_q<dist_oq:prop}, and Theorem~\ref{dist_q=dist_oq:thm}
 hold with $\dist_{\oq}$ and $\dist^R_{\oq}$ replaced by $\dist^*_{\oq}$ and
$\dist^{R*}_{\oq}$;

(4) For any ordinary compact metric space $(X, \rho)$, let $A_X$
be the space of Lipschitz $\Re$-valued functions on $X$ and let
$L_{\rho}$ be the Lipschitz seminorm defined by (\ref{dist to
Lip:eq}). Then $(A_X, L_{\rho})$ is a closed compact quantum
metric space, called  the \emph{associated compact quantum metric
space} of $(X, \rho)$. For any compact metric spaces $(X, \rho_X)$
and $(Y, \rho_Y)$, by \cite[Proposition 4.7]{Rieffel00} and
Theorem~\ref{dist_q=dist_oq:thm} we have $\dist_{\oq}(A_X, A_Y)\le
3\dist_{\q}(A_X, A_Y)\le 3\dist_{\GH}(X, Y)$. Using \cite[Theorem
13.16]{Rieffel00} and Theorems~\ref{GH:thm}, \ref{QGH:thm}, and
\ref{dist_q=dist_oq:thm}, one can see that the distance $(X,
Y)\mapsto \dist_{\oq}(A_X, A_Y)$ determines the same topology on
$\CM$ as does $\dist_{\GH}$.
\end{remark}

As in the discussion for Gromov-Hausdorff
distance in Section~\ref{Prelim:sec},
it suffices to have $V$ to be $A\oplus B$ (equipped
with certain norms) in Definition~\ref{dist_oq:def}.
To this end, for any normed
spaces $V$ and $W$ we call
a norm $\pa \cdot \pa_{V\oplus W}$ on $V\oplus W$ \emph{admissible} if
it extends the norms on $V$ and $W$.

\begin{proposition} \label{V=A+B:prop}
Let $(A, L_A)$ and $(B, L_B)$ be compact quantum metric spaces.
Then
\begin{eqnarray*}
\dist_{\oq}(A, B)= \inf \{\max (\dist^{A\oplus B}_{\rH}(\cD(A), \cD(B)),
\pa r_Ae_A-r_Be_B\pa_{A\oplus B})\},
\end{eqnarray*}
and, for any $R\ge 0$,
\begin{eqnarray*}
\dist^R_{\oq}(A, B)= \inf \{\max (\dist^{A\oplus B}_{\rH}(\cD_{
R}(A), \cD_{ R}(B)),
\pa Re_A-Re_B\pa_{A\oplus B})\},
\end{eqnarray*}
where the infima are taken over all admissible norms $\pa \cdot \pa_{A\oplus B}$ on $A\oplus B$.
\end{proposition}
\begin{proof} We prove the case of $\dist_{\oq}(A, B)$. That of
  $\dist^R_{\oq}$ is similar.
The proof here could be thought of as a dual of Example 5.6 and
Proposition 5.7 in \cite{Rieffel00}. Let $(V, h_A, h_B)$ be as in
Definition~\ref{dist_oq:def}. Let $1>\varepsilon>0$ be given. We will
construct an admissible norm on $V\oplus V$ such that the two
copies of $V$ are $\varepsilon$-close to each other, \ie $\pa
(v,-v)\pa_{V\oplus V}\le \varepsilon \pa v\pa$. Clearly $\pa
(u,v)\pa_{V\oplus V}:=\max(\pa u+v\pa, \varepsilon\pa u\pa,
\varepsilon\pa v\pa)$ satisfies the requirement. Now we identify
$A\oplus B$ with the subspace $h_A(A)\oplus h_B(B)$ of $V\oplus
V$. Then the induced norm on $A\oplus B$ is admissible. And
\begin{eqnarray*}
\dist^{A\oplus B}_{\rH}(\cD(A), \cD(B))
&\le & \dist^V_{\rH}(h_A(\cD(A)),
h_B(\cD(B)))+ \\
& &\dist^{V\oplus V}_{\rH}((h_B(\cD(B)), 0), (0, h_B(\cD(B)))\\
&\le &\dist^V_{\rH}(h_A(\cD(A)), h_B(\cD(B)))+\varepsilon r_B.
\end{eqnarray*}
Similarly, $\pa r_Ae_A-r_Be_B\pa_{A\oplus B} \le \pa
h_A(r_Ae_A)-h_B(r_Be_B)\pa_V+\varepsilon r_B$. This gives the desired
result.
\end{proof}

We start to prove Theorem~\ref{dist_q=dist_oq:thm}.
We prove the triangle inequality first. For this we need
the amalgamation of normed spaces:

\begin{lemma} \label{amal:lemma}
Let $\varphi_j:A \hookrightarrow B_j$ be  linear isometric embeddings
of normed spaces (over $\Re$ or $\Ce$) for $j\in J$, where
$J$ is an index set. Then there is a
normed space $C$ and linear isometric embeddings
$\psi_j:B_j\hookrightarrow C$ such that $\psi_j\circ
\varphi_j=\psi_k\circ \varphi_k$ for all $j, k\in J$.
\end{lemma}
\begin{proof}
Let $\pa \cdot \pa_1$ be the $L_1$-norm on $\oplus_{j\in J}B_j$, \ie
$\pa (u_j)\pa_1=\sum_{j\in J}\pa u_j\pa$. Let
\begin{eqnarray*}
W=\{(u_j):
u_j\in \varphi_j(A) \mbox{ for all } j\in J,
\mbox{ and }\sum_{j\in J}(\varphi_j)^{-1}(u_j)=0\},
\end{eqnarray*}
which is a linear subspace of $\oplus_{j\in J}B_j$.
Let $q:\oplus_{j\in J}B_j\to
(\oplus_{j\in J}B_j)/W$ be the quotient map, and let $\psi_j:B_j\to
(\oplus_{j\in J}B_j)/W$ be the composition of $B_j\to\oplus_{j\in J}B_j$ and
$q$. Then clearly $\psi_j\circ
\varphi_j=\psi_k\circ \varphi_k$ for all $j, k\in J$, and $\psi_j$ is contractive.
For any $u\in B_k$ and $(\varphi_j(v_j))\in W$ we have
\begin{eqnarray*}
\pa u+(\varphi_j(v_j))\pa_1&=&\pa u+\varphi_k(v_k)\pa +\sum_{j\in J,
  j\neq k}\pa v_j\pa\\
&=&\pa u-\varphi_k(\sum_{j\in J,\,
  j\neq k}v_j)\pa +\sum_{j\in J,\,
  j\neq k}\pa v_j\pa\ge \pa u\pa.
\end{eqnarray*}
Therefore $\psi_k$ is isometric.
\end{proof}

Using Lemma~\ref{amal:lemma} one gets immediately the triangle
inequality:
\begin{lemma} \label{triangle dist_oq:lemma}
For any compact quantum metric spaces $(A, L_A)$, $(B,L_B)$, and
$(C, L_C)$ we have
\begin{eqnarray*}
\dist_{\oq}(A, C)\le \dist_{\oq}(A, B)+\dist_{\oq}(B, C).
\end{eqnarray*}
For $R\ge 0$ we also have
\begin{eqnarray*}
\dist^R_{\oq}(A, C)\le \dist^R_{\oq}(A, B)+\dist^R_{\oq}(B, C).
\end{eqnarray*}
\end{lemma}

Next we compare $\dist_{\oq}$ (and $\dist^R_{\oq}$) with $\dist_{\q}$.
For this purpose we express first $\dist_{\q}$ in a form similar to that
of $\dist_{\oq}$. For any compact quantum metric space $(A, L_A)$
denote by $\cE(A)$ the unit ball of $A$ under $L_A$.

\begin{proposition} \label{dist_q:prop}
For any compact quantum metric spaces $(A, L_A)$ and $(B,L_B)$ we
have
\begin{eqnarray*}
\dist_{\q}(A, B)
 =\inf \{\dist^V_{\rH}(\cE(A), \cE(B))\},
\end{eqnarray*}
where the infimum is taken over all order-unit spaces $V$ containing both $A$ and $B$ as
order-unit subspaces. This identity also holds if the infimum is taken over all
normed spaces $V$ containing both $A$ and $B$ such that $e_A=e_B$.
\end{proposition}
\begin{proof} Denote the right hand side of the above identity by $\dist'_{\q}(A, B)$.
Also denote by $\dist''_{\q}(A, B)$ the corresponding term for the infimum being taken over all
normed spaces $V$ containing both $A$ and $B$ such that $e_A=e_B$. Clearly $\dist'_{\q}(A, B)\ge \dist''_{\q}(A, B)$.

Let $L$ be an admissible Lip-norm on $A\oplus B$, and
  set $d=\dist^{\rho_L}_{\rH}(S(A), S(B))$.
Denote by $Z$  the subset of $S(A)\times S(B)$ consisting of pairs
$(p, q)$ with $\rho_L(p, q)\le d$. Since $S(A)$ and $S(B)$ are
compact, the projections $Z\to S(A)$ and $Z\to S(B)$ are
surjective. Think of $A$ and $B$ as subspaces of $C(S(A))$ and
$C(S(B))$ respectively. Then the induced $\Re$-linear maps $A\to C(Z)$ and
$B\to C(Z)$ are unital isometric embeddings. Notice
that for any $a\in A$ and $b\in B$ we have
\begin{eqnarray*}
\pa a-b\pa=\sup\{|p(a)-q(b)|:(p, q)\in Z\}\le L(a, b)d.
\end{eqnarray*}
Let $a\in \cE(A)$.  For any
  $\varepsilon>0$ pick $b\in B$ with $L(a, b)<1+\varepsilon$. Then $\pa
a-b\pa\le L(a, b)d\le (1+\varepsilon)d$, and hence
\begin{eqnarray*}
\pa b\pa \le \pa b-a\pa+\pa a\pa \le
(1+\varepsilon)d+\pa a\pa.
\end{eqnarray*}
Also $L_B(b)\le L(a, b)<1+\varepsilon$. Let
$b'=b/(1+\varepsilon)$.
Then $b'\in \cE(B)$, and
\begin{eqnarray*}
\pa a -b'\pa\le \pa a -b\pa+\pa
b-b'\pa\le (1+\varepsilon)d+\frac{\varepsilon}{1+\varepsilon}\pa b\pa \le (1+2\varepsilon)d+\frac{\varepsilon}{1+\varepsilon}\pa a\pa.
\end{eqnarray*}
Similarly, for any $b\in \cE(B)$ and $\varepsilon>0$
we can find $a'\in \cE(A)$ such that $\pa b-a'\pa \le (1+2\varepsilon)d+\frac{\varepsilon}{1+\varepsilon}\pa b\pa$.
Letting $\varepsilon\to 0$ we get
$\dist'_{\q}(A, B)\le d$. Consequently, $\dist'_{\q}(A, B)\le
\dist_{\q}(A, B)$.

Let $V$ be a normed space $V$ containing both $A$ and $B$ such that $e_A=e_B$,
and set $d=\dist^V_{\rH}(\cE(A), \cE(B))$. Let $\varepsilon>0$ be given.
Define a seminorm $L$ on $A\oplus B$ via $L(a, b)=\max(L_A(a), L_B(b), \pa a-b\pa/(d+\varepsilon))$.
It follows easily from Proposition~\ref{criterion of Lip:prop}
that $L$ is an admissible Lip-norm on $A\oplus B$.
For any $p\in S(A)$, by the Hahn-Banach theorem extend $p$ to a linear functional $\varphi$
on $V$ with $\pa \varphi\pa =1$ and set $q$ to be the restriction of $\varphi$ on $B$.
Since $e_A=e_B$ we have $q(e_B)=1$ and hence $q\in S(B)$. For any $(a, b)\in \cE(A\oplus B)$
we have $|p(a)-q(b)|=|\varphi(a-b)|\le \pa a-b\pa\le d+\varepsilon$. Therefore
$\rho_L(p, q)\le d+\varepsilon$.  Similarly, for any $q'\in S(B)$
we can find $p'\in S(A)$ with $\rho_L(p', q')\le d+\varepsilon$. Thus
$\dist_{\q}(A, B)\le \dist^{\rho_L}_{\rH}(S(A), S(B))\le d+\varepsilon$.
Letting $\varepsilon\to 0$ we get
$\dist_{\q}(A, B)\le d$.
Consequently, $\dist_{\q}(A, B)\le
\dist''_{\q}(A, B)$. This finishes the proof of Proposition~\ref{dist_q:prop}.
\end{proof}

We remark that though $\dist_{\q}$ has a form similar to those of
$\dist_{\oq}$ and $\dist^R_{\oq}$, to prove the criteria
Theorems~\ref{criterion of conv:thm} and \ref{criterion 2 of
conv:thm} we have to use $\dist_{\oq}$ and $\dist^R_{\oq}$ in an
essential way.

\begin{proposition} \label{dist_oq<dist_q:prop}
For any compact quantum metric spaces $(A, L_A)$ and $(B,L_B)$ we
have
\begin{eqnarray}
|r_A-r_B|\le \dist_{\GH}(\cD(A),\cD(B))
&\le &\dist_{\oq}(A, B)\le r_A+r_B, \label{dist_oq 1:eq} \\
|\dist_{\oq}(A, B)-\dist^{r_B}_{\oq}(A, B)|&\le &|r_A-r_B|, \label{dist_oq 2:eq} \\
\dist_{\oq}(A, B)&\le &3\dist_{\q}(A, B).\label{dist_oq 3:eq}
\end{eqnarray}
For $R\ge 0$ we also have
\begin{eqnarray} \label{dist_oq 4:eq}
\dist^R_{\oq}(A, B)\le 2\dist_{\q}(A, B).
\end{eqnarray}
\end{proposition}
\begin{proof} For any compact metric spaces $X$ and $Y$, one has
$|r_X-r_Y|\le \dist_{\GH}(X, Y)$ \cite[Exercise 7.3.14]{BBI01}. Thus
(\ref{dist_oq 1:eq}) is trivial once we notice that $\cD(A)$ has radius $r_A$. To show (\ref{dist_oq 2:eq}) it suffices
to show that $\dist^A_{\rH}(\cD(A), \cD_{r_B}(A))\le |r_A- r_B|$. In fact we have:
\begin{lemma} \label{R, r:lemma}
For any compact quantum metric space $(A, L_A)$ and any $R> r\ge
0$ we have
\begin{eqnarray*}
\dist^A_{\rH}(\cD_{ R}(A), \cD_{
  r}(A))\le R-r.
\end{eqnarray*}
\end{lemma}
\begin{proof}
Notice that $\cD_{ r}(A)$ is a subset of $\cD_{
R}(A)$. For each $a\in \cD_{ R}(A)$ let $a'=\frac{r}{R}a$.
Then $a'\in \cD_{ r}(A)$ and
\begin{eqnarray*}
\pa a-a'\pa = \frac{R-r}{R}\pa a\pa  \le
 R-r.
\end{eqnarray*}
Hence $\dist^A_{\rH}(\cD_{ R}(A), \cD_{ r}(A))\le R-r$.
\end{proof}
Back to the proof of Proposition~\ref{dist_oq<dist_q:prop}. The
inequality (\ref{dist_oq 3:eq}) follows from  (\ref{dist_oq
2:eq}), (\ref{dist_oq 4:eq}), and the fact that $|r_A-r_B|\le
\dist_{\GH}(S(A), S(B))\le \dist_{\q}(A, B)$.
So we are left to prove (\ref{dist_oq 4:eq}).
Let $V$ be a normed space containing both $A$ and $B$ such that
$e_A=e_B$, and set $d=\dist^V_{\rH}(\cE(A), \cE(B))$.
For any $a\in \cD_R(A)$ and $\varepsilon>0$ pick $b\in \cE(B)$
such that $\pa a-b\pa \le d+\varepsilon$. Then $\pa b\pa \le \pa b-a\pa +\pa a\pa
\le d+\varepsilon+R$. By Lemma~\ref{R, r:lemma} we can find
$b'\in \cD_R(B)$ with $\pa b-b'\pa \le d+\varepsilon$. Then
$\pa a-b'\pa \le 2(d+\varepsilon)$. Similarly, for any $b\in \cD_R(B)$
we can find $a'\in \cD_R(A)$ with $\pa a'-b\pa \le 2(d+\varepsilon)$.
It follows that $\dist^R_{\oq}(A, B)\le 2d$. Then (\ref{dist_oq 4:eq})
follows from Proposition~\ref{dist_q:prop}.
\end{proof}

\begin{proposition} \label{dist_q<dist_oq:prop}
Let $(A, L_A)$ and $(B, L_B)$ be compact quantum metric spaces
with $R\ge r_A, r_B$. Then we have
\begin{eqnarray}
\dist_{\q}(A, B)\le \frac{5}{2}\dist^R_{\oq}(A, B), \label{dist_oq 5:eq}\\
\dist_{\q}(A, B)\le 5\dist_{\oq}(A, B). \label{dist_oq 6:eq}
\end{eqnarray}
\end{proposition}
\begin{proof}
Note that (\ref{dist_oq 6:eq}) follows immediately from (\ref{dist_oq
  5:eq}), (\ref{dist_oq 2:eq}), and (\ref{dist_oq 1:eq}).
We prove (\ref{dist_oq 5:eq}). We may assume that both $(A, L_A)$
and $(B, L_B)$ are closed. The case $R=0$ is trivial, so we assume
that $R>0$. Let $V$ be a normed space containing $A$ and $B$, and
let $d=\max(\dist^{V}_{\rH}(\cD_{ R}(A), \cD_{
R}(B)),\, \pa Re_A-Re_B\pa)$. If $d=0$ then it is easy to see from
Lemma~\ref{ball to Lip:lemma} that $(A, L_A)$ and $(B, L_B)$ are
isometric. So we assume that $d>0$. Rieffel used \emph{bridges} in
\cite{Rieffel00} to get upper bounds for $\dist_{\q}(A, B)$. Recall
that a bridge between $(A, L_A)$ and $(B, L_B)$ \cite[Definition
5.1]{Rieffel00} is a seminorm, N, on $A\oplus B$ such that $N$ is
continuous for the order-unit norm on $A\oplus B$, $N(e_A, e_B)=0$
but $N(e_A, 0)\neq 0$, and for any $a\in A$ and $\delta>0$ there
is a $b\in B$ such that $\max(L_B(b), N(a,b))\le L_A(a)+\delta$,
and similarly for $A$ and $B$ interchanged. The importance of
bridges is that the seminorm $L$ on $A\oplus B$ defined by $L(a,
b)=\max(L_A(a), L_B(b), N(a,b))$ is an admissible Lip-norm
\cite[Theorem 5.2]{Rieffel00}. In our situation one natural choice
of $N$ is (the seminorm induced from the quotient map $A\oplus
B\to (A\oplus B)/\Re (e_A, e_B)$ and) the quotient norm on
$(A\oplus B)/\Re (e_A, e_B)$ induced by the norm $\pa (a,
b)\pa_*=\max(\pa a\pa,\, \pa b\pa,\, \pa a-b\pa)$. Let $a\in A$
with $L_A(a)=1$. We can write $a$ as
$a'+\lambda e_A$ with $a'\in \cD_{ R}(A)$ and $\lambda\in
\Re$ by Lemma~\ref{ball to Lip:lemma}. Since $\cD_{ R}(B)$
is compact we can find $b'\in \cD_{ R}(B)$ with $\pa a'
-b'\pa\le d$. If we let $b=b'+\lambda e_B$, then we have $N(a,
b)=N(a', b')$, and we just need $N(a', b')\le 1$. So we need to
replace the norm $\pa \cdot \pa_*$ by $\pa (a, b)\pa_1=\max(\pa
a\pa/R, \, \pa b \pa/R, \, \pa a -b\pa/d)$. Then define $N$ as
$N(a, b)=\inf\{\pa (a, b)+\lambda (e_A, e_B)\pa_1:\lambda \in
\Re\}$. The
  above discussion shows that $N$ is a bridge. Then we have the
  admissible Lip-norm $L$ associated to $N$.

Now let $p\in S(A)$. We need to find $q\in S(B)$ such that
$\rho_L(p, q)\le\frac{5}{2}d$. Let $(a, b)\in A\oplus B$ with
$L(a, b)\le 1$. Adding a scalar multiple of $(e_A, e_B)$, we may
assume that $L_A(a),\, L_B(b),\, \pa (a, b)\pa_1\le 1$. Then $\pa
a-b\pa\le d$, $ \pa a\pa\le R$, and $\pa b\pa\le R$. Hence $a\in
\cD_{ R}(A)$ and $b\in \cD_{ R}(A)$. So we are
looking for $q\in S(B)$ such that $|p(a)-q(b)|\le \frac{5}{2}d$
for all $a\in \cD_{ R}(A), b\in \cD_{ R}(A)$ with
$\pa a-b\pa\le d$. Denote the set of such pairs $(a, b)$ by $X$.
By the Hahn-Banach theorem we can extend $p\in A'$ to a $P\in V'$
with $\pa P\pa_{V'}=\pa p\pa_{A'}=1$. Let $g=P|_B$. Then
$|p(a)-g(b)|=|P(a-b)|\le\pa a -b\pa\le d$ for all $(a, b)\in X$,
and $\pa g\pa_{B'}\le \pa P\pa_{V'}=1$. Also
$|1-g(e_B)|=|P(e_A-e_B)|\le \pa e_A-e_B\pa\le d/R$. Now we need:
\begin{lemma} \label{decom:lemma}
Let $g\in B'$ and $\delta\ge 0$ with $1\ge \pa g\pa\ge g(e_B)\ge 1-\delta>0$. Then
there is a $q\in S(B)$ such that $\pa q-g\pa \le \frac{3}{2}\delta$.
\end{lemma}
\begin{proof}
We use the idea in Lemma 2.1 of \cite{Rieffel99b}. Think of $B$ as
a subspace of $C_{\Re}(S(B))$, the space of $\Re$-valued
continuous functions on $S(B)$. Then by the Hahn-Banach theorem
$g$ extends to a continuous linear functional on $C_{\Re}(S(B))$
with the same norm. Using the Jordan decomposition we can write
$g$ as $\mu-\nu$ with $\pa g\pa=\pa \mu\pa+\pa \nu\pa$, where
$\mu$ and $\nu$ are disjoint non-negative measures on $S(B)$. Then
$1\ge \pa \mu\pa+\pa \nu\pa$ and $\pa \mu\pa -\pa
\nu\pa=\mu(e_B)-\nu(e_B)=g(e_B)\ge 1-\delta$. Consequently $\pa
\mu\pa\ge 1-\delta>0$ and $\pa \nu \pa \le \delta/2$. Note that
$\pa \mu\pa=\mu(e_{B})=\pa \mu\pa_{B'}$. Let $q=\mu/\pa \mu\pa$.
Then $q\in S(B)$ and $\pa g-q\pa\le \pa \nu\pa+\pa \mu-q\pa\le
\delta/2+\pa q(1-\pa \mu\pa)\pa\le \frac{3}{2}\delta$.
\end{proof}
Back to the proof of Proposition~\ref{dist_q<dist_oq:prop}.
Pick $q$ for $g$ and $\delta=d/R$ as in Lemma~\ref{decom:lemma}.
Then $|p(a)-q(b)|\le |p(a)-g(b)|+|q(b)-g(b)|\le
d+\frac{3}{2}(d/R)R=\frac{5}{2}d$ for all $(a, b)\in X$.
Consequently $\dist_{\q}(A, B)\le \dist^{\rho_L}_{\rH}(S(A), S(B))\le \frac{5}{2}d$.
Letting $V$ run over all normed spaces containing $A$ and $B$, we get (\ref{dist_oq 5:eq}).
\end{proof}

Now Theorem~\ref{dist_q=dist_oq:thm} follows from Lemma~\ref{triangle dist_oq:lemma}
and  Propositions~\ref{dist_oq<dist_q:prop}
and \ref{dist_q<dist_oq:prop}.
We do not know whether the constants in Theorem~\ref{dist_q=dist_oq:thm} are the best ones or
not.

\begin{remark} \label{add o-u:remark}
Notice that the terms $\pa h_A(r_Ae_A)-h_B(r_Be_B)\pa$ and $\pa
h_A(Re_A)-h_B(Re_B)\pa$ in Definition~\ref{dist_oq:def} are used
only in the proof of Proposition~\ref{dist_q<dist_oq:prop} (and
hence Theorem~\ref{dist_q=dist_oq:thm}). Denote by $\dist'_{\oq}(A,
B)$ the distance omitting the term $\pa
h_A(r_Ae_A)-h_B(r_Be_B)\pa$. If one can show that any compact
quantum metric spaces $(A, L_A)$ and $(B, L_B)$ with
$\dist'_{\oq}(A, B)=0$ are isometric, \ie $\dist'_{\oq}$ is a metric
on $\CQM$, then it is not hard to use (\ref{dist_oq
1:eq}), (\ref{dist_oq 3:eq}), Lemma~\ref{com equi of dual:lemma},
and Theorems~\ref{GH:thm} and \ref{QGH:thm} to show that $\dist_{\q}$
and $\dist'_{\oq}$ define the same topology on $\CQM$. When
$(A, L_A)$ and $(B, L_B)$ are closed, it is easy to see that
 $\dist'_{\oq}(A, B)=0$ if and only there is an affine isometry from
 $\cD(A)$ onto $\cD(B)$ (which has
 to map $0_A$ to $0_B$). Clearly
 such isometry extends to a linear isometry from $A$ onto $B$.
Thus the question is:
\end{remark}

\begin{question} \label{dist'_oq:ques}
Let $(A, L_A)$ and $(B, L_B)$ be compact quantum metric spaces. If
there is a linear isometry (for the norms) $\varphi$ from $A$ onto
$B$ mapping $\cD(A)$ onto $\cD(B)$, then are $(A, L_A)$ and
$(B, L_B)$ isometric as quantum metric spaces?
\end{question}

Notice that if $\varphi(e_A)=e_B$ then $\varphi$ is an isometry as
quantum metric spaces. A related question is:

\begin{question} \label{o-u:ques}
Let $(A, e_A)$ and $(B, e_B)$ be order-unit spaces. If
they are isometric as normed spaces, then must they be isomorphic as
order-unit spaces?
\end{question}

\section{Quantum compactness theorem}
\label{compact:sec}

In this section we prove Theorem~\ref{QGH 2:thm}, which describes
Rieffel's quantum compactness  theorem in terms of the balls
$\cD(A)$.

The main fact we need is
Corollary~\ref{com equi of dual:coro}, which can be proved directly.
Since  Proposition~\ref{univ embedQ:prop} will be useful at other
places,  we include Proposition~\ref{univ embedQ:prop} here, and
deduce Corollary~\ref{com equi of dual:coro} from it. We shall
need the following well-known fact several times. We omit the
proof.

\begin{lemma} \label{cov of sub:lemma}
Let $(X, \rho)$ be a metric space and $Y$ a subset of $X$. Then for
any $\varepsilon>0$ we have $\Cov_{\rho}(Y, 2\varepsilon)\le \Cov_{\rho}(X,
\varepsilon)$, where $\Cov_{\rho}(X, \varepsilon)$ is the smallest number
of open balls of radius $\varepsilon$ whose union covers $X$.
\end{lemma}

Proposition~\ref{univ embedQ:prop} is the dual version of
the fact that if a subset $\mathcal{S}\subseteq \CM$ satisfies the two
conditions in Theorem~\ref{GH:thm}, then there is a compact metric
space $(Z, \rho)$ such that each $X\in \mathcal{S}$ can be
isometrically embedded into $Z$ \cite[page 65]{Gromov81}.
The proof here is a modification of that for
this fact given in \cite{Gromov81}.

\begin{proposition} \label{univ embedQ:prop}
Let $R\ge 0$. For any compact metric space $(X, \rho)$ let
$C(X)_R:=\{f\in C(X):L_{\rho}(X)\le 1, \pa f\pa \le R\}$, equipped
with the metric induced from the supremum norm in the algebra $C(X)$
of $\Ce$-valued continuous functions on $X$, where
$L_{\rho}$ is the Lipschitz seminorm as defined by (\ref{dist to
Lip:eq}). If a subset $\mathcal{S}\subseteq \CM$
satisfies the condition (2) in Theorem~\ref{GH:thm}, then there
exist a complex Banach space $V$ and a compact convex subset
$Z\subseteq V$ such that for every $(X, \rho)\in \mathcal{S}$
there is a linear isometric embedding $h_X: C(X)\hookrightarrow V$
with
 $h_X(C(X)_R)\subseteq Z$.
\end{proposition}
\begin{proof}
For any $(X, \rho)\in \mathcal{S}$ if we pick a dense sequence in
$X$, then the linear map $C(X)\to \ell^{\infty}$ given by the
evaluations at these points is an isometric embedding. What we
shall do is to choose this dense sequence carefully such that the
image of $C(X)_R$ is contained in some compact $Z\subseteq
\ell^{\infty}$ which does not depend on $(X, \rho)$.

Let $\varepsilon_j=2^{-j}$ for all $j\in \Ne$.  Also let
$K_1=\sup\{\Cov(X, \varepsilon_1):(X, \rho)\in \mathcal{S}\}$ and
$K_j=\sup\{\Cov(X, \frac{\varepsilon_j}{2}): (X, \rho)\in
\mathcal{S}\}$ for all $j>1$.  Denote by $D_j$ the set of all
finite sequences of the form $(n_1, n_2, {\cdots}, n_j), 1\le
n_1\le K_1, 1\le n_2\le K_2, {\cdots}, 1 \le n_j\le K_j$, and
denote by $p_j:D_{j+1}\rightarrow D_j$ the natural projection.

We claim that for each $(X, \rho)\in \mathcal{S}$ there are maps
$I^j_X:D_j\rightarrow X$ with the following properties:

\noindent (a) the image of $I^j_X$ forms an $\varepsilon_j$-net in
$X$, \ie the open $\varepsilon_j$-balls centered at the points of
this image cover $X$;

\noindent (b) for each $\omega \in D_{j+1}, j=1, 2, {\cdots}$, the
point $I^{j+1}_X(\omega)$ is contained in the open
$\varepsilon_j$-ball centered at $I^j_X(p_j(\omega))$.

These maps are constructed as follows. Notice that $K_1\ge \Cov(X,
\varepsilon_1)$. So we can cover $X$ by $K_1$ open balls of radius
$\varepsilon_1$, and we take any bijective map from $D_1$ onto the
set of centers of these balls. This is our map $I^1_X$. For any
$\varepsilon_1$-ball $\mathcal{B}$, by Lemma~\ref{cov of sub:lemma},
we have $K_2\ge \Cov(X, \frac{\varepsilon_2}{2})\ge \Cov(\mathcal{B},
\varepsilon_2)$. So we can cover each open $\varepsilon_1$-ball by $K_2$
balls of radius $\varepsilon_2$ and map $D_2$ onto the set of centers
of these $\varepsilon_2$-balls so that $(n_1, n_2)$ goes to the
center of a ball which we used to cover the $\varepsilon_1$-ball with
center at $I^1_X((n_1))$. This is our $I^2_X$. Then we cover each
$\varepsilon_2$-ball by $K_3$ open balls of radius $\varepsilon_3$ and
map $D_3$ onto the set of centers of these $\varepsilon_3$-balls so
that $(n_1, n_2, n_3)\in D_3$ goes to the center of a ball which
was used in covering the $\varepsilon_2$-ball with center at
$I^2_X((n_1, n_2))$, and so on.

Denote by $D$ the union $\cup^{\infty}_{j=1}D_j$, and let $V$ be
the space of all bounded $\Ce$-valued functions on $D$. Then $V$
is a Banach space under the supremum norm $\pa \cdot \pa$. Denote
by $Z\subseteq V$ the set which consists of all functions $f$
satisfying the following inequalities:

if  $\omega \in D_1\subseteq D$, then $|f(\omega)|\le R$,

if $\omega \in D_j$ and $j>1$, then $|f(\omega)-f(p_{i-1}(\omega))|\le
\varepsilon_{j-1}$.

Clearly $Z$ is a closed convex subset of $V$.
 We show
that $Z$ is totally bounded. For any $\varepsilon>0$, pick $k$ such
that $\varepsilon_k<\varepsilon$. Let $P_k$ be the map restricting
functions on $D$ to $\cup^k_{j=1}D_j$. From the inequalities above
we see that $|f(\omega)|\le
R
+\sum^{i-1}_{j=1}\varepsilon_j$ for each $f\in Z$ and $\omega \in
D_i$. So $P_k(Z)$ is contained in $F_k:=\{g\in C(\cup^k_{j=1}D_j):
\pa g\pa\le R+\sum^{k-1}_{j=1}\varepsilon_j\}$. Hence $P_k(Z)$ is
totally bounded. Pick $f_1, {\cdots}, f_m$ in $Z$ such that the
open $\varepsilon$-balls around $P_k(f_1), {\cdots}, P_k(f_m)$ cover
$P_k(Z)$. Then for any $f\in Z$, there is some $1\le l\le m$ so
that $\pa P_k(f)-P_k(f_l)\pa<\varepsilon$. This means that
$|f(\omega)-f_l(\omega)|<\varepsilon$ for all $\omega \in
\cup^{k}_{j=1}D_j$. In particular,
$|f(\omega)-f_l(\omega)|<\varepsilon$ for all $\omega \in D_k$. From
the second inequality above we see that $|f(\omega)-f_l(\omega)|<
\varepsilon+\sum^{\infty}_{j=k}\varepsilon_j=\varepsilon+2\varepsilon_k<3\varepsilon$
for all $\omega \in D\setminus \cup^{\infty}_{j=k+1}D_j$. So $\pa
f-f_l\pa<3\varepsilon$. Therefore $f_1, {\cdots}, f_m$ is a
$3\varepsilon$-net of $Z$, and hence $Z$ is totally bounded. So $Z$
is compact.

Denote by $I_X:D\rightarrow X$ the map corresponding to all
$I^j_X, j=1,2, \cdots$. Then we can define $h_X:C(X)\rightarrow V$
as the pull back of $I_X$:
\begin{eqnarray*}
(h_X(f))(\omega)=f(I_X(\omega)), & f \in C(X), & \omega \in D.
\end{eqnarray*}
Clearly $h_X$ is linear. The property (a) implies that $I_X(D)$ is dense in
$X$. Thus the map $h_X$ is isometric. For each $f \in C(X)_R$
and $\omega \in D_1$ ,
we have $|(h_X(f))(\omega)|\le \pa h_X(f)\pa =\pa f\pa\le R$.
If $\omega \in D_j$, and $j>1$, then by property (b), we have
\begin{eqnarray*}
|(h_X(f))(\omega)-(h_X(f))(p_{j-1}(\omega))|
&=&|f(I_X(\omega))-f(I_X(p_{j-1}(\omega)))|
\\ &\le & L_{\rho}(f)\rho(I_X(\omega),I_X(p_{j-1}(\omega)))\le \varepsilon_{j-1}.
\end{eqnarray*}
So $h_X(f)\in Z$.
Therefore $h_X(C(X)_R)$ is contained in $Z$.
\end{proof}

\begin{corollary} \label{com equi of dual:coro}
Let the notation and hypothesis be as in Proposition~\ref{univ embedQ:prop}.
Then the set $\{C(X)_R:(X, \rho)\in
\mathcal{S}\}$ satisfies the condition
(2) in Theorem~\ref{GH:thm}.
\end{corollary}
\begin{proof} This is a direct consequence of Proposition~\ref{univ
    embedQ:prop} and Lemma~\ref{cov of sub:lemma}.
\end{proof}

\begin{lemma}
 \label{com equi of dual:lemma}
Let $\mathcal{S}$ be a subset of $\CQM$. Pick a closed
representative for each element in $\mathcal{S}$. Then the set
$\{S(A): (A, L)\in \mathcal{S}\}$ satisfies the conditions (1) and
(2) in Theorem~\ref{GH:thm} if and only if the set
$\{\cD(A):(A, L)\in \mathcal{S}\}$ does. If $R\ge
\sup\{r_A: (A, L)\in \mathcal{S}\}$, then the set $\{S(A): (A,
L)\in \mathcal{S}\}$ satisfies the condition (2) in
Theorem~\ref{GH:thm} if and only if the set $\{\cD_{
R}(A):(A, L)\in \mathcal{S}\}$ does.
\end{lemma}
\begin{proof}
We prove the first equivalence. The proof for the second one is
similar. Notice that the radius of $\cD(A)$ is
  exactly $r_A$. Thus $\{S(A): (A, L)\in \mathcal{S}\}$
satisfies the condition (1) in Theorem~\ref{GH:thm} if and only if
$\{\cD(A):(A, L)\in \mathcal{S}\}$ does.

Assume that the condition (1) is satisfied now. Let $R> \sup
\{r_A:(A, L)\in \mathcal{S}\}$. Notice that the natural inclusion
$\cD(A)\to C(S(A))$ is isometric, and has image in
$C(S(A))_R$, where $C(S(A))_R$ is defined as in
Proposition~\ref{univ embedQ:prop}. Then the ``only if'' part
follows from Corollary~\ref{com equi of dual:coro} and
Lemma~\ref{cov of sub:lemma}.

Notice that the natural pairing between $A$ and $A'$ gives a map
$\psi:S(A)\to C(\cD(A))$. Clearly $\psi$ maps
$S(A)$ into $C(\cD(A))_R$. From Lemma~\ref{ball to
Lip:lemma} it is easy to see that $\psi$ is isometric. Then the
``if'' part also follows from Corollary~\ref{com equi of
dual:coro} and Lemma~\ref{cov of sub:lemma}.
\end{proof}

Combining Lemma~\ref{com equi of dual:lemma} and Theorem~\ref{QGH:thm} together we get:

\begin{theorem} \label{QGH 2:thm}
 A subset
$\mathcal{S}\subseteq \CQM$ is totally bounded if and
only if

(1') there is a constant $D'$ such that
$\diam(\cD(A))\le D'$ for all $(A, L)\in
\mathcal{S}$;

(2') for any $\varepsilon>0$, there exists a constant
$K'_{\varepsilon}>0$ such that $\Cov(\cD(A),
\varepsilon)\le K'_{\varepsilon}$ for all $(A, L)\in \mathcal{S}$.
\end{theorem}

One can also give a direct proof of Theorem~\ref{QGH 2:thm} (see \cite[Remark 4.10]{Li12}).

In Section~\ref{comp&bm:sec} we shall use Theorem~\ref{QGH 2:thm}
to prove Theorem~\ref{compact to bounded:thm}, which tells us when
a family of compact quantum metric spaces induced from ergodic
actions of a fixed compact group is totally bounded.

\section{Continuous fields of compact quantum metric spaces}
\label{ConField:sec}

In this section we define continuous fields of compact quantum
metric spaces, a framework we shall use in
Section~\ref{CritConv:sec} to discuss the continuity of families
of compact quantum metric spaces with respect to $\dist_{\oq}$. The
main results of this section are Theorem~\ref{lower cont of
radius:thm} and Proposition~\ref{exist section:prop}. We refer the
reader to \cite[Sections 10.1 and 10.2]{Dixmier77} for basic
definitions and facts about continuous fields of Banach spaces.

We first define continuous fields of order-unit spaces. To reflect the
continuity of the order structures, clearly we should require that the
order-unit section is continuous.

\begin{definition}  \label{cont field of order-unit:def}
Let $T$ be a locally compact Hausdorff space.
A \emph{continuous field of order-unit spaces} over $T$
is a continuous field $(\{A_t\},\Gamma)$  of Banach spaces over
$T$, each $A_t$ being a complete order-unit space with its order-unit norm,
and the unit section $e$
given by $e_t=e_{A_t}, t\in T$ being in the space $\Gamma$ of continuous sections.
\end{definition}

\begin{remark} \label{field of order-unit:remark}
Not every continuous field of Banach spaces consisting of
order-unit spaces is continuous as a field of order-unit spaces.
For a trivial example, let $T=[0,1]$, and let $(\{A_t\}, \Gamma)$ be
the trivial field over $T$ with fibres $(A_t, e_{A_t})=(\Re, 1)$.
For each $f\in \Gamma$ define a section $f^*$ as $f^*_t=f_t$ for
$0\le t<1$ and $f^*_1=-f_1$. Then $\Gamma^*=\{f^*:f\in \Gamma\}$
defines a continuous field of Banach spaces over $T$ with the same
fibres, but with the section $t\mapsto 1$ no longer being
continuous.
\end{remark}

Before we define continuous fields of compact quantum metric
spaces, let us take a look at one example:

\begin{example}[Quotient Field of a Compact Quantum Metric Space]
\label{quotient field of CQM:eg} Let $(B, L_B)$ be a closed
compact quantum metric space. Let $T$ be the set of all nonempty
convex closed subsets of $S(B)$. Notice that for any compact
metric space $(X, \rho)$, the space $\SUB(X)$ of closed
nonempty subsets of $X$ is compact equipped with the metric $\dist^X_{\rH}$
\cite[Proposition 7.3.7]{BBI01}. It is easy to see that $(T,
\dist^{S(B)}_{\rH})$ is a closed subspace of $(\SUB(S(B)),$
$\dist^{S(B)}_{\rH})$, and hence is a compact metric space. Now each
$t\in T$ is a closed convex subset of $S(B)$. Let $(A_t, L_t)$ be
the corresponding quotient of $(B, L_B)$ (see the discussion right
after Proposition~\ref{criterion of Lip:prop}). Then
$\overline{A_t}=\Af_{\Re}(t)$. Let $\pi_t:\Af_{\Re}(S(B))\rightarrow
\Af_{\Re}(t)$ be the restriction map. Since each $w\in
\Af_{\Re}(S(B))$ is uniformly continuous over $S(B)$, clearly the
function $t\mapsto \pa \pi_t(w)\pa$ is continuous over $T$. Hence
the sections $\pi(w)=\{\pi_t(w)\}$ for all $w\in  \Af_{\Re}(S(B))$
generate a continuous field of Banach spaces over $T$ with fibres
$\Af_{\Re}(t)=\overline{A_t}$. Notice that the unit section is just
$\pi(e_B)$. So this is a continuous field of order-unit spaces. We
shall call it the \emph{quotient field} of $(B, L_B)$. According
to \cite[Proposition 5.7]{Rieffel00} we have that $\dist_{\q}(A_t,
A_{t_0})\to 0$ as
 $t\to t_0$ for any $t_0\in T$.
\end{example}

Certainly the above example deserves to be called a continuous
field of compact quantum metric spaces. In general, we start with
a continuous field of order-unit spaces $(\{A_t\}, \Gamma)$ over
some $T$ and a Lip-norm $L_t$ on (a dense subspace of) each $A_t$.
These $L_t$'s should satisfy certain continuity conditions for the
field to be called a continuous field of compact quantum metric
spaces. If we look back at the definition of continuous fields of
Banach spaces, we see that the main ingredient is that there are
enough continuous sections. Thus one may want to require that
there are enough sections $f$ with the functions $t\mapsto
L_t(f_t)$ being continuous. But in the above example of the
quotient field, clearly $t\mapsto L_t(\pi_t(w))$ is always
lower semi-continuous, and there are no obvious $w$'s except the
scalars for which the functions $t\mapsto L_t(\pi_t(w))$ are
continuous. Thus this requirement is too strong. Now there are two
weaker ways to explain ``enough continuous sections''. The first
one is that the structure (which is $L_{t_0}$ in our case) at
$A_{t_0}$ should be determined by the sections "continuous at
$t_0$". Let $\Gamma^L_{t_0}$ be the set of sections $f$ in
$\Gamma$ such that $t\mapsto L_t(f_t)$ is continuous at $t_0$.
Then when $L_{t_0}$ is closed, one wants every $a\in A_{t_0}$ to
have a lifting in $\Gamma^L_{t_0}$. When $L_{t_0}$ is not closed
(which happens in a lot of natural examples), recalling how the
closure $L^c_{t_0}$ is defined, one wants $L^c_{t_0}$ to be
determined by these $f_{t_0}$ when $f$ runs over $\Gamma^L_{t_0}$.
The second way to think of ``enough continuous sections'' is that
there should be enough continuous sections to connect the fibres.
Then one wants that for every $a\in A_{t_0}$ there is some $f\in
\Gamma$ such that $f_{t_0}=a$ and $L_t(f_t)\le L_{t_0}(a)$ for all
$t\in T$. This implies that for every $a\in A_{t_0}$ there is some
$f\in \Gamma$ with $f_{t_0}=a$ and $t\mapsto L_t(f_t)$ being
upper semi-continuous at $t_0$. This is weaker than what we get
above in the first way. However, it turns out that this condition
is strong enough for us to prove some properties of continuous
fields of compact quantum metric spaces (see Theorem~\ref{lower
cont of radius:thm} and Proposition~\ref{exist section:prop}),
especially the criteria for continuity under the order-unit
quantum distance (Theorems~\ref{criterion of conv:thm} and
\ref{criterion 2 of conv:thm}).

\begin{definition} \label{cont field of CQM:def}
Let $T$ be a locally compact Hausdorff space, and let $(A_t, L_t)$
be a compact quantum metric space for each $t\in T$, with
completion $\overline{A_t}$. Let $\Gamma$ be the set of continuous
sections of
 a continuous
field of order-unit spaces over $T$ with fibres $\overline{A_t}$.
For each $t_0\in T$ set
$$\Gamma^L_{t_0}=\{f\in \Gamma:
\mbox{ the function } t\mapsto L_t(f_t) \mbox{ is upper
semi-continuous at } t_0\},$$ where we use the convention that
$L_t=+\infty$ on $\overline{A_t}\setminus A_t$. We call $(\{(A_t,
L_t)\}, \Gamma)$ a {\it continuous field of compact quantum metric
spaces} over $T$ if for any $t\in T$ the restriction of $L_t$ to
$\{f_t: f\in \Gamma^L_t\}$ determines the closure of $L_t$, \ie
for any $a\in A_t$ and $\varepsilon>0$ there exists $f\in \Gamma^L_t$
such that $\pa f_t-a\pa<\varepsilon$ and $L_t(f_t)<L_t(a)+\varepsilon$.
Sections in $\Gamma^L_t$ are called \emph{Lipschitz sections} at
$t$. If every $L_t$ is closed, we say that $(\{(A_t, L_t)\},
\Gamma)$ is {\it closed}.
\end{definition}

\begin{remark} \label{appro Lip:remark}
At first sight, for the restriction of $L_t$ to  $\{f_t: f\in
\Gamma^L_t\}$ to determine the closure of $L_t$, we should require
that for any $a\in \overline{A_t}$ and $\varepsilon>0$ there exists
$f\in \Gamma^L_t$ such that $\pa f_t-a\pa<\varepsilon$ and $L_t(f_t)<
\overline{L_t}(a)+\varepsilon$. This seems stronger than  the
condition we put in Definition~\ref{cont field of CQM:def}. In
fact they are equivalent. By the definition of $\overline{L_t}$ we
can find $a'\in A_t$ with $\pa a'-a\pa<\frac{1}{2}\varepsilon$ and
$L_t(a')< \overline{L_t}(a)+\frac{1}{2}\varepsilon$. Assume that the
condition in
 Definition~\ref{cont field of CQM:def} holds. Then there exists $f\in \Gamma^L_t$
such that $\pa f_t-a'\pa<\frac{1}{2}\varepsilon$ and
$L_t(f_t)<L_t(a')+\frac{1}{2}\varepsilon$.
Consequently, $\pa f_t-a\pa<\varepsilon$ and $L_t(f_t)<\overline{L_t}(a)+\varepsilon$.
\end{remark}

As one would expect, the fiber-wise closure of a continuous field
of compact quantum metric spaces is still such a field:

\begin{proposition} \label{closure field:prop}
If $(\{(A_t, L_t)\}, \Gamma)$ is a continuous field of compact
quantum metric spaces over a locally compact Hausdorff space $T$,
then so is $(\{(A^c_t, L^c_t)\}, \Gamma)$. If $(\{(A_t, L_t)\},
\Gamma)$ is closed, then $\overline{A_t}=\{f_t:f\in \Gamma^L_t\}$
for every $t\in T$.
\end{proposition}
\begin{proof} Let $t_0\in T$, and let $a\in A^c_{t_0}$.
We need to
  find $g\in \Gamma$ such that $g_{t_0}=a$ and
$L^c_{t_0}(a)\ge \limsup_{t\to t_0}L^c_t(g_t)$. Take a section
$f\in \Gamma$ with $f_{t_0}=a$. For each $n\in \Ne$ by
Remark~\ref{appro Lip:remark} we can find an $f_n\in
\Gamma^L_{t_0}$ such that $\pa (f_n)_{t_0}-a\pa <\frac{1}{n}$ and
$L_{t_0}((f_n)_{t_0})< L^c_{t_0}(a)+\frac{1}{n}$. There is an open
neighborhood $\mathcal{U}_n$ of $t_0$ with compact closure such
that $\pa (f_n)_t-f_t\pa <\frac{1}{n}$ and $L_t((f_n)_t)<
L^c_{t_0}(a)+\frac{1}{n}$ for all $t\in \mathcal{U}_n$. By
shrinking these neighborhoods we may assume that
$\overline{\mathcal{U}_{n+1}}\subseteq \mathcal{U}_n$ for all $n$.
By Urysohn's lemma \cite[page 115]{Kelley75} we can find a
continuous function $w_n$ on $T$ with $0\le w_n\le 1$,
$w_n|_{T\setminus \mathcal{U}_n}=0$, and
$w_n|_{\mathcal{U}_{n+1}}=1$. Define a section $g$ by
$g_t=(f_1)_t$ for $t\in T\setminus \mathcal{U}_1$, $g_t=
w_n(t)(f_{n+1})_t+(1-w_n(t))(f_n)_t$ for $t\in
\mathcal{U}_n\setminus \mathcal{U}_{n+1}$, and $g_t=f_t$ for $t\in
\cap^{\infty}_{n=1}\mathcal{U}_n$. Clearly $g\in \Gamma$,
$g_{t_0}=a$, and $L^c_{t_0}(a)\ge \limsup_{t\to t_0}L^c_t(g_t)$.
\end{proof}

\begin{example}[Pull back] \label{pull back:eg}
Let $(\{(A_t, L_t)\}, \Gamma)$ be a continuous field of compact
quantum metric spaces over a locally compact Hausdorff space $T$.
Let $T'$ be another locally compact Hausdorff space, and let
$\Phi:T'\rightarrow T$ be a continuous map. Set $(A_{t'},
L_{t'})=(A_{\Phi(t')}, L_{\Phi(t')})$ for each $t'\in T'$. For
each $f\in \Gamma$ , define a section $\Phi^{*}(f)$ over $T'$ by
$(\Phi^{*}(f))_{t'}=f_{\Phi(t')}$. Then the set
$\Phi^{*}(\Gamma)$ of all these sections generates a continuous
field of Banach spaces over $T'$ with fibres
$\overline{A_{t'}}=\overline{A_{\Phi(t)}}$. This is called the
{\it pull back} of the continuous field $(\{\overline{A_t}\},
\Gamma)$. Let $\overline{\Phi^{*}(\Gamma)}$ be the set of all
continuous sections of this field. Notice that the pull back of
the unit section is exactly the unit section on $T'$. So the pull
back is a continuous field of order-unit spaces. Clearly for each
$t'_0\in T'$ and $f\in \Gamma^L_{\Phi(t'_0)}$ the function
$t'\mapsto L_{t'}((\Phi^{*}(f))_{t'})=L_{\Phi(t')}(f_{\Phi(t')})$
is upper semi-continuous at $t'_0$. Hence
$\Phi^{*}(\Gamma^L_{\Phi(t'_0)})\subseteq
\overline{\Phi^{*}(\Gamma)}^L_{t'_0}$. Then it is easy to see
that $(\{(A_{t'}, L_{t'})\}, \overline{\Phi^{*}(\Gamma)})$ is a
continuous field of compact quantum metric spaces over $T'$. We
shall call it the {\it pull back} of  $(\{(A_t, L_t)\}, \Gamma)$.
\end{example}

\begin{example}[Quotient Field continued]
\label{quotient field of CQM 2:eg} Let the notation be as in
Example~\ref{quotient field of CQM:eg}. For each $t_0\in T$ and
$a\in A_{t_0}$, the proof of \cite[Proposition 3.3]{Rieffel00}
shows that we can find $b\in B$ with $\pi_{t_0}(b)=a$ and
$L_B(b)=L_{t_0}(a)$. Then obviously $\pi(b)$ is in
$\Gamma^L_{t_0}$. Therefore $(\{(A_t, L_t)\}, \Gamma)$ is a closed
continuous field of compact quantum metric spaces.
\end{example}

In fact, we can say more about the Lip-norms in the quotient field:

\begin{proposition} \label{Lip of quotient field:prop}
Let $(B, L_B)$ be a closed compact quantum metric space. Let
$(\{(A_t, L_t)\}, \Gamma)$ be the corresponding quotient field of
compact quantum metric spaces. Then for any $t_0\in T$ and $a\in
\overline{A_{t_0}}$ we have that
\begin{eqnarray*}
L_{t_0}(a)&=&\inf\{\limsup_{t\to t_0}L_t(f_t):f\in \Gamma, \,
f_{t_0}=a
\}\\
&=& \inf\{\liminf_{t\to t_0}L_t(f_t):f\in \Gamma, \, f_{t_0}=a \}.
\end{eqnarray*}
\end{proposition}
\begin{proof}
By Proposition~\ref{closure field:prop} we have
 $L_{t_0}(a)\ge \inf\{\limsup_{t\to t_0}L_t(f_t):f\in \Gamma,
 f_{t_0}=a
\}$.
So we just need to show that for any $a\in \overline{A_{t_0}}$ and
$f\in \Gamma$ with $f_{t_0}=a$
we have
$L_{t_0}(a)\le \liminf_{t\to t_0}L_t(f_t)$. If $L_{t_0}(a)=0$, this is
trivial. So we assume that $L_{t_0}(a)>0$. We prove the case
$L_{t_0}(a)<+\infty$. The proof for the case $L_{t_0}(a)=+\infty$ is
similar.

Let $\rho=\rho_{L_B}$.
By \cite[Proposition 3.3]{Rieffel00} we can find $b\in B$ with
$\pi_{t_0}(b)=a$.
For
any $\varepsilon>0$, since $L_{t_0}$ coincides with the Lip-norm induced
by $\rho|_{t_0}$, we can pick distinct points $p_1, p_2$ in $t_0$ with
$\frac{|a(p_1)-a(p_2)|}{\rho(p_1, p_2)}\ge L_{t_0}(a)-\varepsilon$.
For any $\delta>0$ and $t\in T$ with $\dist^{S(B)}_{\rH}(t, t_0)<\delta$,
we can find $q_1, q_2\in t$ with $\rho(p_j, q_j)<\delta$. Since
$b$ is uniformly continuous on $S(B)$, when
$\delta$ is small enough, we have
\begin{eqnarray*}
|\frac{|(\pi(b))_t(q_1)-(\pi(b))_t(q_2)|}{\rho(q_1, q_2)}
-\frac{|a(p_1)-a(p_2)|}{\rho(p_1, p_2)}|
&=&|\frac{|b(q_1)-b(q_2)|}{\rho(q_1, q_2)}
-\frac{|b(p_1)-b(p_2)|}{\rho(p_1, p_2)}|\\
&<&\varepsilon.
\end{eqnarray*}
Then $\frac{|(\pi(b))_t(q_1)-(\pi(b))_t(q_2)|}{\rho(q_1,
  q_2)}\ge L_{t_0}(a)-2\varepsilon$.  Now $f, \pi(b)\in \Gamma$ and
$f_{t_0}=(\pi(b))_{t_0}=a$. This  implies that
$\pa f_t-(\pi(b))_t\pa \to 0$ as $t\to t_0$. For
  $\delta<\frac{1}{3}\rho(p_1, p_2)$, we have $\rho(q_1,
  q_2)\ge \frac{\rho(p_1, p_2)}{3}$.  Hence when
  $t$ is close enough to $t_0$ we have that
\begin{eqnarray*}
|\frac{|f_t(q_1)-f_t(q_2)|}{\rho(q_1, q_2)}
-\frac{|(\pi(b))_t(q_1)-(\pi(b))_t(q_2)|}{\rho(q_1, q_2)}|
<\varepsilon.
\end{eqnarray*}
Therefore
\begin{eqnarray*}
L_t(f_t)\ge \frac{|f_t(q_1)-f_t(q_2)|}{\rho(q_1, q_2)}
\ge \frac{|(\pi(b))_t(q_1)-(\pi(b))_t(q_2)|}{\rho(q_1,
  q_2)}-\varepsilon
\ge  L_{t_0}(a)-3\varepsilon.
\end{eqnarray*}
Thus $L_{t_0}(a)\le \liminf_{t\to t_0}L_t(f_t)$.
\end{proof}

\begin{example} \label{conv field:eg}
Let $(B, L_B)$ be a closed compact quantum metric space, and let
$(\{(A_t, L_t)\},\Gamma)$ be the corresponding quotient field of
compact quantum metric spaces. Let $\{t_n\}_{n\in \Ne}$ be a
sequence of closed convex subsets of $S(B)$ converging to some
closed convex subset $t_0$ under $\dist^{S(B)}_{\rH}$. Set
$T'=\{\frac{1}{n}:n\in \Ne\}\cup\{0\}$ and define
$\Phi:T'\rightarrow T$ by $\Phi(\frac{1}{n})=t_n$ and $
\Phi(0)=t_0$. Then $\Phi$ is continuous. By Example~\ref{pull
back:eg} we have the pull back continuous field $(\{(A_{t'},
L_{t'})\}, \overline{\Phi^{*}(\Gamma)})$ of compact quantum metric
spaces over $T'$. We will call it {\it the continuous field
corresponding to $t_n\to t_0$}.  Clearly the fibres are
$A_{\frac{1}{n}}=A_{t_n}$ and $A_{0}=A_{t_0}$. The field of Banach
spaces $(\{\overline{A_{t'}}\}, \overline{\Phi^{*}(\Gamma)})$ is
generated by the restrictions of functions in $\Af_{\Re}(S(B))$.
\end{example}

In the same way as for Proposition~\ref{Lip of quotient field:prop},
one can show:

\begin{proposition} \label{Lip of limit:prop}
Let $(B, L_B)$ be a closed compact quantum metric space. Let
$(\{(A_t, L_t)\},$ $\Gamma)$ be the corresponding quotient field of
compact quantum metric spaces. Let $\{t_n\}_{n\in \Ne}$ be a
sequence of closed convex subsets of $S(B)$ converging to some
closed convex subset $t_0$ under $\dist^{S(B)}_{\rH}$, and let
$(\{(A_{t'}, L_{t'})\}, \overline{\Phi^{*}(\Gamma)})$ be the pull
back field as in Example~\ref{conv field:eg}. Then for any $a\in
\overline{A_{t_0}}$ we have that
\begin{eqnarray*}
L_{t_0}(a)&=&\inf\{\limsup_{n\to \infty}L_{t_n}(f_{\frac{1}{n}}):
f\in \overline{\Phi^{*}(\Gamma)},\, f_0=a
\}\\
&=&\inf\{\liminf_{n\to \infty}L_{t_n}(f_{\frac{1}{n}}): f\in
\overline{\Phi^{*}(\Gamma)},\, f_0=a \}.
\end{eqnarray*}
\end{proposition}

\begin{theorem} \label{lower cont of radius:thm}
Let $(\{(A_t, L_t)\}, \Gamma)$ be a continuous field of compact
quantum metric spaces over a locally compact Hausdorff space $T$.
Then the radius function $t\mapsto r_{A_t}$ is lower
semi-continuous over $T$.
\end{theorem}

\begin{lemma} \label{quotient field:lemma}
Let $(\{A_t\}, \Gamma)$ be a continuous field of real Banach spaces
over a locally compact Hausdorff space $T$, and let $f$ be a
nowhere-vanishing section in $\Gamma$. Let $\pa \cdot \pa^{\sim}_t$ be the
quotient norm in $A_t/\Re f_t$. Then for any $g\in \Gamma$ the
function $t\mapsto \pa {\tilde g_t}\pa^{\sim}_t$ is continuous over $T$.
\end{lemma}
\begin{proof}
Replacing $f$ by $t\mapsto \frac{f_t}{\pa f_t\pa}$, we may assume
that $\pa f_t\pa=1$ for all $t\in T$.
For every $t\in T$ pick $c_t\in \Re$ with $\pa g_t-c_tf_t\pa
=\pa {\tilde g_t}\pa^{\sim}_t$.
Let $t_0\in T$ and $\varepsilon>0$ be given.
Since $g-c_{t_0}f\in \Gamma$, the function
$t\mapsto \pa g_t-c_{t_0}f_t\pa$ is continuous
over $T$. So there is a neighborhood $\mathcal{U}$ of $t_0$ such that
for any $t\in \mathcal{U}$,
we have $\pa g_t-c_{t_0}f_t\pa<\pa {\tilde
  g_{t_0}}\pa^{\sim}_{t_0}+\varepsilon$,
and hence
$\pa {\tilde g_t}\pa^{\sim}_t<\pa {\tilde g_{t_0}}\pa^{\sim}_{t_0}+\varepsilon$.
This shows that the function $t\mapsto \pa {\tilde g_t}\pa^{\sim}_t$
is upper semi-continuous over $T$.

We proceed to show that the function $t\mapsto \pa {\tilde g_t}\pa^{\sim}_t$
is lower
semi-continuous over $T$. We may assume that the neighborhood $\mathcal{U}$ in
the above is compact. Let $M:=\sup \{\pa g_t\pa:t \in
\mathcal{U}\}<\infty$.
Then for every $t\in \mathcal{U}$ we have that
\begin{eqnarray*}
|c_t|=\pa c_tf_t\pa
\le \pa g_t\pa +
\pa g_t-c_tf_t\pa
\le M+\pa {\tilde g_t} \pa^{\sim}_t<M+\pa {\tilde g_{t_0}}
    \pa^{\sim}_{t_0} +\varepsilon.
\end{eqnarray*}
Let $I=[-(M+\pa {\tilde g_{t_0}}
    \pa^{\sim}_{t_0} +\varepsilon),\, M+\pa {\tilde g_{t_0}}
    \pa^{\sim}_{t_0} +\varepsilon]$.
Clearly the function $(c, t)\mapsto \pa g_t-cf_t\pa$ is continuous
over $I\times \mathcal{U}$. Since $I$ is compact, we can find a neighborhood
    $\mathcal{U}_1\subseteq \mathcal{U}$ of $t_0$ so that $|\pa g_t-cf_t\pa-\pa
    g_{t_0}-cf_{t_0}\pa|<\varepsilon$ for all $(c, t)\in I\times \mathcal{U}_1$.
Then for any $t\in \mathcal{U}_1$ we have that
\begin{eqnarray*}
\pa {\tilde g_{t_0}}\pa^{\sim}_{t_0}\le  \pa g_{t_0}-c_tf_{t_0}\pa
<\pa g_t-c_tf_t\pa+\varepsilon = \pa {\tilde g_t}\pa^{\sim}_t+\varepsilon
\end{eqnarray*}
So the function $t\mapsto \pa {\tilde g_t}\pa^{\sim}_t$ is lower
semi-continuous, and hence continuous, over $T$.
\end{proof}

Taking $f$ in Lemma~\ref{quotient field:lemma} to be the unit section,
we get immediately:

\begin{lemma}  \label{quotient norm:lemma}
Let $(\{A_t\}, \Gamma)$ be a continuous field of order-unit spaces
over a locally compact Hausdorff space $T$.
Let $\pa \cdot \pa^{\sim}_t$ be the
quotient norm in ${\tilde A_t}=A_t/\Re e_{A_t}$. Then for any $f\in
\Gamma$,
the
function $t\mapsto \pa {\tilde f_t}\pa^{\sim}_t$ is continuous over $T$.
\end{lemma}

We are ready to prove Theorem~\ref{lower cont of radius:thm}.

\begin{proof}[Proof of Theorem~\ref{lower cont of radius:thm}]
 By Proposition~\ref{closure field:prop} we may assume
  that $(\{(A_t, L_t)\}, \Gamma)$ is closed.
Let $t_0\in T$ and $\varepsilon>0$ be given.
If $A_{t_0}=\Re e_{A_{t_0}}$, then $r_{A_{t_0}}=0$ and the radius
function is obviously lower semi-continuous at $t_0$.  So we may
assume that $A_{t_0}\neq \Re e_{A_{t_0}}$. Then $r_{A_{t_0}}=\sup
\{\pa {\tilde a}\pa^{\sim}_{t_0}:{\tilde a} \in {\tilde A_{t_0}}
\mbox{ with } {\tilde L_{t_0}}({\tilde a})=1\}$ by
Proposition~\ref{criterion of Lip:prop}. Pick $a\in
A_{t_0}$ with $L_{t_0}(a)={\tilde L_{t_0}}({\tilde a})=1$ and $\pa {\tilde
a}\pa^{\sim}_{t_0}>r_{A_{t_0}}-\varepsilon$.
By Proposition~\ref{closure
  field:prop} we can
find $f\in \Gamma^L_{t_0}$ with $f_{t_0}=a$.
Then the function $t\mapsto L_{t}(f_t)$
is upper semi-continuous at $t_0$.
By Lemma~\ref{quotient norm:lemma} the function $t\mapsto
\pa \widetilde{f_t}\pa^{\sim}_t$ is continuous over $T$.
So there is some neighborhood $\mathcal{U}$ of $t_0$ in
$T$ such that
$\pa \widetilde{f_t}\pa^{\sim}_{t}/L_{t}(f_t)>
\pa \widetilde{f_{t_0}}\pa^{\sim}_{t_0}
-\varepsilon
> r_{A_{t_0}}-2\varepsilon$ for all $t\in \mathcal{U}$.
Then $r_{A_t}> r_{A_{t_0}}-2\varepsilon$ for all $t\in \mathcal{U}$
by Proposition~\ref{criterion of Lip:prop}. So the radius function
is lower semi-continuous at $t_0$.
\end{proof}

Next we show that there are enough Lipschitz sections to connect the
fibres. The next lemma is probably known, but we can not find a
reference, so we include a proof.

\begin{lemma} \label{cont for lower cont:lemma}
Let $T$ be a locally compact Hausdorff space, and let $w$ be a nonnegative
function on $T$. If $w$ is lower semi-continuous at some point $t_0\in T$,
then there is a continuous nonnegative function $w'$ over $T$ with
$w'(t_0)=w(t_0)$
and $w'\le w$ on $T$.
\end{lemma}
\begin{proof} If $w(t_0)=0$, we may take $w'=0$.
So assume $w(t_0)>0$. Replacing $w$ by $\frac{w}{w(t_0)}$, we may
assume
$w(t_0)=1$.

Since $w$ is lower semi-continuous at $t_0$, for each $n\in \Ne$ we
can find an open neighborhood $\mathcal{U}_n$ of $t_0$
such that $w\ge 1-2^{-n}$ on $\mathcal{U}_n$.
By shrinking $\mathcal{U}_n$ we may assume that the closure of
$\mathcal{U}_n$ is compact and contained in $\mathcal{U}_{n-1}$
for each $n\in \Ne$, where $\mathcal{U}_0=T$. By Urysohn's lemma
\cite[page 115]{Kelley75}, we can find a continuous function
$w'_n$ over $T$ with $0\le w'_n\le 2^{-n}$, $w'_n|_{T\setminus
\mathcal{U}_n}=0$, and
$w'_n|_{\overline{\mathcal{U}_{n+1}}}=2^{-n}$. Then
$w'=\sum^{\infty}_{n=1}w'_n$ is continuous over $T$.

If $t\in T\setminus \mathcal{U}_1$, then $w'_j(0)=0$ for all $j$
and hence $w'(0)=0\le w(t)$. If $t\in \mathcal{U}_n\setminus
\mathcal{U}_{n+1}$ for some $n\in \Ne$, then $w'_j(t)=2^{-j}$ for
all $1\le j<n$, $0\le w'_n(t)\le \frac{1}{2^n}$, and $w'_j(t)=0$
for all $j>n$. So $w'(t)\le \sum^{n}_{j=1}2^{-j}=1-2^{-n}\le
w(t)$. If $t\in \cap^{\infty}_{n=1}\mathcal{U}_n$, then $w(t)\ge
1$ according to the construction of $\mathcal{U}_n$. In this case,
$w'_j(t)=2^{-j}$ for all $j$. So $w'(t)=1\le w(t)$. In particular,
we see that $w'(t_0)=1=w(t_0)$. So $w'$ satisfies our requirement.
\end{proof}

\begin{proposition} \label{exist section:prop}
Let $(\{(A_t, L_t)\}, \Gamma)$ be a continuous field of compact
quantum metric spaces over a locally compact Hausdorff space $T$.
Then for any $t_0\in T$ and $a\in \{f_{t_0}:f\in
\Gamma^L_{t_0}\}\cap \cD(A_{t_0})$, there
exists $f\in \Gamma^L_{t_0}$ with $f_t\in \cD(A_t)$ for all $t\in T$
and $f_{t_0}=a$. In particular,
when $(\{(A_t, L_t)\}, \Gamma)$ is closed, such $f$ exists for every
$a\in \cD(A_{t_0})$.
\end{proposition}
\begin{proof} If $a=0$, we can pick $f=0$. So suppose that
$a\neq 0$. Take $g\in \Gamma^L_{t_0}$ with $g_{t_0}=a$. Then
$0<\pa g_{t_0}\pa=\pa a\pa\le r_{A_{t_0}}$. Since $t\mapsto \pa
g_t\pa$ is continuous on $T$, there is some neighborhood
$\mathcal{U}$ of $t_0$ such that $\pa g_t\pa>0$ for all $t\in
\mathcal{U}$. Define a nonnegative function $w$ on $T$ by
$w(t_0)=1$, $w(t)=\frac{r_{A_t}}{\pa g_t\pa}$ for $t\in
\mathcal{U}\setminus \{t_0\}$, and $w(t)=0$ for $t\in T\setminus
\mathcal{U}$. By Theorem~\ref{lower cont of radius:thm} the radius
function $t\mapsto r_{A_t}$ is lower semi-continuous over $T$. Then
it is easy to see that the function $w$ is lower semi-continuous at
$t_0$. According to Lemma~\ref{cont for lower cont:lemma} we can
find a continuous nonnegative function $w'$ on $T$ such that
$w'(t_0)=1$ and $w'(t)\le w(t)$ for all $t\in T$. Then $w'(t)\le
\frac{r_{A_t}}{\pa g_t\pa}$ for all $t\in \mathcal{U}$, and
$w'(t)=0$ for $t\in T\setminus \mathcal{U}$. Set $h_t=w'(t)g_t$.
Then $h\in \Gamma^L_{t_0}$. Also, $h_{t_0}=a$ and $\pa h_t\pa\le
r_{A_t}$ for all $t\in T$.

If $L_{t_0}(a)<1$, then $L_{t_0}(h_{t_0})<1$. Since $h\in
\Gamma^L_{t_0}$, there is an open neighborhood $\mathcal{U}_1$ of
$t_0$ with compact closure such that $L_t(h_t)<1$ for all $t\in
\mathcal{U}_1$. Take a continuous function $w''$ on $T$ with $0\le
w''\le 1$, $w''(t_0)=1$, and $w''|_{T\setminus \mathcal{U}_1}=0$.
Define a section $f$ by $f_t=w''(t)h_t$. Then $f$ is in
$\Gamma^L_{t_0}$, and satisfies $L_t(f_t)\le 1$, $\pa f_t\pa \le
\pa h_t\pa \le r_{A_t}$ for all $t\in T$. So $f_t\in \cD(A_t)$ for all $t\in T$ . Also
$f_{t_0}=w''(t_0)h_{t_0}=a$. Hence $f$ satisfies our requirement.

Now suppose that $L_{t_0}(a)=1$. Then $L_{t_0}(h_{t_0})=1$. Define
a nonnegative function $w_1$ on $T$ as $w_1(t)=\min
(\frac{1}{L_t(h_t)}, 1)$ for all $t\in T$, where
$\frac{1}{L_t(h_t)}=\infty$ if $L_t(h_t)=0$. Then $w_1(t_0)=1$.
Since $h\in \Gamma^L_{t_0}$, it is easy to see that $w_1$ is lower
semi-continuous at $t_0$. According to Lemma~\ref{cont for lower
cont:lemma} we can find a continuous nonnegative function $w'_1$
on $T$ such that $w'_1(t_0)=1$ and $w'_1(t)\le w_1(t)$ for all
$t\in T$. Then $w'_1\le w_1\le 1$ on $T$. Define a section $f$ by
$f_t=w'_1(t)h_t$ for all $t\in T$. Then $f\in \Gamma^L_{t_0}$ and
$f_{t_0}=w'_1(t_0)h_{t_0}=a$. Clearly $L_t(f_t)=w'_1(t)L_t(h_t)\le
w_1(t)L_t(h_t)\le 1$ and $\pa f_t\pa \le \pa h_t\pa \le r_{A_t}$
for all $t\in T$. So $f_t\in \cD(A_t)$ for all
$t\in T$. Hence $f$ satisfies our requirement.

The assertion about closed $(\{(A_t, L_t)\}, \Gamma)$ follows from
Proposition~\ref{closure field:prop}.
\end{proof}

\section{Criteria for metric convergence}
\label{CritConv:sec}

In this section we prove Theorems~\ref{criterion of conv:thm} and \ref{criterion 2 of conv:thm}.

When applying Theorem~\ref{criterion of conv:thm}, usually we need
to show that the radius function $t\mapsto r_{A_t}$ is continuous
at $t_0$. This is often quite difficult. However, sometimes we can
show easily that the radii are bounded (for example, compact
quantum metric spaces induced by ergodic actions of compact groups
on complete order-unit spaces of finite multiplicity, see
Theorem~\ref{finite to Lip:thm}). In these cases, the next
criterion is more useful. This is also the reason we introduced
$\dist^R_{\oq}$.

\begin{theorem}  \label{criterion 2 of conv:thm}
Let $(\{(A_t, L_t)\}, \Gamma)$ be a continuous field of compact
quantum metric spaces over a locally compact Hausdorff space $T$.
Let $R\ge 0$. Let $t_0\in T$, and let $\{f_n\}_{n\in \Ne}$ be a
sequence in $\Gamma$ such that $(f_n)_{t_0}\in \cD_{
R}(A_{t_0})$ for each $n\in \Ne$ and the set $\{(f_n)_{t_0}:n\in
\Ne\}$ is dense in $\cD_{R}(A_{t_0})$. Then the following
are equivalent:

(i) $\dist^R_{\oq}(A_t, A_{t_0})\to 0$ as $t\to t_0$;

(ii) $\dist_{\GH}(\cD_{ R}(A_t), \cD_{
R}(A_{t_0}))\to 0$ as $t\to t_0$;

(iii) for any $\varepsilon>0$, there is an $N$ such that the open
$\varepsilon$-balls in $A_t$ centered at $(f_1)_t, {\cdots}, (f_N)_t$
cover $\cD_{ R}(A_t)$ for all $t$ in some neighborhood
$\mathcal{U}$ of $t_0$.

\end{theorem}

The proof of
Theorem~\ref{criterion 2 of
conv:thm} is
similar to that of Theorem~\ref{criterion of conv:thm}.
So we shall prove
only Theorem~\ref{criterion of conv:thm}. We need some
preparation. The next lemma generalizes Lemma~\ref{amal:lemma} to
deal with ``almost amalgamation'':
\begin{lemma} \label{almost amal:lemma}
Let $A$ and $B$ be normed spaces (over $\Re$ or $\Ce$). Let $X$ be a
linear subspace of $A$, and let $\varepsilon\ge 0$.
Let $\varphi:X\to B$ be a linear map with
$(1-\varepsilon)\pa x\pa \le \pa \varphi(x)\pa\le (1+\varepsilon)\pa x\pa$
for all $x\in X$.
Then there are a normed space $V$ and linear isometric
embeddings $h_A:A\hookrightarrow V$ and $h_B:B\hookrightarrow V$
such that
$\pa h_A(x)-(h_B\circ
\varphi)(x)\pa\le \varepsilon\pa x\pa$ for all $x\in X$.
\end{lemma}
\begin{proof}
We define a seminorm, $\pa \cdot \pa_*$, on $A\oplus B$ by
\begin{eqnarray*}
\pa (a, b)\pa_*:=\inf\{\pa a-x\pa+\pa b+\varphi(x)\pa +\varepsilon
\pa x\pa:x\in X\}.
\end{eqnarray*}
We claim that $\pa \cdot \pa_*$ extends the norm of $A$. Let $a\in
A$. Taking $x=0$ we get $\pa (a, 0)\pa_*\le \pa a\pa$. For any
$x\in X$ we have $\pa a-x\pa+\pa 0+\varphi(x)\pa +\varepsilon \pa
x\pa \ge \pa a-x\pa+(1-\varepsilon)\pa x\pa +\varepsilon \pa x\pa \ge
\pa a\pa$. So $\pa(a, 0)\pa_*\ge \pa a\pa$,  and hence $\pa
(a,0)\pa_*=\pa a\pa$. Similarly, $\pa \cdot \pa_*$ extends the
norm of $B$. For any $x\in X$ we have $\pa (x,
-\varphi(x))\pa_*\le \pa x-x\pa+\pa
-\varphi(x)+\varphi(x)\pa+\varepsilon \pa x\pa=\varepsilon\pa x\pa$. Let
$N$ be the null space of $\pa \cdot \pa_*$, and let $V$ be
$(A\oplus B)/N$. Then $\pa\cdot \pa_*$ induces a norm on $V$, and
the natural maps $A,\, B\rightarrow V$ satisfy the requirement.
\end{proof}

A \emph{subtrivialization} \cite[page 133]{Blanchard98} of a
continuous field of $C^*$-algebras $(\{\mathcal{A}_t\}, \Gamma)$
over a locally compact Hausdorff space $T$ is a faithful
$*$-homomorphism $h_t$ of each  $\mathcal{A}_t$ into a common
$C^*$-algebra $\mathcal{A}$ such that for every $f\in \Gamma$ the
$\mathcal{A}$-valued function $t\mapsto h_t(f_t)$ is continuous over $T$.
Not every continuous field of $C^*$-algebras can be subtrivialized
\cite[Remark 5.1]{KW96}. We can talk about the subtrivialization
of continuous fields of Banach spaces similarly by requiring
$h_t$'s to be linear isometric embeddings into some common Banach
space.
One natural question is:

\begin{question} \label{subtrivial:ques}
Can every continuous field of Banach spaces over a locally compact
Hausdorff space be subtrivialized?
\end{question}

Blanchard and Kirchberg gave affirmative answer for separable continuous fields of
complex Banach spaces over compact metric spaces \cite[Corollary 2.8]{BK}.
However, to prove Theorem~\ref{criterion of conv:thm} we have to deal with
continuous fields of real separable Banach spaces
over general locally compact Hausdorff spaces.
For us the following weaker answer to Question~\ref{subtrivial:ques} is sufficient:

\begin{proposition} \label{subtrivial:prop}
Let $(\{A_t\}, \Gamma)$ be a continuous field of Banach spaces (over
$\Re$ or $\Ce$) over a locally compact Hausdorff space $T$.
Let $t_0\in T$ with $A_{t_0}$ separable. Then there are a normed space
$V$ and linear isometric embeddings $h_t:A_t\hookrightarrow V$ such
that for every $f\in \Gamma$ the $V$-valued map $t\mapsto h_t(f_t)$ is
continuous at $t_0$.
\end{proposition}
\begin{proof}
We prove the case where $A_{t_0}$ is infinite-dimensional. The
case where $A_{t_0}$ is finite-dimensional is similar and easier.
Since $A_{t_0}$ is separable we can find a linearly independent
sequence $x_1, x_2, \cdots$ in $A_{t_0}$ such that the linear span
of $\{x_k\}_{k\in \Ne}$ is dense in $A_{t_0}$. For each $k$ pick a
section $f_k\in \Gamma$ with $(f_k)_{t_0}=x_k$.  Then for each $t$
the map $\varphi_t:x_k\mapsto (f_k)_t$, $k=1,2, \cdots$, extends
uniquely to a linear map from $\spn\{x_k:k\in \Ne\}$ to $A_t$,
which we still denote by $\varphi_t$. Let $X_n=\spn\{x_1, \cdots,
x_n\}$, and let $\varphi_{n,t}$ be the restriction of $\varphi_t$
on $X_n$. Notice that for each $x\in X_n$ the section $t\mapsto
\varphi_{n,t}(x)$ is in $\Gamma$. Then using a standard
compactness argument we can find a neighborhood $\mathcal{U}_n$ of
$t_0$ such that $1-\frac{1}{n}\le \pa \varphi_{n,t}(x)\pa \le
1+\frac{1}{n}$ for all $x$ in the unit sphere of $X_n$ and $t\in
\mathcal{U}_n$. We may assume that $\mathcal{U}_1\supseteq
\mathcal{U}_2\supseteq\cdots$. We shall find a normed space $V_t$
containing both $A_t$ and $A_{t_0}$ for each $t\in \Gamma$ such
that $A_t$ and $A_{t_0}$ are kind of close to each other inside of
$V_t$. If $t\not \in \mathcal{U}_1$ we let $V_t$ simply be
$A_t\oplus A_{t_0}$ equipped with any admissible norm. If  $t\in
\mathcal{U}_n\setminus \mathcal{U}_{n+1}$, then by
Lemma~\ref{almost amal:lemma} we can find a normed space $V_t$
containing both $A_t$ and $A_{t_0}$  such that $\pa
x-\varphi_{n,t}(x)\pa \le \frac{1}{n}\pa x\pa$ for all $x\in X_n$.
If $t\in \cap^{\infty}_{n=1}\mathcal{U}_n$, then $\varphi_{n,t}$
is an isometric embedding for all $n$, and hence $\varphi_t$
extends to a linear isometric embedding from $A_{t_0}$ into $A_t$.
So  for $t\in \cap^{\infty}_{n=1}\mathcal{U}_n$ we can identify
$A_{t_0}$ with $\varphi_t(A_{t_0})$, and let $V_t=A_t$.
 Now by Lemma~\ref{amal:lemma}
we can find a normed space $V$ containing
all these $V_t$'s such that the copies of $A_{t_0}$ are identified.
Let $h_t$ be the composition $A_t\hookrightarrow V_t\hookrightarrow
V$. Then for each $x\in \spn\{x_k:k\in \Ne\}$ clearly the map
$t\mapsto h_t(\varphi_t(x))$ is continuous at $t_0$. Now it is easy to see that
for every section $f\in \Gamma$ the map $t\mapsto h_t(f_t)$ is
continuous  at $t_0$.
\end{proof}

\begin{remark} \label{trivial:remark}
The $C^*$-algebraic analogue of Proposition~\ref{subtrivial:prop}
is not true, \ie for a continuous field $(\{\mathcal{A}_t\},
\Gamma)$ of $C^*$-algebras, in general we can not find a
$C^*$-algebra $\mathcal{B}$ and faithful $*$-homomorphisms
$h_t:\mathcal{A}_t\hookrightarrow \mathcal{B}$ such that for every
$f\in \Gamma$ the map $t\mapsto h_t(f_t)$ is continuous at $t_0$.
The reason is that such $\mathcal{B}$ and $h_t$'s will imply that
for any $C^*$-algebra $\mathcal{C}$ and any $\sum_jf_j\otimes c_j$
in the algebraic tensor product $\Gamma\otimes_{alg} \mathcal{C}$
the function $t\mapsto \pa \sum_j(f_j)_t\otimes
c_j\pa_{\mathcal{A}_t\otimes \mathcal{C}}$ is continuous at
  $t_0$, where $\mathcal{A}_t\otimes \mathcal{C}$ is the minimal
  tensor product. But there are examples \cite[Proposition 4.3]{KW96}
  where $t\mapsto \pa \sum_j(f_j)_t\otimes
c_j\pa_{\mathcal{A}_t\otimes \mathcal{C}}$ is not continuous, even
when $T$ is simply the one-point compactification of $\Ne$. Notice
that in the proof of Proposition~\ref{subtrivial:prop} we used
only Lemmas~\ref{amal:lemma} and \ref{almost amal:lemma}. The
$C^*$-algebraic analogue of Lemma~\ref{amal:lemma} has been proved
by Blackadar \cite[Theorem 3.1]{Blackadar80}. Thus the
$C^*$-algebraic analogue of Lemma~\ref{almost amal:lemma} (with
$\varphi$ still being a linear map, but $h_A$ and $h_B$ being
faithful $*$-homomorphisms) is not true.
\end{remark}

Recall that for a metric space $X$ and $\varepsilon>0$ the \emph{packing
  number} $P(X, \varepsilon)$ is the maximal cardinality of an
  $\varepsilon$-separated (\ie $\rho_X(x, x')>\varepsilon$ if $x\neq x'$) subset in
  $X$. When $X$ is compact, $P(X, \varepsilon)$ is finite. In fact
  clearly $P(X, \varepsilon)\le \Cov(X, \frac{1}{2}\varepsilon)$.

\begin{lemma} \label{GH to ball:lemma}
Let $X$  be a compact metric space, and let $\varepsilon>0$.
For any closed subset $Y$ of
$X$ if $\dist_{\GH}(X, Y)< \frac{1}{4P(X, \varepsilon/2)}\varepsilon$,
then the open $\varepsilon$-balls centered at points of $Y$ cover $X$.
\end{lemma}
\begin{proof} Let $N=P(X, \varepsilon/2)$.
Let $h_X:X\rightarrow Z$ and $h_Y:Y\rightarrow Z$ be isometric embeddings into some metric space $Z$
such that $\frac{1}{4N}\varepsilon>\dist^Z_{\rH}(h_X(X), h_Y(Y))$.
For each $x\in X$ pick $\varphi(x)\in Y$ with
$\frac{1}{4N}\varepsilon>\rho_Z(h_X(x), h_Y(\varphi(x)))$. Let $x\in X$.
Define $x_n$ inductively by $x_0=x$ and $x_{n}=\varphi(x_{n-1})$.
Then for any $m>n\ge 1$ we have that
\begin{eqnarray*}
& &\rho_X(x_n, x_m) \\
&=&\rho_Y(x_n, x_m)\\
&\ge &
 \rho_X(x_{n-1}, x_{m-1})-\rho_Z(h_X(x_{n-1}), h_Y(x_n))-
\rho_Z(h_X(x_{m-1}), h_Y(x_m))\\
&\ge &\rho_X(x_{n-1},
x_{m-1})-\frac{1}{2N}\varepsilon.
\end{eqnarray*}
Consequently $\rho_X(x_n, x_m)\ge \rho_X(x_0,
x_{m-n})-\frac{n}{2N}\varepsilon\ge \rho_X(x,
Y)-\frac{n}{2N}\varepsilon$ for all $m>n\ge 0$. Therefore $x_0, x_1,
{\cdots}, x_N$ are $(\rho_X(x, Y)-\frac{1}{2}\varepsilon)$-separated.
Thus $\rho_X(x, Y)-\frac{1}{2}\varepsilon<\frac{1}{2}\varepsilon$. Then
$\rho_X(x, Y)<\varepsilon$ follows.
\end{proof}

\begin{lemma} \label{packing:lemma}
Let $X$ and $Y$ be compact metric spaces, and let $\varepsilon>0$.
If $\dist_{\GH}(X, Y)\\< \frac{1}{4}\varepsilon$ then $P(X, \varepsilon)\le
P(Y, \frac{1}{2}\varepsilon)$.
\end{lemma}
\begin{proof}
Let $\rho$ be an admissible metric on $X\coprod Y$ with
$\dist^{\rho}_{\rH}(X, Y)<\frac{1}{4}\varepsilon$ (see the discussion
preceding Theorem~\ref{GH:thm}). Let $\{x_1, {\cdots}, x_n\}$ be
an $\varepsilon$-separated set in $X$. For each $k$ pick $y_k\in Y$
such that $\rho(x_k, y_k)<\frac{1}{4}\varepsilon$. Then clearly
$\{y_1, {\cdots}, y_n\}$ is $\frac{1}{2}\varepsilon$-separated in
$Y$. Therefore $n\le P(Y, \frac{1}{2}\varepsilon)$.
\end{proof}

Now we are ready to prove Theorem~\ref{criterion of conv:thm}.

\begin{proof}[Proof of Theorem~\ref{criterion of conv:thm}]
We claim first that (iii) does not depend on the choice of the
sequence $\{f_n\}_{n\in \Ne}$. Suppose that $\{f'_n\}_{n\in \Ne}$
is another sequence in $\Gamma$ satisfying the conditions in the
theorem. If (iii) holds for $\{f_n\}_{n\in \Ne}$, then for any
$\varepsilon>0$, we can find $N$ and a neighborhood $\mathcal{U}$ as
in (iii). Since $\{(f'_n)_{t_0}:n\in \Ne\}$ is dense in
$\cD(A_{t_0})$, there is some $N'\in \Ne$
so that for each $1\le n\le N$, there is some $1\le \sigma(n)\le
N'$ with $\pa
(f_n)_{t_0}-(f'_{\sigma(n)})_{t_0}\pa_{t_0}<\varepsilon$. Then we can
find a neighborhood $\mathcal{U}'\subseteq \mathcal{U}$ of $t_0$
such that $\pa (f_n)_t-(f'_{\sigma(n)})_t \pa_{t_0}<2\varepsilon$ for
all $1\le n\le N$ and all $t\in \mathcal{U}'$. It is clear that
the open $3\varepsilon$-balls in $A_t$ centered at $(f'_1)_t,
{\cdots}, (f'_{N'})_t$ cover $\cD(A_t)$ for all
$t\in \mathcal{U}'$. So (iii) is also satisfied for
$\{f'_n\}_{n\in \Ne}$, and hence it does not depend on the choice
of the sequence $f_n$.

Since $\cD(A_t)$ is dense in $\cD(A^c_t)$ for every $t$, (iii) does not depend on whether we
take  $(\{(A_t, L_t)\}, \Gamma)$ or its closure. Clearly neither
does (i) nor (ii). So we may assume that $(\{(A_t, L_t)\}, \Gamma)$
is closed. Take a dense sequence $\{a_n\}_{n\in \Ne}$ in
$\cD(A_{t_0})$. According to
Proposition~\ref{exist section:prop} we can find $f_n\in \Gamma$
with $(f_n)_{t_0}=a_n$ and $(f_n)_t\in \cD(A_t)$ for all $t\in T$. Then $\{f_n\}_{n\in \Ne}$
satisfies the condition in the theorem. In the rest of the  proof
we will use this sequence $\{f_n\}_{n\in \Ne}$.

Since $S(A_{t_0})$ is a compact metric space, $A_{t_0}\subseteq
C(S(A_{t_0}))$ is separable. So by
Proposition~\ref{subtrivial:prop} we can find a normed space $V$
containing all $V_t$'s such that for every $f\in \Gamma$ the map
$t\mapsto f_t$ from $T$ to $V$ is continuous at $t_0$. For any
$n\in \Ne$ and $t\in T$ let $Y_{n,t}=\{(f_1)_t, {\cdots},
(f_n)_t\}$. Also let $X_t=\cD(A_t)$ for all
$t\in T$.

We show first that (iii)$\Rightarrow$(i). Let $\varepsilon>0$ be
given. Pick $N$ and a neighborhood $\mathcal{U}$ of $t_0$ for
$\varepsilon$ as in (iii). Then $\dist^V_{\rH}(Y_{N, t}, X_t)\le \varepsilon$
throughout $\mathcal{U}$. By shrinking $\mathcal{U}$ we may assume
that $\pa e_{A_t}-e_{t_0}\pa\le \varepsilon$ and $\pa
(f_k)_t-(f_k)_{t_0}\pa \le \varepsilon$ for all $t\in \mathcal{U}$
and $1\le k\le N$. Let $t\in \mathcal{U}$. Then $\dist^V_{\rH}(Y_{N,
t}, Y_{N, t_0})\le \varepsilon$. Hence $\dist^V_{\rH}(X_t, X_{t_0})\le
\dist^V_{\rH}(X_t,Y_{N, t}) +\dist^V_{\rH}(Y_{N, t},Y_{N, t_0})+
\dist^V_{\rH}(Y_{N, t_0}, X_{t_0})\le 3\varepsilon$. Therefore
$\dist_{\oq}(A_t, A_{t_0})\le 3\varepsilon$. This proves
(iii)$\Rightarrow$(i).

(i)$\Rightarrow$(ii) follows from (\ref{dist_oq 1:eq}). So we are
left
   to show that (ii)$\Rightarrow$(iii). Let
   $\delta=\min((12P(X_{t_0},\frac{1}{4}\varepsilon))^{-1}, \frac{1}{8})\varepsilon$.
   Take $N$ so that the set $Y_{N, t_0}$ is
   $\delta$-dense in $X_{t_0}$. Then we have
   $\dist^V_{\rH}(Y_{N, t_0},X_{t_0})\le \delta$.
   Similarly as above there is some neighborhood $\mathcal{U}$ of $t_0$
   such that $\dist^V_{\rH}(Y_{N, t}, Y_{N,t_0})\le \delta$
   for all $t\in \mathcal{U}$.
By shrinking $\mathcal{U}$ we may assume that $\dist_{\GH}(X_t, X_{t_0})< \delta\le \frac{1}{8}\varepsilon$ for
all $t\in \mathcal{U}$. Let $t\in \mathcal{U}$. Then
\begin{eqnarray*}
\dist_{\GH}(Y_{N, t}, X_t) &\le &
  \dist_{\GH}(Y_{N, t}, Y_{N,t_0})+\dist_{\GH}(Y_{N,
   t_0},X_{t_0})+\dist_{\GH}(X_{t_0}, X_t)\\
&< & 3\delta.
\end{eqnarray*}
    Also $P(X_t, \frac{1}{2}\varepsilon)\le P(X_{t_0}, \frac{1}{4}\varepsilon)$ by
    Lemma~\ref{packing:lemma}.
    So $\dist_{\GH}(Y_{N, t}, X_t)< \varepsilon/(4P(X_t,
    \frac{1}{2}\varepsilon))$. Then Lemma~\ref{GH to ball:lemma} tells us that the
   open $\varepsilon$-balls centered at points of $Y_{N, t}$ cover
   $X_t$.
\end{proof}

We give one example to illustrate how to apply
Theorem~\ref{criterion 2 of conv:thm}. Later in
Sections~\ref{ContFieldCQMbyG:sec} and \ref{Contdeform:sec}
 we shall use
Theorem~\ref{criterion 2 of conv:thm} to study the continuity of
compact quantum metric
  spaces induced by ergodic actions (Theorem~\ref{criterion of cont field of action:thm})
and the continuity of $\theta$-deformations (Theorem~\ref{theta-deform cont:thm}).

\begin{example} \label{finite dim cont field:eg}
Let $V$ be a finite-dimensional real vector space equipped with a
distinguished element $e$. Let $T$ be a locally compact Hausdorff
space. Suppose that for each $t\in T$ there is an order-unit space
structure on $V$ with unit $e$. Denote the order-unit space for
$t$ by $(V_t, e)$ and the norm by $\pa \cdot \pa_t$. If the
function $t\mapsto \pa v\pa_t$ is continuous on $T$ for each $v\in
V$, this is called a {\it continuous field of finite-dimensional
order-unit spaces} by Rieffel \cite[Section 10]{Rieffel00}.
Clearly this fits into our Definition~\ref{cont field of
  order-unit:def}.
If there is also a Lip-norm $L_t$ for each $t$ such that $t\mapsto
L_t(v)$ is continuous on $T$ for each $v\in V$, then $\{L_t\}$ is
called a {\it continuous field of Lip-norms} \cite[Section
11]{Rieffel00}. Again, this fits into our Definition~\ref{cont
field of CQM:def}. Rieffel proved that for a continuous field of
Lip-norms, $\dist_{\q}(V_t, V_{t_0})\to 0$ as $t\to t_0$ for each
$t_0\in T$ \cite[Theorem 11.2]{Rieffel00}. By
Theorem~\ref{dist_q=dist_oq:thm} this is equivalent to saying that
$\dist_{\oq}(V_t, V_{t_0})\to 0$ as $t\to t_0$. We use
Theorem~\ref{criterion 2 of conv:thm} to give the latter a new
proof.

For later use we consider a more general case. We want to
allow $V$ to be
infinite-dimensional. To still get the continuity under $\dist_{\oq}$ we
need stronger conditions.

\begin{definition} \label{unif:def}
Let $V$ be a real vector space equipped with
a distinguished element $e$. Let $T$ be a locally compact Hausdorff space.
Suppose that for each $t\in T$ there is an order-unit space structure on
$V$ with unit $e$, for which we denote the order-unit norm by
$\pa \cdot \pa_t$.
We call  $(V, e, \{\pa \cdot \pa_t\})$
{\it a uniformly continuous field of order-unit spaces} if for
any $t_0\in T$  and $\varepsilon>0$ there is a neighborhood $\mathcal{U}$ of $t_0$
such that $(1-\varepsilon)\pa \cdot \pa_{t_0}\le \pa \cdot\pa_t\le
(1+\varepsilon)\pa \cdot \pa_{t_0}$ throughout $\mathcal{U}$.
Let $L_t$ be a Lip-norm on $(V, e, \pa \cdot\pa_t)$ for  each $t\in T$.
We call $\{L_t\}$
{\it a uniformly continuous field of Lip-norms} if for any
 $v\in V$ the function $t\mapsto L_t(v)$ is continuous,
 and if for
any $t_0\in T$  and $\varepsilon>0$ there is a neighborhood $\mathcal{U}$ of $t_0$
such that $(1-\varepsilon)L_{t_0}\le L_t$ throughout $\mathcal{U}$.
\end{definition}

Notice that we do not need $L_t\le (1+\varepsilon)L_{t_0}$.
For a continuous field of finite-dimensional order-unit spaces
(resp. finite-dimensional Lip-norms),
by a standard compactness argument we can find a neighborhood $\mathcal{U}$ of $t_0$
such that $1-\varepsilon\le \pa v\pa_t\le
1+\varepsilon$ (resp. $1-\varepsilon\le L^{\sim}_t(\tilde{v})$)
for all $t\in \mathcal{U}$ and $v$ (resp. $\tilde{v}$) in the unit sphere
of $(V, \pa \cdot \pa_{t_0})$ (resp. $(\tilde{V}, L^{\sim}_{t_0})$).
Therefore continuous fields of finite-dimensional order-unit spaces
and Lip-norms are uniformly continuous.
The assertion that $\dist_{\oq}(V_t, V_{t_0})\to 0$ as $t\to
t_0$ follows directly from  Theorem~\ref{criterion 2 of conv:thm} and the
next lemma:
\end{example}

\begin{lemma} \label{unif to cover:lemma}
Let $(V, e, \{\pa \cdot \pa_t\}, \{L_t\})$ be a uniformly continuous
field of order-unit spaces and Lip-norms over $T$. Denote the
order-unit space for $t$ by $(V_t, e)$. Then the radius function
$t\mapsto r_{V_t}$ is upper semi-continuous over $T$. Let $t_0\in
T$, and let $R> 0$. Let $\{v_n\}_{n\in \Ne}$ be a sequence dense
in $\cD_{ R}(V_{t_0})$.  Then for any $\varepsilon>0$, there
is an $N$ such that  the open $\varepsilon$-balls in $V_t$ centered
at $v_1, {\cdots}, v_N$ cover $\cD_{ R}(V_t)$ throughout
some neighborhood of $t_0$.
\end{lemma}
\begin{proof}
Let $1>\varepsilon>0$ be given.
Let $\mathcal{U}$ be a neighborhood of
$t_0$ such that $(1-\varepsilon)\pa \cdot \pa_{t_0}\le \pa \cdot\pa_t\le
(1+\varepsilon)\pa \cdot \pa_{t_0}$ and $(1-\varepsilon)L_{t_0}\le L_t$
throughout $\mathcal{U}$. Let $t\in
\mathcal{U}$. Then $(1-\varepsilon)\pa \cdot \pa^{\sim}_{t_0}\le \pa \cdot\pa^{\sim}_t\le
(1+\varepsilon)\pa \cdot \pa^{\sim}_{t_0}$ and
$(1-\varepsilon)L^{\sim}_{t_0}\le L^{\sim}_t$. By Proposition~\ref{criterion of
  Lip:prop} $\pa \cdot\pa^{\sim}_{t_0}\le r_{V_{t_0}}
L^{\sim}_{t_0}$. Thus $\pa \cdot\pa^{\sim}_t\le
(1+\varepsilon)\pa \cdot \pa^{\sim}_{t_0}\le (1+\varepsilon)r_{V_{t_0}}
L^{\sim}_{t_0}\le
\frac{1+\varepsilon}{1-\varepsilon}r_{V_{t_0}}L^{\sim}_t$. Applying Proposition~\ref{criterion of
  Lip:prop} again we see that $r_{V_t}\le
\frac{1+\varepsilon}{1-\varepsilon}r_{V_{t_0}}$. This shows that the
radius function $t\mapsto r_{V_t}$ is upper semi-continuous. Pick
$N$ such that the open $\varepsilon$-balls in $V_{t_0}$ centered at
$v_1, {\cdots}, v_N$ cover $\cD_{ R}(V_{t_0})$. Let $v\in
\cD_{ R}(V_t)$. Then $\pa v\pa_{t_0}\le
\frac{1}{1-\varepsilon}\pa v\pa_t\le \frac{1}{1-\varepsilon}R$ and
$L_{t_0}(v)\le \frac{1}{1-\varepsilon}L_t(v)\le
\frac{1}{1-\varepsilon}$. Thus $(1-\varepsilon)v\in \cD_{
  R}(V_{t_0})$. Then $\pa (1-\varepsilon)v-v_n\pa_{t_0}<\varepsilon$ for some $1\le n\le
N$. Consequently
\begin{eqnarray*}
\pa v-v_n\pa_t&\le &\pa v-(1-\varepsilon)v\pa_t+\pa
(1-\varepsilon)v-v_n\pa_t \\
&\le &\varepsilon \pa v\pa_t+
(1+\varepsilon)\pa (1-\varepsilon)v-v_n\pa_{t_0}\le \varepsilon(R+1+\varepsilon).
\end{eqnarray*}
Thus the open $(R+2)\varepsilon$-balls in $V_t$ centered at $v_1,
{\cdots}, v_N$ cover $\cD_{ R}(V_t)$.
\end{proof}


\section{Lip-norm and finite multiplicity}
\label{Lip&FinMul:sec}

In this section we prove Theorem~\ref{Lip to finite:thm} to
determine when an ergodic action of a compact group induces
a Lip-norm.

Throughout the rest of this paper $G$ will be a nontrivial compact group with
identity $e_G$, endowed with the normalized Haar measure.
Denote by $\hat{G}$ the dual of $G$, and by $\gamma_0$
the class of trivial representations.
For any $\gamma\in \hat{G}$ let  $\chi_{\gamma}$ be the corresponding character on
$G$, and let $\bar{\gamma}$ be the contragradient representation.
 For any $\gamma
\in \hat{G}$ and any representation of $G$ on some complex vector space
$V$, we denote by $V_{\gamma}$ the $\gamma$-isotypic component of $V$. If
$\mathcal{J}$ is a finite subset of $\hat{G}$, we also let
$V_{\mathcal{J}}=\sum_{\gamma\in \mathcal{J}}V_{\gamma}$, and let
$\bar{\mathcal{J}}=\{\bar{\gamma}:\gamma \in \mathcal{J}\}$.
For a strongly continuous action $\alpha$ of $G$ on a complete order-unit space
$(\bar{A}, e)$ as automorphisms, we endow $\bar{A}^{\Ce}=\bar{A}\otimes_{\Re}\Ce=\bar{A}+i\bar{A}$
with the diagonal action $\alpha^{\Ce}:=\alpha\otimes I$. We
say that $\alpha$ is  of {\it finite multiplicity}
if $\mul(\bar{A}^{\Ce}, \gamma)<\infty$ for all $\gamma
\in \hat{G}$, and that $\Gamma$ is \emph{ergodic} if the only $\alpha$-invariant
elements are the scalar multiples of $e$.

We also fix a length function on $G$, \ie a continuous
real-valued function, $\mathnormal{l}$, on $G$ such that
\begin{eqnarray*}
\mathnormal{l}(xy)&\le &\mathnormal{l}(x)+\mathnormal{l}(y)\mbox{ for
  all } x, y\in G \\
\mathnormal{l}(x^{-1})&=& \mathnormal{l}(x) \mbox{ for all }x\in G \\
\mathnormal{l}(x)&=& 0 \mbox{ if and only if } x=e_G.
\end{eqnarray*}

\begin{remark} \label{length func:remark}
One can verify easily that
a length function $\mathnormal{l}$ on $G$ is equivalent to
a left invariant metric $\rho$ on $G$ under the correspondence
$\rho(x, y)=\mathnormal{l}(x^{-1}y)$. Since every metric on a compact
group could be integrated to be a left invariant one,
we see that a compact group $G$ has a
length function if and only if it is metrizable.
\end{remark}

Let $\mathcal{A}$ be a unital
$C^*$-algebra, and let $\alpha$ be a strongly continuous ergodic action of $G$ on
$\mathcal{A}$ by automorphisms. In \cite{Rieffel98b} Rieffel defined a (possibly $+\infty$-valued)
seminorm $L$ on $\mathcal{A}$ by
\begin{eqnarray} \label{def of L:eq}
L(a)=\sup \{\frac{\pa \alpha_x(a)-a\pa}{\mathnormal{l}(x)}: x\in G, x\neq e_G\}.
\end{eqnarray}
He showed that the set $A=\{a\in \mathcal{A}_{\sa}:L(a)<\infty\}$
is a dense subspace of $\mathcal{A}_{\sa}$ containing the identity
$e$, and that $(A, L|A)$ is a closed compact quantum metric space \cite[Theorem 2.3]{Rieffel98b}.
In fact, the proof
there shows more:

\begin{theorem}\cite[Theorem 2.3]{Rieffel98b} \label{finite to Lip:thm}
Let $\alpha$ be a strongly continuous isometric action of $G$ on a
(real or complex) Banach space $\bar{V}$. Define a (possibly
$+\infty$-valued) seminorm $L$ on $\bar{V}$ by (\ref{def of L:eq}).
Then $V:=\{v \in \bar{V}:L(v)<\infty\}$ is always a dense subspace
of $\bar{V}$. If $\bar{V}=(\bar{A}, e)$ is a complete order-unit
space and $G$ acts as automorphisms of $\bar{A}$, then $A=\{a \in
\bar{A}:L(a)<\infty\}$ also contains $e$, and hence we can
identify $S(A)$ with $S(\bar{A})$.  If furthermore $\alpha$ is
ergodic and of finite multiplicity, then $A$ with the restriction
of $L$ is a closed compact quantum metric space, and $r_A\le
\int_G\mathnormal{l}(x)\, dx$.
\end{theorem}

The aim of this section is to show that the converse of
Theorem~\ref{finite to Lip:thm} is also true:

\begin{theorem} \label{Lip to finite:thm}
Let $\alpha$ be an ergodic strongly continuous action of $G$ on a
complete unit-order space $(\bar{A}, e)$. Define $L$ and $A$
as in Theorem~\ref{finite to Lip:thm}. If the restriction
of $L$ on $A$ is a Lip-norm, then $\alpha$ is of finite multiplicity.
\end{theorem}

The intuition is that covering numbers of $\cD_r(A)$
increase (fast) as the multiplicities $\mul(\bar{A}^{\Ce}, \Gamma)$
increase. Thus the compactness of  $\cD_r(A)$
in Proposition~\ref{criterion of Lip:prop} forces $\mul(\bar{A}^{\Ce}, \Gamma)$ to be finite.

\begin{lemma} \label{bound Lip for omega:lemma}
For any finite subset $\mathcal{J}$ of $\hat{G}$ and
any map $\omega:\mathcal{J}\rightarrow
\Ne\cup \{0\}$, there is a constant $M_{\mathcal{J}, \omega}> 0$ such that
for any strongly continuous isometric action $\alpha$ of $G$ on a
finite-dimensional
complex Banach space $V$ with $\mul(V, \gamma)\le \omega(\gamma)$ for
$\gamma \in \mathcal{J}$ and $\mul(V, \gamma)=0$ for
$\gamma \in \hat{G}\setminus \mathcal{J}$, we have $L\le
M_{\mathcal{J}, \omega}\pa \cdot \pa$ on $V$.
\end{lemma}
\begin{proof}
Let $X$ be the set of all functions $\omega_1:\mathcal{J}\rightarrow
\Ne \cup \{0\}$ with $\omega_1\le \omega$. Set
$\omega_1(\gamma)=\mul(V, \gamma)$
for all $\gamma \in \mathcal{J}$.
Let $\pa \cdot \pa_V$ be the norm on $V$, and let $N=\dim
  V$.
Let $W$ be a Hilbert space with dimension $N$.
It is a theorem of John \cite[Proposition 9.12]{Tomczak89}
that there is a linear isomorphism $\phi:V\to W$ such that
$\pa \phi\pa, \, \pa \phi^{-1}\pa\le \root 4 \of N$.
Define an inner product, $<, >_*$, on $V$ by
$<u, v>_*=\int_G<\phi(\alpha_x(u)), \phi(\alpha_x(v))>\, dx$.
Then $<, >_*$ is $G$-invariant. Let $\pa \cdot \pa_*$ be the
corresponding norm. Since $\alpha$ is isometric with respect to $\pa
\cdot \pa_V$, for any $v\in V$ we have
$\pa v\pa^2_*=\int_G\pa \phi(\alpha_x(v))\pa^2\, dx\le
\pa \phi\pa^2 \int_G\pa \alpha_x(v)\pa^2_V\, dx=\pa \phi\pa^2\pa
v\pa^2_V\le \sqrt{N}\pa v\pa^2_V$. Thus $\pa \cdot \pa_*\le \root 4 \of N \pa \cdot \pa_V$.
Similarly, $\pa \cdot \pa_*\ge \frac{1}{\root 4 \of N}\pa \cdot \pa_V$.

For each $\gamma \in \hat{G}$ fix a Hilbert space $H_{\gamma}$
with an irreducible unitary action of type $\gamma$.
Let $H$ be the Hilbert space direct sum $\oplus_{\gamma \in
  \mathcal{J}}\oplus^{\omega_1(\gamma)}_{j=1}H_{\gamma}$ equipped
  with the natural action $\beta$ of $G$. By Theorem~\ref{finite to Lip:thm} $L$ is finite
on a dense subspace of $H$. Since
$H$ is finite dimensional, $L_{\rH}$ is finite on the
whole $H$. Clearly there is
a constant $M$ such that $L_{\rH}\le M \pa \cdot \pa_{\rH} $
on $H$.

Since $(V, <, >_*)$ has the same
  multiplicities as does $H$, there is a $G$-equivariant
  unitary map $\varphi: (V, <, >_*)\to H$.
Then for any $v\in V$ and $x\in G$ we have that
\begin{eqnarray*}
\pa v-\alpha_x(v)\pa &\le & \root 4 \of N \pa v-\alpha_x(v)\pa_*=\root
4 \of N \pa
\varphi(v)-\varphi(\alpha_x(v))\pa_{\rH} \\
&=&\root 4 \of N \pa
\varphi(v)-\beta_x(\varphi(v))\pa_{\rH}
\le M\cdot \root 4 \of N \cdot \mathnormal{l}(x)\pa
\varphi(v)\pa_{\rH} \\
&=& M\cdot \root 4 \of N \cdot \mathnormal{l}(x)\pa
v\pa_*
\le M\cdot \sqrt{N}\cdot \mathnormal{l}(x)\pa v\pa_V.
\end{eqnarray*}
Let $M'_{\mathcal{J},\omega_1}=M\cdot \sqrt{N}$.
Then $L_V\le M_{\mathcal{J}, \omega}\pa \cdot \pa_V$ on $V$.
Thus $M_{\mathcal{J},\omega}:=1+\max\{M'_{\mathcal{J},\omega_1}:\omega_1
  \in X\}$ satisfies the requirement.
\end{proof}

In particular, let $M_{\gamma}$ be the constant $M_{\mathcal{J}, \omega}$ for
$\mathcal{J}=\{\gamma\}$ and $\omega(\gamma)=1$.

We shall need a well-known
fact (cf. the discussion at the end of page 217
of \cite{RS90}, noticing
that $\delta$ can be $1$ in Lemma 2 when $E_1$ is  finite-dimensional):

\begin{lemma} \label{seq of dist 1:lemma}
For any (real or complex) normed space $V$ there is a sequence
$p_1, p_2, \cdots
$ in the unit sphere of $V$,
with length $\dim(V)$ when $V$ is finite-dimensional or length $\infty$
otherwise,
satisfying that $\pa p_m-q\pa \ge 1$ for all
$m$  and all $q\in \spn\{p_1, {\cdots}$, $p_{m-1}\}$.
\end{lemma}

\begin{lemma}  \label{far points:lemma}
Let $\alpha$ be a strongly continuous isometric action
 of $G$ on a complex Banach space $V$, and let
$\gamma \in \hat{G}$.
Then there is a subset $X\subseteq \cD_{1/M_{\gamma}}(V)=\{v\in V:L(v)\le 1, \pa v\pa \le 1/M_{\gamma}\}$,
with $\mul(V, \gamma)$ many elements when $\mul(V, \gamma)<\infty$ or
 infinitely many elements otherwise,
such that any two distinct points in $X$ have distance no less than $1/M_{\gamma}$.
\end{lemma}
\begin{proof} Fix a Banach space $H_{\gamma}$ with an irreducible action of type $\gamma$.
We have $V_{\gamma}=\oplus^{\mul(V, \gamma)}_{j=1}V_j$ with
G-equivariant isomorphisms $\varphi_j:H_{\gamma}\to V_j$.
Take a nonzero $u$ in $H_{\gamma}$, and let $W=\spn\{ \varphi_j(u)\}_j$.
Then for any nonzero $v$ in $V'$, $v$ is "purely" of type $\gamma$,
\ie the action of $G$ on
$\spn\{\alpha_x(v)\}_{x\in G}$ is an irreducible one of type $\gamma$.
By Lemma~\ref{seq of dist 1:lemma} we can  find a subset $Y$ in the
unit sphere of $W$, with $\dim(W)$ many elements when $W$ is  finite-dimensional or
 infinitely many elements otherwise,
such
that any two distinct points in $Y$ have distance no less than $1$.
According to Lemma~\ref{bound Lip for omega:lemma} any element in the
unit sphere of $W$ has $L$ no bigger than $M_{\gamma}$.
Therefore $Y/M_{\gamma}\subseteq \cD_{1/M_{\gamma}}(V)$.
\end{proof}

For a set $S$ and  a subset $X$ of $S$ we say that $X$ is an $n$-subset
if $X$ consists of $n$ elements. For $q_1, q_2\ge 2$
let $N(q_1, q_2;2)$ be the Ramsey number \cite{VW92}, \ie the minimal
number $n$ such that for any set $S$ with
at least $n$ elements, if the set of all $2$-subsets of $S$ is
 divided into $2$ disjoint families $Y_1$ and $Y_2$ (``colors''), then
there are a $j$ and some $q_j$-subset of $S$ for which every $2$-subset
is in $Y_j$. Consequently, for any set $S$ with
at least $n$ elements, if the set of all $2$-subsets of $S$ is
 the union of $2$ (not necessarily disjoint) families $T_1$ and $T_2$, then
there are a $j$ and some $q_j$-subset of $S$ for which every $2$-subset
is in $T_j$.

\begin{lemma}  \label{Cov bound:lemma}
Let $\alpha$ be a strongly continuous action of $G$ on a
 complete order-unit space $(\bar{A}, e)$ as automorphisms.
Suppose that for some $\gamma \in \hat{G}\setminus \{\gamma_0\}$ and $q>0$,
$\mul(\bar{A}^{\Ce}, \gamma)\ge N(q, q;2)$. Then
for any $0< \varepsilon <1/(4M_{\gamma})$ we have that
$\Cov(\cD_{1/M_{\gamma}}(A), \varepsilon)\ge q$.
\end{lemma}
\begin{proof} Since $\bar{A}=\Af_{\Re}(S(\bar{A}))$, we can identify
$\bar{A}^{\Ce}$ with $\Af_{\Ce}(S(\bar{A}))$, the space of
$\Ce$-valued continuous affine functions on $S(\bar{A})$ equipped with the supremum norm,
 and hence $\bar{A}^{\Ce}$
becomes a complex Banach space whose norm extends that of $\bar{A}$.
Notice that the action $\alpha$ corresponds to an action
$\alpha'$ of $G$ on $S(\bar{A})$. Then clearly $\alpha^{\Ce}$ is
isometric and strongly continuous with respect to this norm.

According to Lemma~\ref{far points:lemma} we can find a subset
$X\subseteq \cD_{1/M_{\gamma}}(\bar{A}^{\Ce})$,
with $\mul(\bar{A}^{\Ce}, \gamma)$ many elements when
$\mul(\bar{A}^{\Ce}, \gamma)<\infty$ or infinitely many elements otherwise,
such that any two distinct points in $X$ have distance at least $1/M_{\gamma}$.
 Then $|X|\ge N(q,q;2)$.
For any distinct $a_1+ia_2,\, b_1+ib_2\in X$ with $a_j, b_j\in \bar{A}$,
we have that $\pa (a_1+ia_2)-(b_1+ib_2)\pa \ge 1/M_{\gamma}$, and hence
$|a_1-b_1|\ge 1/(2M_{\gamma})$ or $|a_2-b_2|\ge 1/(2M_{\gamma})$.
Denote by $X^{(2)}$ the set of all $2$-subsets of $X$. Let
$T_j=\{\{a_1+ia_2,b_1+ib_2\}\in X^{(2)}:|a_j-b_j|\ge 1/(2M_{\gamma})\}$.
 Then
$T_1\cup T_2=X^{(2)}$. So there are a $j$ and some $q$-subset $X'$
of $X$ for which every $2$-subset is in $T_j$.
Let $Y=\{a_j: a_1+ia_2\in X'\}$.
Clearly $Y$ is contained in $\cD_{1/M_{\gamma}}(A)$
and any two distinct points in $Y$ have distance at least
$1/(2M_{\gamma})$.

Let $0< \varepsilon <1/(4M_{\gamma})$. Suppose that
$p_1, {\cdots}, p_k\in \cD_{1/M_{\gamma}}(A)$ and the open
$\varepsilon$-balls centered at $p_1, {\cdots}, p_k$ cover
$\cD_{1/M_{\gamma}}(A)$. Then each such open ball could
contain
at most one point in $Y$. So $k\ge |Y|=q$, and hence
$\Cov(\cD_{1/M_{\gamma}}(A),\varepsilon)\ge q$.
\end{proof}

\begin{proof}[Proof of Theorem~\ref{Lip to finite:thm}]
Suppose that $\mul(\bar{A}^{\Ce}, \gamma)=\infty$ for some $\gamma \in
\hat{G}$. For any $0< \varepsilon
<1/(4M_{\gamma})$ by Lemma~\ref{Cov bound:lemma} we have
$\Cov(\cD_{1/M_{\gamma}}(A), \varepsilon)=\infty$. Therefore $\cD_{1/M_{\gamma}}(A)$
is not totally bounded. By Proposition~\ref{criterion of
Lip:prop} $L$ is not a Lip-norm on $A$.
\end{proof}

It should be pointed out that there do exist examples of ergodic
strongly continuous action of $G$ on a complete unit-order space $(\bar{A}, e)$,
for which $\mul(\bar{A}^{\Ce}, \gamma)=\infty$ for some
$\gamma \in \hat{G}$.
 We shall give such an example in Section~\ref{comp&bm:sec}.

\section{Compactness and bounded multiplicity}
\label{comp&bm:sec}

In this section we investigate when a family of compact quantum
metric spaces induced from ergodic actions of $G$ is totally
bounded.

In the discussion after Theorem 13.5 of \cite{Rieffel00} Rieffel
observed that the set of all isometry classes of compact quantum
metric spaces (for given $\mathnormal{l}$) induced from ergodic
actions of $G$ on unital $C^*$-algebras is totally bounded under
$\dist_{\q}$. In fact, the argument there works for general ergodic
actions of $G$ on complete order-unit spaces:

\begin{theorem}\cite[Section 13]{Rieffel00}  \label{bounded to compact:thm}
Let $\mathcal{S}$ be a set of compact quantum metric spaces $(A,
L)$ induced by ergodic actions $\alpha$ of $G$ on $\bar{A}$ for a
fixed $\mathnormal{l}$. Let $\mul(\mathcal{S}, \gamma)=\sup
\{\mul(\bar{A}^{\Ce}, \gamma):(A, L)\in \mathcal{S}\}$ for each
$\gamma \in \hat{G}$. If $\mul(\mathcal{S}, \gamma)<\infty$ for all
$\gamma \in \hat{G}$, then $\mathcal{S}$ is totally bounded under
$\dist_{\q}$.
\end{theorem}

We show that the converse is also true:

\begin{theorem} \label{compact to bounded:thm}
Let $\mathcal{S}$ be a set of compact quantum metric spaces $(A,
L)$ induced by ergodic actions $\alpha$ of $G$ on $\bar{A}$ for a
fixed $\mathnormal{l}$. If $\mathcal{S}$ is totally bounded under
$\dist_{\q}$, then $\mul(\mathcal{S}, \gamma)< \infty$ for all $\gamma
\in \hat{G}$.
\end{theorem}
\begin{proof}
Suppose that $\mathcal{S}$ is totally bounded. For any $R>\int_G
 \mathnormal{l}(x)\, dx$,
by Theorems~\ref{QGH 2:thm}, \ref{finite to Lip:thm}, and Lemma~\ref{com equi of dual:lemma}
 we have that $\sup\{\Cov(\cD_{ R}(A), \varepsilon): (A, L)\in
 \mathcal{S}\}<\infty$ for all $\varepsilon>0$.
Taking $R>\max(\int_G
 \mathnormal{l}(x)\, dx, 1/M_{\gamma})$, by Lemma~\ref{cov of
 sub:lemma} we get $\sup \{\Cov(\cD_{1/M_{\gamma}}(A),\varepsilon):
(A, L)\in \mathcal{S}\}<\infty$ for all $\varepsilon>0$.
Let
$$M=\sup \{\Cov(\cD_{1/M_{\gamma}}(A),1/(5M_{\gamma})):
(A, L)\in \mathcal{S}\}.$$
Then Lemma~\ref{Cov bound:lemma} tells us that
$\mul(\mathcal{S}, \gamma)<N(M+1, M+1;2)$.
\end{proof}

It should be pointed out that there do exist ergodic actions of $G$
on complete order-unit spaces with big multiplicity:

\begin{example}  \label{large mul:eg}
 Let $\{(A_j, L_j)\}_{j\in I}$ be a family of compact quantum
metric spaces induced by ergodic actions $\alpha_j$ of $G$ on
$(\bar{A_j}, e_j)$.
Also let $\eta_j(a)=\int_G (\alpha_j)_x(a)\, dx$ be the unique
$G$-invariant state on $\bar{A_j}$ and $V_j=\ker (\eta_j)$.
Then $\bar{A_j}=\Re e_j\oplus V_j$ as vector spaces.
As in Section 12 of \cite{Rieffel00},
consider $\prod^b\bar{A_j}$ the subspace of the full product
which consists of sequences $\{a_j\}$ for which $\pa
  a_j\pa$ is bounded.
This is a complete order-unit space with unit $e=\{e_j\}$.
Consider the reduced product $\prod^{rb}\bar{A_j}
=\{\{a_j\}\in \prod^b\bar{A_j}:\eta_j(a_j)=\eta_k(a_k) \,
\mbox{ for all } j,k\in I\}$.
Then $\prod^{rb}\bar{A_j}$ is a closed subspace in
$\prod^b\bar{A_j}$, and is also a complete order-unit space.
Clearly $\prod^{rb}A_j=\Re e\oplus \prod^{b}V_j$ as vector spaces.
The actions $\alpha_j$ of $G$ on the components $\bar{A_j}$ give an isometric
action on  $\prod^{b}\bar{A_j}$, which we denote by $\alpha$.
Although $\alpha$ is not ergodic on $\prod^{b}\bar{A_j}$,
it is on $\prod^{rb}\bar{A_j}$ because of the above decomposition
as a direct sum.
By the natural $G$-equivariant embedding
$\bar{A_j}\hookrightarrow \prod^{b}\bar{A_j}$,
we see that $\mul((\prod^{rb}\bar{A_j})^{\Ce}, \gamma)\ge
\sum_{j\in I}\mul(\bar{A_j}^{\Ce}, \gamma)$ for
every $\gamma \in \hat{G}\setminus \{\gamma_0\}$.

In general, $\alpha$ may not be strongly continuous on
$\prod^{rb}\bar{A_j}$. But there are two special cases in which it is strongly continuous:

(1). When $(\bar{A_j}, e_j)$  and $\alpha_j$
are all the same and finite dimensional,
say $(\bar{A_j}, e_j)=(\bar{A},e)$.
Then $\bar{A}=A$ and there is some constant $M>0$ such that
$L\le M \pa \cdot \pa$ on $A$.
Therefore  $\pa a-(\alpha_j)_x(a)\pa \le L(a)\mathnormal{l}(x)\le
M\pa a\pa\mathnormal{l}(x)$ for all $a\in A$ and $x\in G$.
Then it is easy to see that $\alpha$ is  strongly continuous on
$\prod^{rb}A_j$.
It is standard \cite{Sugiura90} that in the left regular representation
of $G$ on $C(G)$ the multiplicity of each
$\gamma \in \hat{G}$
equals $\dim(\gamma)$. For a finite subset
$\mathcal{J}\subseteq \hat{G}$
and any $a+ia'\in C(G)$, clearly $a+ia'\in
(C(G))_{\mathcal{J}}$ if and only if $a-ia'\in (C(G))_{\bar{\mathcal{J}}}$.
Therefore if $\mathcal{J}=\bar{\mathcal{J}}$ and $\gamma_0\in
\mathcal{J}$, then
$(C(G))_{\mathcal{J}}$ is closed under the involution and
  contains the constant functions. Hence $((C(G))_{\mathcal{J}})_{\sa}$ is a complete
finite-dimensional order-unit space and $(C(G))_{\mathcal{J}}
=(((C(G))_{\mathcal{J}})_{\sa})^{\Ce}$.
Taking $I=\Ne$ and $A=((C(G))_{\mathcal{J}})_{\sa}$, we get
$\mul((\prod^{rb}A_j)^{\Ce}, \gamma)=\infty$
for every $\gamma \in\mathcal{J}\setminus \{\gamma_0\}$ as promised at the end
of Section~\ref{Lip&FinMul:sec}.

(2). When $I$ is finite, $\alpha$ is always strongly continuous on
$\prod^{rb}\bar{A_j}$.
In particular, take $\bar{A}=(C(G))_{\sa}$,
 and $\bar{A_j}=\bar{A}$ for all $j\in I$.
Since $\mul(\bar{A}^{\Ce}, \gamma)$ equals
$\dim(\gamma)$
for all $\gamma \in \hat{G}$,
we see that $\mul((\prod^{rb}\bar{A_j})^{\Ce}, \gamma)$
could be as big as we want for any $\gamma\neq \gamma_0$.
\end{example}


\section{Continuous fields of compact quantum metric spaces induced by
ergodic compact group actions} \label{ContFieldCQMbyG:sec}

In this section we study continuous fields of compact quantum
metric spaces induced by ergodic actions, and prove
Theorem~\ref{criterion of cont field of action:thm}.
In Examples~\ref{tori 2:eg} and \ref{Berezin 2:eg} we use
Theorem~\ref{criterion of cont field of action:thm} to give a
unified treatment of the continuity of noncommutative tori and
integral coadjoint orbits, which were studied by Rieffel before.

Rieffel has defined continuous fields of actions of a locally
compact group on $C^*$-algebras \cite[Definition 3.1]{Rieffel89b}.
We adapt it to actions on order-unit spaces:

\begin{definition}  \label{cont field of action:def}
Let $(\{\overline{A_t}\}, \Gamma)$ be a continuous field
of order-unit spaces over a locally compact Hausdorff space $T$,
and let $\alpha_t$ be a strongly continuous action of $G$ on $\overline{A_t}$
for each $t\in T$. We say that $\{\alpha_t\}$ is a {\it continuous field
of strongly continuous actions} of $G$ on $(\{\overline{A_t}\}, \Gamma)$
if the action of $G$ on $\Gamma_{\infty}$
is strongly continuous, where $\Gamma_{\infty}:=\{f\in \Gamma:\mbox{ the function } t\mapsto\pa f_t\pa
\mbox{ vanishes at } \infty\}$ is the space of continuous sections
vanishing at $\infty$.
If each $\alpha_t$ is ergodic, we say that this is a \emph{field of ergodic
actions}. If each $\alpha_t$ is of finite multiplicity,
we say that this is a \emph{field of finite actions}.
\end{definition}

\begin{remark} \label{norm of field:remark}
For a continuous field $(\{\overline{A_t}\}, \Gamma)$ of order-unit
spaces over a compact Hausdorff space $T$, it is easy to see that
$\Gamma$ is a complete order-unit space with the unit $e$ and
the order defined as
$f\ge g$ if and only if $f_t\ge g_t$ for all $t\in T$. Then the
natural projections $\Gamma\rightarrow \overline{A_t}$ become
order-unit space quotient maps. According to the discussion right after Proposition~\ref{criterion of Lip:prop}
we may identify
$S(\overline{A_t})$ with a closed convex subset of $S(\Gamma)$.
From the definition of the order in $\Gamma$ it is easy to see that
the convex hull of the union of all the $S(\overline{A_t})$'s is dense in $S(\Gamma)$.
If we identify $\overline{A_t}^{\Ce}$ and
$\Gamma^{\Ce}$  with $\Af_{\Ce}(S(\overline{A_t}))$ and
$\Af_{\Ce}(S(\Gamma))$ respectively, then they are endowed with
complex vector space norms and $\pa f+ig\pa=\sup\{\pa
f_t+ig_t\pa:t\in T\}$ for all $f+ig\in \Gamma^{\Ce}$. When we talk
about $\overline{A_t}^{\Ce}$ and $\Gamma^{\Ce}$ as complex Banach
spaces, we always mean these norms. If $(\alpha_t)$ is a continuous field of
strongly continuous actions of $G$ on $(\{\overline{A_t}\}, \Gamma)$,
then the action of $G$ on $\Gamma^{\Ce}$ is easily seen to be
also strongly continuous.
\end{remark}

For a continuous field of strongly continuous finite ergodic
actions of $G$ on order-unit spaces, obviously we get a field of
compact quantum metric spaces.
Theorem~\ref{criterion of cont field of action:thm} indicates
that this is indeed
a continuous field, as one may
expect. However, as Theorems~\ref{criterion of conv:thm} and
 ~\ref{criterion 2 of conv:thm} indicate, as $t\to t_0$ the
corresponding compact quantum metric spaces do not necessarily
converge to that at $t_0$. We give a trivial example here:

\begin{example} \label{not conv:eg}
Take a complete order-unit space $\bar{A}$ with a nontrivial
ergodic action of $G$ with
finite multiplicity (for example, $(C(G))_{\sa}$ with the left regular
representation of $G$). Let $T=[0,1]$. Then we have
the trivial field $(\{\overline{A^1_t}\}, \Gamma^1)$ with
$\overline{A^1_t}=\bar{A}$ for all $t\in T$. The action of
$G$ on $\Gamma^1$ is clearly strongly continuous. Now we take the subfield
$((\overline{A_t}), \Gamma)$ with $\overline{A_t}=\bar{A}$ for all $0<t\le 1$
 and $\overline{A_0}=\Re e_A$. The action of $G$ restricted on
$\Gamma$ is still strongly continuous. But $r_{A_0}=0$ and
$r_{A_t}=r_A>0$ for all $0<t\le 1$. So $\dist_{\oq}(A_t, A_0)=\dist_{\oq}(A, \Re
e_A)$
does not converge to $0$ as $t\to 0$.
\end{example}

Notice that in the above example the multiplicities degenerate at
$t_0=0$. Theorem~\ref{criterion of cont field of action:thm} tells us that
this is exactly why we do not get $\dist_{\oq}(A_t, A_0)\to 0$.

We start to prove Theorem~\ref{criterion of cont field of action:thm}.
We show first that the multiplicity function $t\mapsto \mul(\overline{A_t}^{\Ce},\gamma)$
is lower semi-continuous. This shows
(ii)$\Longrightarrow$(i) in Theorem~\ref{criterion of cont field of action:thm}.

\begin{lemma} \label{Ind open:lemma}
Let $(\{V_t\}, \Gamma)$ be a continuous field of (real or complex)
Banach spaces over a locally compact Hausdorff space $T$.
For any $f_1, {\cdots}, f_m\in \Gamma$
the set $\Ind(f_1, {\cdots}, f_m)=
\{t \in T: ((f_1)_t, {\cdots}, (f_m)_t) \mbox{ are linearly independent}\}$
is open.
\end{lemma}
\begin{proof}
Let $t_0\in \Ind(f_1, {\cdots}, f_m)$. Then for any $t\in T$ the map
$(f_j)_{t_0}\mapsto (f_j)_t, \, j=1, {\cdots}, m$, extends uniquely to a
linear map $\varphi_t$ from $W:=\spn\{(f_1)_{t_0}, {\cdots}, (f_m)_{t_0}\}$
to $V_t$. Notice that for any $v\in W$
the section $t\mapsto \varphi_t(v)$ is in $\Gamma$.
A standard compactness argument shows that there is a neighborhood
$\mathcal{U}$ of $t_0$ such that $\frac{1}{2}<\pa \varphi_t(v)\pa$ for
all $t\in \mathcal{U}$ and $v$ in the unit sphere of $W$. In particular,
$\varphi_t$ is injective throughout $\mathcal{U}$. Thus $(f_1)_t,
{\cdots}, (f_m)_t$ are linear independent throughout $\mathcal{U}$.
\end{proof}

We shall need the following well-know fact
several times. We omit the proof.
\begin{lemma} \label{proj TVS:lemma}
Let $G$ be a compact group.
Let $\alpha$ be a continuous action of $G$ on a complex
Banach space
 $V$.  For a continuous
$\Ce$-valued function $\varphi$ on $G$ let
\begin{eqnarray*}
\alpha_{\varphi}(v)=\int_G \varphi(x)\alpha_x(v)\, dx
\end{eqnarray*}
for $v\in V$. Then $\alpha_{\varphi}:V\rightarrow V$ is a continuous linear
map. If $\mathcal{J}$ is a finite subset of $\hat{G}$ and if $\varphi$ is
a linear combination of the characters of $\gamma\in
\bar{\mathcal{J}}$, then $\alpha_{\varphi}(V)\subseteq V_{\mathcal{J}}$.
Let
\begin{eqnarray*}
\alpha_{\mathcal{J}}=\alpha_{\sum_{\gamma \in \mathcal{J}}\dim(\gamma)
\overline{\chi_{\gamma}}}.
\end{eqnarray*}
(When $\mathcal{J}$ is a one-element set $\{ \gamma\}$, we will simply
write $\alpha_{\gamma}$ for $\alpha_{\{\gamma\}}$.)
Then
$\alpha_{\mathcal{J}}(v)=v$ for all $v\in V_{\mathcal{J}}$, and
$\alpha_{\mathcal{J}}(v)=0$ for all $v\in V_{\gamma}$ with $\gamma \in
\hat{G}\setminus \mathcal{J}$.
\end{lemma}

\begin{lemma} \label{ls of mul:lemma}
Let $\{\alpha_t\}$ be a continuous field of strongly continuous actions of
$G$ on a continuous field of order-unit spaces
$(\{\overline{A_t}\}, \Gamma)$
over a compact Hausdorff space $T$. Then for any $\gamma \in
\hat{G}$
the multiplicity function (possibly $+\infty$-valued)
$t\mapsto \mul(\overline{A_t}^{\Ce}, \gamma)$ is lower semi-continuous over $T$.
For any finite subset $\mathcal{J}$  of $\hat{G}$ and
$v\in (\overline{A_t}^{\Ce})_{\mathcal{J}}$ we can lift $v$ to $f$ in
$\Gamma^{\Ce}_{\mathcal{J}}$.
If furthermore
$\mathcal{J}=\overline{\mathcal{J}}$ and $v$ is in
$A_t$, we may take
$f$ to be in $\Gamma$.
\end{lemma}
\begin{proof} Let $\alpha$ be the action of $G$ on $\Gamma$.
Suppose that $v_1, {\cdots}, v_m$ are linearly independent vectors in
$(\overline{A_t}^{\Ce})_{\mathcal{J}}$. Let $f_1, {\cdots}, f_m$ be
lifts of $v_1, {\cdots}, v_m$ in $\Gamma^{\Ce}$. Let
$\alpha^{\Ce}_{\mathcal{J}}$ and $(\alpha^{\Ce}_t)_{\mathcal{J}}$
be the maps for $\alpha^{\Ce}$ and $\alpha^{\Ce}_t$ as defined in
Lemma~\ref{proj TVS:lemma}. Let
$\tilde{f}_j=\alpha^{\Ce}_{\mathcal{J}}(f_j)$. Then
$\tilde{f_1}, {\cdots}, \tilde{f_m}$ are in
$\Gamma^{\Ce}_{\mathcal{J}}$. Since the projection
$\Gamma^{\Ce}\to \overline{A_t}^{\Ce}$ is $G$-equivariant, by
\cite[Lemma 3.3]{Li9} we have
$(\alpha^{\Ce}_{\mathcal{J}}(f_j))_t=(\alpha^{\Ce}_t)_{\mathcal{J}}((f_j)_t)$.
Thus
$(\tilde{f_j})_t=(\alpha^{\Ce}_{\mathcal{J}}(f_j))_t=(\alpha^{\Ce}_t)_{\mathcal{J}}((f_j)_t)=
(\alpha^{\Ce}_t)_{\mathcal{J}}(v_j)=v_j$.

By Lemma~\ref{Ind open:lemma}  $(\tilde{f_1})_{t'},
{\cdots},(\tilde{f_m})_{t'}$ are linearly independent in
some open neighborhood $\mathcal{U}$ of $t$. Since $\tilde{f_i}\in
\Gamma^{\Ce}_{\mathcal{J}}$, $(\tilde{f_i})_{t'}\in
(\overline{A_{t'}}^{\Ce})_{\mathcal{J}}$ for all $t'\in T$.
Taking $\mathcal{J}=\{\gamma\}$, we get the lower semi-continuity of the
multiplicity function.

If furthermore $\mathcal{J}=\bar{\mathcal{J}}$ and
$v_1, {\cdots}, v_m$ are all in
$A_{t}$, we may take
$f_1, {\cdots}, f_m$ to be all in $\Gamma$.
Notice that $\alpha^{\Ce}_{\mathcal{J}}=\alpha^{\Ce}_{\varphi}$, where
$\varphi=\sum_{\gamma \in \mathcal{J}}\dim(\gamma)
\overline{\chi_{\gamma}}$.
Since the function $x\mapsto (\sum_{\gamma \in
  \mathcal{J}}\dim(\gamma) \overline{\chi_{\gamma}})(x)$ is real-valued in this
case, we see that $\tilde{f_1},{\cdots}, \tilde{f_m}$ are also in $\Gamma$.
\end{proof}

Next we show that there are enough Lipschitz sections.
Recall that a vector $f\in \Gamma^{\Ce}$ is called {\it $G$-finite} if
the linear span of its orbit under $\alpha$ is finite dimensional.
We will show that the Lip-norm function $t\mapsto L_t(f_t)$ is continuous
for $G$-finite $f$.
The case of this for quantum tori is proved by Rieffel in Lemma 9.3 of
\cite{Rieffel00}. Our proof
for the general case follows the way given there.

\begin{lemma} \label{cont of Lip:lemma}
Let $(\alpha_t)$ be a continuous field of strongly continuous actions of
$G$ on a continuous field of
order-unit spaces
$(\{\overline{A_t}\}, \Gamma)$
over a compact Hausdorff space $T$.
Then for any $G$-finite $f\in \Gamma^{\Ce}$ the function
$t \to L_t(f_t)$ takes finite values and
is continuous on $T$.
\end{lemma}
\begin{proof} Let $L$ be the seminorm on $\Gamma^{\Ce}$ defined
by (\ref{def of L:eq}).
Let $f\in \Gamma^{\Ce}$ be $G$-finite.
By Theorem~\ref{finite to Lip:thm} $L$ takes finite value on
a dense subspace of $\spn\{\alpha_x(f):x\in G\}$.
Since $\spn\{\alpha_x(f):x\in G\}$ has finite dimension,
$L$ is finite on the whole $\spn\{\alpha_x(f):x\in G\}$.
It is clear that $L(f)=\sup_{t\in T}L_t(f_t)$.
So $t\mapsto L_t(f_t)$ is a real-valued function.

Let
\begin{eqnarray*}
D_f=\{ (\alpha_x(f)-f)/\mathnormal{l}(x): x\neq e_G\}.
\end{eqnarray*}
Then $\sup_{g\in D_f}\pa g\pa =L(f)$. So $D_f$ is a bounded subset of
the finite dimensional space $\spn\{\alpha_x(f):x\in G\}$, and
hence totally bounded.

For $g\in \Gamma^{\Ce}$, let $F_g(t)=\pa g_t \pa$ on $T$.
Then clearly $\pa F_g-F_h\pa_{\infty} \le \pa g-h\pa$, and hence $F$ as a map from
 $\Gamma^{\Ce}$ to $C(T)$ is Lipschitz.
Therefore $F(D_f)$ is totally bounded in $C(T)$. By the
Arzela-Ascoli theorem \cite{Conway90} $F(D_f)$ is equicontinuous.
Since the supremum of a family of equicontinuous functions is
continuous, we see that the function $(t\mapsto
L_t(f_t))=\sup_{g\in D_f}F(g)$ is continuous on $T$.
\end{proof}

The next lemma generalizes Lemmas 8.3 and 8.4 of \cite{Rieffel00}.

\begin{lemma} \label{finite approx:lemma}
For any $\varepsilon>0$ there is a finite subset
$\mathcal{J}=\bar{\mathcal{J}}$ in
$\hat{G}$, containing the class of the trivial representations,
depending only
on $\mathit{l}$ and $\varepsilon$, such that for any strongly continuous
action $\alpha$ on a complete order-unit space $\bar{A}$ and for any
$a\in \bar{A}$, there is some $a'\in A_{\mathcal{J}}:=A\cap
(\bar{A}^{\Ce})_{\mathcal{J}}$ with
\begin{eqnarray*}
\pa a'\pa\le \pa a\pa, \quad L(a')\le
L(a), \quad \mbox{ and }\pa a-a'\pa\le \varepsilon L(a).
\end{eqnarray*}
\end{lemma}
\begin{proof}
 The complex conjugation is
an isometric involution invariant under $\alpha$.
By \cite[Lemma 4.4]{Li9} it suffices to show that
for any linear combination $\varphi$ of finitely many characters
    on $G$ we have $L\circ \alpha_{\varphi}\le \pa \varphi\pa_1\cdot
    L$ on $\bar{A}^{\Ce}$, where
    $\alpha_{\varphi}$ is the linear map on
    $\bar{A}^{\Ce}$ defined in Lemma~\ref{proj TVS:lemma}.
Notice that $\varphi$ is central. Then it is easy to see
  that $\alpha_{\varphi}$ is $G$-equivariant.
Thus for
  any $b\in \bar{A}^{\Ce}$ and $x\in G$ we have
\begin{eqnarray*}
\pa \alpha_{\varphi}(b)-\alpha_x(\alpha_{\varphi}(b))\pa &=&
\pa \alpha_{\varphi}(b)-\alpha_{\varphi}(\varphi_x(b))\pa \le
\pa \varphi\pa_1\cdot \pa b-\varphi_x(b)\pa \\
&\le &
\mathnormal{l}(x)\pa \varphi\pa_1L(b).
\end{eqnarray*}
Consequently, $L(\alpha_{\varphi}(b))\le \pa \varphi\pa_1L(b)$.
\end{proof}

We are ready to prove Theorem~\ref{criterion of cont field of
  action:thm}.

\begin{proof}[Proof of Theorem~\ref{criterion of cont field of
  action:thm}]
Since the conditions in Definition~\ref{cont field of CQM:def} and
(i)-(iii) in Theorem~\ref{criterion of cont field of
  action:thm} are all local statements, we may assume that $T$ is
compact. By Lemmas~\ref{ls of mul:lemma} and \ref{cont of
Lip:lemma} the set of $G$-finite elements in $A_t$ is contained in
$\Gamma^L_t$. Lemma~\ref{finite approx:lemma} tells us that the
restriction of $L_t$ on the set of $G$-finite elements determines
the whole of $L^c_t=L_t$. Thus the induced field $(\{(A_t, L_t)\},
\Gamma)$ is a continuous field of compact quantum metric spaces.
(ii)$\Longrightarrow$(i) follows from Lemma~\ref{ls of mul:lemma}.
Let $R>\int_G\mathnormal{l}(x)\, dx$. By Theorem~\ref{finite to
  Lip:thm}, $R\ge r_{A_t}$ for all $t\in T$. We will pick
a special sequence $f_n$ in $\Gamma$ in order to apply Theorem~\ref{criterion 2 of
  conv:thm}.

As indicated in Remark~\ref{length func:remark},
$G$ is metrizable and
hence $L^2(G)$ is separable.
Since every $\gamma\in \hat{G}$ appears in the left regular
 representation, $\hat{G}$ is countable. Then we can take an
increasing sequence of finite subsets $\mathcal{J}_1\subseteq
\mathcal{J}_2\subseteq \cdots
$ of $\hat{G}$ such that $\gamma_0\in \mathcal{J}_1$,
$\cup^{\infty}_{k=1}\mathcal{J}_k=\hat{G}$, and
$\overline{\mathcal{J}_k}=\mathcal{J}_k$ for all $k$. Clearly
$a+ia'\in (\overline{A_t}^{\Ce})_{\mathcal{J}_k}$ if and only
$a-ia'\in (\overline{A_t}^{\Ce})_{\mathcal{J}_k}$. Let
$(A_t)_{\mathcal{J}_k}=A_t \cap
(\overline{A_t}^{\Ce})_{\mathcal{J}_k}$. Then
$(A_t)_{\mathcal{J}_k}$ spans
$(\overline{A_t}^{\Ce})_{\mathcal{J}_k}$ as a complex vector
space. Let $m_k=\dim((\overline{A_{t_0}}^{\Ce})_{\mathcal{J}_k})$.
Let $A^{\fin}_t$ and $\Gamma^{\fin}$ be the set of $G$-finite
elements
 in $A_t$ and $\Gamma$ respectively.
Then we can find a basis $(b_1, b_2, \cdots )$ of $A^{\fin}_{t_0}$
 such that $b_1=e_{A_{t_0}}$ and $(b_1, {\cdots}, b_{m_k})$ is a basis of
 $(A_{t_0})_{\mathcal{J}_k}$ for all $k$. Let $\pi_t$ be the projection
$\Gamma\to \overline{A_t}$.

\begin{lemma} \label{G-lift:lemma}
There exists a $G$-equivariant linear (probably unbounded) map $\varphi:A^{\fin}_{t_0}\to
\Gamma^{\fin}$ such that $\varphi(e_{A_{t_0}})=e$ and $\varphi$ is a
right inverse of $\pi_{t_0}$, \ie  $\pi_{t_0}\circ
\varphi$ is the identity map on $A^{\fin}_{t_0}$.
\end{lemma}
\begin{proof}
Let $\phi:A^{\fin}_{t_0}\to \Gamma$ be a linear right inverse map
of $\pi_{t_0}$ with $\phi(e_{A_{t_0}})=e$. For any $a\in
A^{\fin}_{t_0}$ its $G$-orbit $\{(\alpha_{t_0})_x(b):x\in G\}$ is
contained in a finite dimensional subspace of $A^{\fin}_{t_0}$.
Thus $\int_G (\alpha_x\circ \phi\circ
(\alpha_{t_0})_{x^{-1}})(a)\, dx$ makes sense. It is standard
\cite[page 77]{Bt95} that $\varphi=\int_G \alpha_x\circ \phi\circ
(\alpha_{t_0})_{x^{-1}}\, dx$ is a $G$-equivariant linear map from
$A^{\fin}_{t_0}$ to $\Gamma$. Clearly $\varphi(e_{A_{t_0}})=e$. For
any $b\in A^{\fin}_{t_0}$ we have
\begin{eqnarray*}
\pi_{t_0}(\varphi(b))&=&
\pi_{t_0}(\int_G (\alpha_x\circ \phi\circ(\alpha_{t_0})_{x^{-1}})(b)\,dx)
=\int_G (\pi_{t_0}\circ \alpha_x\circ\phi\circ(\alpha_{t_0})_{x^{-1}})(b)\,dx \\
&=&\int_G ((\alpha_{t_0})_x\circ \pi_{t_0}\circ \phi\circ(\alpha_{t_0})_{x^{-1}})(b)\,dx
=\int_G ((\alpha_{t_0})_x\circ (\alpha_{t_0})_{x^{-1}})(b)\,dx \\
&=& b.
\end{eqnarray*}
Thus $\varphi$ is also a right inverse of $\pi_{t_0}$.
\end{proof}

Take a dense sequence $\{a_n:n\in \Ne\}$
in $\cD_{ R}(A^{\fin}_{t_0})$
 such that $\{a_n:n\in \Ne\}\cap
\cD_{ R}((A_{t_0})_{\mathcal{J}_k})$ is dense in
$\cD_{ R}((A_{t_0})_{\mathcal{J}_k})$ for all $k$. By
Lemma~\ref{finite
  approx:lemma}
the set $\cD_{ R}(A^{\fin}_{t_0})$ is dense in
$\cD_{ R}(A_{t_0})$. Consequently, so is $\{a_n:n\in
\Ne\}$. For each $n$ let $f_n=\varphi(a_n)$.
 Then $f_n\in \Gamma$
and $(f_n)_{t_0}=a_n$. Let $\varphi_t=\pi_t\circ \varphi$. Let
$g_n=\varphi(b_n)$, and let
$(V_t)_k=\varphi_t((A_{t_0})_{\mathcal{J}_k})$. Since $\varphi_t$
is $G$-equivariant, $(V_t)_k$ is contained in
$(A_t)_{\mathcal{J}_k}$. Now we apply Theorem~\ref{criterion 2 of
  conv:thm} to these $f_n$'s.

Suppose that (i) holds. Let $\varepsilon>0$ be given. By
Lemma~\ref{finite approx:lemma} there is some $k$ such that
$\cD_{ R}((A_t)_{\mathcal{J}_k})$ is
$\frac{\varepsilon}{2}$-dense in $\cD_{ R}(A_t)$ for all
$t\in T$. Then there is a neighborhood $\mathcal{U}_1$ of $t_0$
such that $\dim((\overline{A_t}^{\Ce})_{\mathcal{J}_k})=
\dim((\overline{A_{t_0}}^{\Ce})_{\mathcal{J}_k})$ and hence
$\dim((A_t)_{\mathcal{J}_k})=\dim((A_{t_0})_{\mathcal{J}_k})$ for
all $t\in \mathcal{U}_1$. According to Lemma~\ref{Ind open:lemma}
there is some compact neighborhood $\mathcal{U}_2\subseteq
\mathcal{U}_1$ of $t_0$ such that $(g_1)_t,{\cdots}, (g_{m_k})_t$
are linearly independent in $A_t$ for all $t\in \mathcal{U}_2$.
Then $(g_1)_t, {\cdots}, (g_{m_k})_t$ is a basis of
$(A_t)_{\mathcal{J}_k}$ for all $t\in \mathcal{U}_2$. Take a real
vector space $V$ of dimension $m_k$ with a fixed basis $v_1,
{\cdots}, v_{m_k}$. For each $t\in \mathcal{U}_2$ let
$\psi_t:V\rightarrow (A_t)_{\mathcal{J}_k}$ be the linear
isomorphism determined by $\varphi_t(v_j)=(g_j)_t$ for all $1\le
j\le m_k$. Then $V$ gets an order-unit space structure and a
Lip-norm for each $t\in\mathcal{U}_2$
 by identifying $V$ and $(A_t)_{\mathcal{J}_k}$ via $\psi_t$.
Lemma~\ref{cont of Lip:lemma} tells us that this is a continuous
field of finite-dimensional order-unit spaces and Lip-norms as
defined in Example~\ref{finite dim cont field:eg}.
Let
$\{f_{n_s}\}_{s\in \Ne}$ be the subsequence of
$\{f_n\}_{n\in \Ne}$ whose image under $\pi_{t_0}$ is contained in
$(A_{t_0})_{\mathcal{J}_k}$. Then
$\{f_{n_s}\}_{s\in \Ne}$ is dense in $\cD_{
R}((A_{t_0})_{\mathcal{J}_k})$, and we can apply Lemma~\ref{unif
to
  cover:lemma} to
$\{f_{n_s}\}_{s\in \Ne}$.
  So
there are a neighborhood $\mathcal{U}\subseteq \mathcal{U}_2$ of
$t_0$ and some $S\in \Ne$ such that the open
$\frac{\varepsilon}{2}$-balls in $(A_t)_{\mathcal{J}_k}$ centered at
$(f_{n_1})_t, {\cdots}, (f_{n_S})_t$ cover $\cD_{
R}((A_t)_{\mathcal{J}_k})$ for all $t\in \mathcal{U}$.
Consequently, the open $\varepsilon$-balls in $A_t$ centered at
$(f_1)_t, {\cdots}, (f_{n_S})_t$ cover $\cD_{ R}(A_t)$ for
all $t\in \mathcal{U}$. By Theorems~\ref{criterion 2 of
  conv:thm} and \ref{dist_q=dist_oq:thm} we
get (iii).

We proceed to show that (iii)$\Longrightarrow$(ii). Suppose that
(iii) holds. Let $\gamma\in \hat{G}$. Say, $\gamma\in
\mathcal{J}_k$. Let
$(\alpha_t)_{\mathcal{J}_k}:\overline{A_t}^{\Ce}\to
(\overline{A_t}^{\Ce})_{\mathcal{J}_k}$ be the continuous map
defined in Lemma~\ref{proj TVS:lemma}. Then $\pa
(\alpha_t)_{\mathcal{J}_k}\pa \le M_1:= \sum_{\gamma'\in
\mathcal{J}_k}\dim(\gamma')\pa \chi_{\overline{\gamma'}}\pa_1$.
Since $\mathcal{J}_k=\overline{\mathcal{J}_k}$,
$(\alpha_t)_{\mathcal{J}_k}$ maps $\overline{A_t}$ into
$(A_t)_{\mathcal{J}_k}$. Let $\varepsilon$ be a positive number which
we shall choose later. By Theorems~\ref{criterion 2 of
  conv:thm} and \ref{dist_q=dist_oq:thm} there is a neighborhood
  $\mathcal{U}$ of $t_0$ and some $N$
such that the open $\varepsilon$-balls in $A_t$ centered at $(f_n)_t,
n=1, {\cdots}, N$, cover $\cD_{ R}(A_t)$ for all $t\in
\mathcal{U}$. Then the open $M_1\varepsilon$-balls in $A_t$ centered
at $(\alpha_t)_{\mathcal{J}_k}((f_n)_t), n=1, {\cdots}, N$, cover
$\cD_{ R}((A_t)_{\mathcal{J}_k})$ for all $t\in
\mathcal{U}$. Notice that
$(\alpha_t)_{\mathcal{J}_k}((f_n)_t)=\varphi_t((\alpha_{t_0})_{\mathcal{J}_k}(a_n))$
is contained in $(V_t)_k$ for all $n$.

Suppose that $\mul(\overline{A_t}^{\Ce}, \gamma)>
\mul(\overline{A_{t_0}}^{\Ce}, \gamma)$ for some $t$. Since
$\varphi_t$ is $G$-equivariant, we have $\mul((V_t)^{\Ce},
\gamma)\le \mul(\overline{A_{t_0}}^{\Ce}, \gamma)$. So we can find
a $u$ in $\overline{A_t}^{\Ce}\setminus (V_t)^{\Ce}$ such that the
complex linear span of $\{(\alpha_t)_x(u):x\in G\}$, the $G$-orbit
of $u$, is irreducible of type $\gamma$. Say $u=u'+iu''$ with $u',
u''\in (A_{t_j})_{\mathcal{J}_k}$. Let $W$ be the sum of $(V_t)_k$
and the real linear span of the $G$-orbits of $u'$ and $u''$. Then
$W\varsupsetneqq (V_t)_k$. Clearly $\mul(W^{\Ce}, \gamma')=0$ for
all $\gamma' \in \hat{G}\setminus \mathcal{J}_k$ and $\mul(W^{\Ce},
\gamma')\le \mul(\overline{A_{t_0}}^{\Ce}, \gamma')+2$ for all
$\gamma'$ in $\mathcal{J}_k$. Let $\omega:\mathcal{J}_k\to \Ne$ be
the function $\omega(\gamma')= \mul(\overline{A_{t_0}}^{\Ce},
\gamma')+2$ for all $\gamma'$ in $\mathcal{J}_k$. Let $M_2$ be the
constant $M_{\mathcal{J}_k, \omega}$ in Lemma~\ref{bound Lip for
  omega:lemma}. Then
$L_t\le M_2\pa \cdot \pa$ on $W^{\Ce}$. Pick a vector in
$W/(V_t)_k$ with norm $\min(\frac{1}{M_2}, R)$ and lift it up to a
vector $v$ in $W$ with the same norm. Then $\pa v-a\pa\ge
\min(\frac{1}{M_2}, R)$ for all $a\in (V_t)_k$ and $L_t(v)\le
M_2\pa v\pa \le 1$. So $v\in \cD_{
  R}(A_t)$.
Thus if we choose $\varepsilon$ small enough so
that $\min(\frac{1}{M_2}, R)>M_1\cdot \varepsilon$, then
$\mul(\overline{A_t}^{\Ce}, \gamma)\le
\mul(\overline{A_{t_0}}^{\Ce}, \gamma)$ throughout $\mathcal{U}$.
This completes our proof of Theorem~\ref{criterion of cont field of
  action:thm}.
\end{proof}

\begin{remark} \label{proof of criterion:remark}
Based on Lemmas~\ref{cont of Lip:lemma} and ~\ref{finite
  approx:lemma}, one can also prove (i)$\Longrightarrow$(iii) along
the lines Rieffel used to prove the continuity of quantum tori
\cite[Theorem 9.2]{Rieffel00}.
\end{remark}

\begin{example}[Quantum Tori] \label{tori 2:eg}
Fix $n\ge 2$. Denote by $\Theta$ the space of all real skew-symmetric
$n\times n$ matrices. For $\theta \in \Theta$, let
$\mathcal{A}_{\theta}$ be the corresponding  quantum torus
\cite{Rieffel88, Rieffel90}.  It could be described as follows.
Let $\omega_{\theta}$ denote the skew bicharacter on $\Ze^n$
defined by
\begin{eqnarray*}
\omega_{\theta}(p, q)=e^{i\pi p\cdot \theta q}.
\end{eqnarray*}
For each $p\in \Ze^n$ there is a unitary $u_p$ in
$\mathcal{A}_{\theta}$. And $\mathcal{A}_{\theta}$ is generated by
these unitaries with the relation
\begin{eqnarray*}
u_pu_q=\omega_{\theta}(p, q)u_{p+q}.
\end{eqnarray*}
So one may think of vectors in $\mathcal{A}_{\theta}$ as some kind of
functions on $\Ze^n$ .  The
$n$-torus $\mathbb{T}^n$ has a canonical ergodic action
$\alpha_{\theta}$ on $\mathcal{A}_{\theta}$. Notice that $\Ze^n$ is
the dual group of $\mathbb{T}^n$. We denote the duality by $\left<p,
x\right>$ for $x\in \mathbb{T}^n$ and $p\in \Ze^n$. Then
$\alpha_{\theta}$ is determined by
\begin{eqnarray*}
\alpha_{\theta, x}(u_p)=\left<p, x\right>u_p.
\end{eqnarray*}
Fix a length function on $G=\mathbb{T}^n$.  Let $L_{\theta}$ and
$A_{\theta}$ be as in Theorem~\ref{finite to Lip:thm}
for the order-unit space $((\mathcal{A}_{\theta})_{\sa}, e_{\mathcal{A}_{\theta}})$.
Then $(A_{\theta}, L_{\theta})$ is a compact quantum metric space.
Rieffel showed that for each $\theta_0 \in \Theta$ we have
$\dist_{\q}(A_{\theta}, A_{\theta_0})\to 0$ as $\theta\to \theta_0$
\cite[Theorem 9.2]{Rieffel00}.
Here we give a new proof using Theorems~\ref{criterion of cont
field of action:thm} and \ref{dist_q=dist_oq:thm}.
By \cite[Corollary 2.8]{Rieffel89b}
 the sections $\theta\mapsto u_p$, where $p$
runs through $\Ze^n$, generate a continuous field of
$C^*$-algebras $(\{\mathcal{A}_{\theta}\}, \Gamma)$ over $\Theta$.
Notice that for any $x\in \mathbb{T}^n$ and $p\in \Ze^n$ the
section $\theta\mapsto \alpha_{\theta, x}(u_p)=\left<p,
  x\right>u_p$ is also in $\Gamma$. Then it is easy to see that
$\{\alpha_{\theta}\}$ is a continuous field of strongly continuous
ergodic actions. For each $p\in \Ze^n=\widehat{\mathbb{T}^n}$ and
$\theta\in \Theta$, the multiplicity of $p$ in $\alpha_{\theta}$
is one. Then Theorems~\ref{criterion of cont field of action:thm}
and \ref{dist_q=dist_oq:thm} imply that $\dist_{\q}(A_{\theta},
A_{\theta_0})\to 0$ as $\theta\to \theta_0$ for all $\theta_0\in
\Theta$.
\end{example}

\begin{example}[Integral Coadjoint Orbits] \label{Berezin 2:eg}
Let $G$ be a compact connected Lie group with a fixed length
function. Choose a maximal torus of $G$ and a Cartan-Weyl basis
of the complexification of the Lie algebra of $G$. Then there are
bijective correspondences between equivalence classes of
irreducible unitary representations of $G$, dominant weights, and integral coadjoint
orbits of $G$ \cite{Bt95} \cite[Section IV]{Landsman98a}.
Let $\mathcal{O}_{\lambda}$ be an integral coadjoint
orbit corresponding to a dominant weight $\lambda$.
Then the restriction of the coadjoint action of $G$
on $\mathcal{O}_{\lambda}$ is transitive and hence the induced
action $\alpha_0$ on $C(\mathcal{O}_{\lambda})$, the algebra of
$\Ce$-valued continuous functions on $\mathcal{O}_{\lambda}$,
is ergodic.
So we have the compact quantum metric space $(A_0, L_0)$ defined
as in Theorem~\ref{finite to Lip:thm}.  Also let $H_n$ be the
carrier space of the irreducible representation of $G$ with
highest weight $n\lambda$. Then the conjugate action
$\alpha_{1/n}$ of $G$ on $B(H_n)$, the algebra of bounded
operators on $H_n$, is ergodic. Let $(A_{1/n}, L_{1/n})$ be the
corresponding compact quantum metric space defined as in
Theorem~\ref{finite to Lip:thm}.
Using the Berezin quantization, Rieffel proved that when $G$ is
semisimple, $\dist_{\q}(A_{1/n}, A_0)\to 0$ as $n\to \infty$
\cite[Theorem 3.2]{Rieffel01}. This means that the matrix algebras
$B(H_n)$ converges to the coadjoint orbit $\mathcal{O}_{\lambda}$
as $n\to \infty$.
Here we give a new proof using Theorems~\ref{criterion of cont
field of action:thm} and \ref{dist_q=dist_oq:thm}. Let $P_n$ be
the rank-one projection of $B(H_n)$ corresponding to the highest
weight $n\lambda$. For any $a\in B(H_n)$ its Berezin covariant
symbol \cite{Berezin75, Perelomov86}, $\sigma_a$, is defined by
\begin{eqnarray*}
\sigma_a(x)=\tr(a\alpha_{1/n, x}(P_n)),
\end{eqnarray*}
where $x\in G$ and $\tr$ denotes the usual (un-normalized) trace on
$B(H_n)$. There is a natural $G$-equivariant homeomorphism from
the orbit $GP_n$ of $P_n$ (in the projective space)
under $\alpha_{1/n}$ onto the coadjoint orbit
$\mathcal{O}_{n\lambda}$
\cite[Proposition 4]{Landsman98a}. Dividing everything in $\mathcal{O}_{n\lambda}$
by $n$, we may identify $GP_n$ with $\mathcal{O}_{\lambda}$.
It is evident that $\sigma_a$ could be
viewed as a continuous function on $GP_n=\mathcal{O}_{\lambda}$.
One can check easily that $a\mapsto \sigma_a$ gives a
unital, completely
positive, $G$-equivariant linear map $\sigma_n$ from $B(H_n)$ to
$C(\mathcal{O}_{\lambda})$. Endow $\mathcal{O}_{\lambda}$
with the image of the Haar measure on $G$, which is a probability
measure invariant under $\alpha_0$. Then
$C(\mathcal{O}_{\lambda})$ has an inner product as usual.
Clearly this inner product is invariant under $\alpha_0$. Using
the normalized trace on $B(H_n)$, which is invariant under
$\alpha_{\frac{1}{n}}$, $B(H_n)$ has the Hilbert-Schmidt inner
product. Then $\sigma_n$ has an adjoint operator,
$\hat{\sigma}_n$, from $C(\mathcal{O}_{\lambda})$ to
$B(H_n)$. For any $a\in B(H_n)$ a function $f\in
C(\mathcal{O}_{\lambda})$ with $\hat{\sigma}_n(f)=a$ is
called a Berezin contravariant symbol \cite{Berezin75,
Perelomov86} for $a$. It is easy to see that $\hat{\sigma}_n$ is
unital, completely
positive and $G$-equivariant.
Since unital completely positive maps are norm-nonincreasing \cite[Lemma 5.3]{Lance95},
$\hat{\sigma}_n$ is norm-nonincreasing.
In \cite{Landsman98a}
Landsman proved that the sections given by these
$\hat{\sigma}_n(f)$'s, where $f$ runs through
$C(\mathcal{O}_{\lambda})$, generate a continuous field of
$C^*$-algebras over $T'=\{\frac{1}{n}:n\in \Ne\}\cup \{0\}$ with
fibres $B(H_n)$ at $\frac{1}{n}$ and $C(\mathcal{O}_{\lambda})$ at
$0$. In fact, Landsman proved that this is a strict quantization
of the canonical symplectic structure on $\mathcal{O}_{\lambda}$,
though we do not need this fact here. Using the fact that
$\hat{\sigma}_n$ is $G$-equivariant and norm-nonincreasing, it is
easy to check that the $\alpha_{\frac{1}{n}}$'s and $\alpha_0$ are
a continuous field of strongly continuous ergodic actions of $G$.
When $G$ is semisimple, it is known that the maps $\sigma_n$ are
all injective \cite{Perelomov86} \cite[Lemma A.2.1]{Simon80} \cite[Theorem 3.1]{Rieffel01}.
Then for each $\gamma\in \hat{G}$ we see that
$\mul(B(H_n), \gamma)\le \mul(C(\mathcal{O}_{\lambda}), \gamma)$. So
Theorem~\ref{criterion of cont field of action:thm} and
\ref{dist_q=dist_oq:thm} tell us that $\dist_{\q}(A_{1/n}, A_0)\to
0$ as $n\to \infty$.
\end{example}

\section{Continuity of $\theta$-deformations}
\label{Contdeform:sec}

In this section we prove Theorem~\ref{theta-deform cont:thm}.

We use the notation in \cite[Sections 3 and 5]{Li9}.
Let us explain first some convention used in the statement of
Theorem~\ref{theta-deform cont:thm}. $C^{\infty}(M_{\theta})$ is a
locally convex $*$-algebra, and has a natural $*$-homomorphism
$\Psi_{\theta}$ into $C(M_{\theta})$ (see the discussion after
Definition 3.9 in \cite{Li9}). Let $W_{\theta}$ be the image of
$(C^{\infty}(M_{\theta}))_{\sa}$ under the map $\Psi_{\theta}$.
\cite[Theorem 1.1]{Li9} tells us that $(C(M_{\theta}),
L_{\theta})$ is a $C^*$-algebraic compact quantum metric space.
This means that $(W_{\theta}, L_{\theta}|_{W_{\theta}})$ is a
compact quantum metric space. Since the map $\Psi_{\theta}$ is
injective \cite[Lemma 3.10]{Li9}, we may identify
$C^{\infty}(M_{\theta})$ with its image
$\Psi_{\theta}(C^{\infty}(M_{\theta}))$. In this way
$(C^{\infty}(M_{\theta}))_{\sa}$ is identified with $W_{\theta}$,
and $L_{\theta}$ has a restriction on $C^{\infty}(M_{\theta})$,
which we still denote by $L_{\theta}$ in the statement of
Theorem~\ref{theta-deform cont:thm}. In order to make the argument
clear, in the rest of this section we shall still distinguish
$C^{\infty}(M_{\theta})$ (resp. $(C^{\infty}(M_{\theta}))_{\sa}$)
and $\Psi_{\theta}(C^{\infty}(M_{\theta}))$ (resp. $W_{\theta}$).

Let $(\{\mathcal{A}_{\theta}\}, \Gamma)$ be the continuous field of
$C^*$-algebras over $\Theta$ in Example~\ref{tori 2:eg}.
We shall also see later, after Lemma~\ref{subtrivial:lemma}, that
the elements in $C(M, Cl^{\Ce}M)\otimes_{alg}\Gamma$ generate
a continuous field of $C^*$-algebras \cite{Dixmier77} over $\Theta$ with fibres
$C(M, Cl^{\Ce}M)\otimes \mathcal{A}_{\theta}$. Let $(\{C(M_{\theta})\},
\Gamma^M)$ be the subfield with fibers
$C(M_{\theta})$.

\begin{lemma} \label{subtrivial:lemma}
There exist a $C^*$-algebra $\mathcal{B}$ and faithful $*$-homomorphisms
$\phi_{\theta}:\mathcal{A}_{\theta}\hookrightarrow \mathcal{B}$
such that for every $f\in \Gamma$,
the $\mathcal{B}$-valued function $\theta \mapsto \phi_{\theta}(f_{\theta})$
is continuous over $\Theta$.
\end{lemma}
\begin{proof}
Every unital $C^*$-algebra admitting an ergodic action of
$\mathbb{T}^n$ is nuclear \cite[Lemma 6.2]{OPT80} \cite[Proposition 3.1]{DLRZ02}.
Thus  $\mathcal{A}_{\theta}$ is
nuclear.

Notice that $\mathcal{A}_{\theta}$ is isomorphic to
$\mathcal{A}_{\theta+M}$ naturally for any skew-symmetric $n\times
n$ matrix $M$ with even integer entries, by identifying the
corresponding $u_q$. So we may think of $(\{\mathcal{A}_{\theta}\},
\Gamma)$ as a continuous field over the quotient space of $\Theta$
by all the skew-symmetric $n\times n$ matrices with even integer
entries. This quotient space is just a torus of dimension
$\frac{n(n-1)}{2}$. Now our assertion follows from the result of
Blanchard \cite[Theorem 3.2]{Blanchard98} that every separable unital
continuous field of nuclear $C^*$-algebras over a
compact metric space has a faithful $*$-homomorphism from each
fibre into the Cuntz algebra $\mathcal{O}_2$
such that the global
continuous sections become continuous paths in $\mathcal{O}_2$.
\end{proof}

Via identifying $\mathcal{A}_{\theta}$ with
$\phi_{\theta}(\mathcal{A}_{\theta})$ we see that the continuous field
$(\{\mathcal{A}_{\theta}\}, \Gamma)$ becomes a subfield of the trivial
field over $\Theta$ with fibre $\mathcal{B}$. Then
the elements of $C(M, Cl^{\Ce}M)\otimes_{alg}\Gamma$ are continuous
sections of the trivial field with fibre \\$C(M,
Cl^{\Ce}M)\otimes \mathcal{B}$. So they generate a
subfield
of the trivial field with fibres $ C(M,
Cl^{\Ce}M)\otimes \mathcal{A}_{\theta}$.

Now we need to distinguish the norms for elements of the form
$\sum^k_{j=1}y_{q_j}\otimes u_{q_j}$ at different $\theta$. For this
we let $\pa \cdot \pa_{\theta}$ denote the norm of
$C(M, Cl^{\Ce}M)\otimes \mathcal{A}_{\theta}$.

Let $\mathcal{J}$ be a finite subset of
$\Ze^n=\widehat{\mathbb{T}^n}$ such that
$\mathcal{J}=\bar{\mathcal{J}}$ and $\gamma_0\in \mathcal{J}$. For
any $q\in \Ze^n$ let $(C(M_{\theta}))_q$ be the $q$-isotypic
component of $C(M_{\theta})$ under the action $\alpha=I\otimes
\tau$, and let $(C^{\infty}(M))_q$ be the $q$-isotypic component
of $C^{\infty}(M)$ under the action $\sigma$ as in \cite[Section
6]{Li9}, where $\tau$ and $\sigma$ are the actions of
$\mathbb{T}^n$ on $\mathcal{A}_{\theta}$ and $C(M)$ respectively.
Similarly we define
$(C(M_{\theta}))_{\mathcal{J}}$ and $(C^{\infty}(M))_{\mathcal{J}}$.
By \cite[Lemma 6.2]{Li9}
we have
\begin{eqnarray*}
(C(M_{\theta}))_{\mathcal{J}}\cap W_{\theta}=
(\sum_{q\in
  \mathcal{J}}(C^{\infty}(M))_q\otimes u_q)_{\sa}.
\end{eqnarray*}
Let $V=(\sum_{q\in
  \mathcal{J}}(C^{\infty}(M))_q\otimes u_q)_{\sa}$, and
let $e=1_M\otimes u_{\gamma_0}$. Then $(V, e)$ gets an order-unit space
structure from $((C(M_{\theta}))_{\sa}, e)$. Clearly the restriction of
$\pa \cdot \pa_{\theta}$ is exactly the order-unit norm.
Denote by $(V_{\theta}, e)$ this order-unit space.
For each $\theta\in \Theta$, by Proposition~\ref{criterion of
  Lip:prop}
the restriction of $L_{\theta}$ to $V_{\theta}$ is a Lip-norm with
$r_{V_{\theta}}\le r_{W_{\theta}}$.

\begin{lemma} \label{unif:lemma}
$(V, e, \{\pa \cdot\pa_{\theta}|_V\}, \{L_{\theta}|_V\})$ is a
uniformly continuous field of order-unit spaces and Lip-norms over
$\Theta$ (see Definition~\ref{unif:def}). For any $v\in V$ the
function $\theta \mapsto L_{\theta}(v)$ is continuous over
$\Theta$.
\end{lemma}
\begin{proof}
Let $v\in V$. Say $v=\sum_{q\in \mathcal{J}}v_q\otimes u_q$.
By \cite[Corollary 5.7]{Li9} we have
\begin{eqnarray} \label{Ltheta to LD:eq}
L_{\theta}=L^D
\end{eqnarray}
on $C(M_{\theta})$, where $L^D$ was defined in \cite[Definition
5.3]{Li9}. Recall that for any $f\in C^{\infty}(M)$ we have
\cite[Lemma II.5.5]{LM89}
\begin{eqnarray} \label{D df:eq}
[D, f]=df  \mbox{ as linear maps on } C^{\infty}(M, S),
\end{eqnarray}
where $df\in C^{\infty}(M, T^*M^{\Ce})\subseteq C^{\infty}(M, Cl^{\Ce}M)$
acts on $C^{\infty}(M, S)$ via the left $C^{\infty}(M,Cl^{\Ce}M)$-module
structure of $C^{\infty}(M, S)$.
Then
\begin{eqnarray} \label{L to d 2:eq}
L_{\theta}(v)
&=& L_{\theta}(\sum_{q\in \mathcal{J}}v_q\otimes u_q)
\overset{(\ref{Ltheta to LD:eq})}
= L^D(\sum_{q\in \mathcal{J}}v_q\otimes
u_q)
= \pa [D^{L^2},\sum_{q\in \mathcal{J}}v_q\otimes
u_q]\pa_{\theta} \nonumber \\
&=& \pa  \sum_{q\in \mathcal{J}}[D, v_q]\otimes
u_q\pa_{\theta}
\overset{(\ref{D df:eq})}= \pa  \sum_{q\in \mathcal{J}}(dv_q)\otimes
u_q\pa_{\theta}.
\end{eqnarray}
Therefore the function $\theta \mapsto L_{\theta}(v)$ is
continuous over $\Theta$.
As in the proof of \cite[Lemma 4.6]{Li9}
we have that
\begin{eqnarray*}
L_{\theta}(v_q\otimes u_q)\le L_{\theta}(v) \quad \mbox{ and }
\quad \pa v_q\otimes u_q\pa_{\theta} \le \pa v\pa_{\theta}
\end{eqnarray*}
for all $q\in \Ze^n$.

Let $\theta_0\in \Theta$, and let $\varepsilon>0$ be given.
Since for each $q$ the map $\theta\mapsto \phi_{\theta}(u_q)$ from
$\Theta$ to $\mathcal{B}$ is continuous,
there is some neighborhood $\mathcal{U}$ of $\theta_0$ such that
$\sum_{q\in \mathcal{J}}\pa \phi_{\theta}(u_q)-\phi_{\theta_0}(u_q)\pa
<\varepsilon$ throughout $\mathcal{U}$. Let $\theta\in \mathcal{U}$.
Then for any $z_q$'s in $C(M, Cl^{\Ce}M)$ with $q\in \mathcal{J}$ and $\pa z_q\pa \le 1$ we have
\begin{eqnarray*}
& &\pa (I\otimes \phi_{\theta})
(\sum_{q\in \mathcal{J}}z_q\otimes u_q) -
(I\otimes \phi_{\theta_0})
(\sum_{q\in \mathcal{J}}z_q\otimes u_q)\pa \\
&\le &
\sum_{q\in \mathcal{J}}\pa z_q\pa \cdot \pa
\phi_{\theta}(u_q)-\phi_{\theta_0}(u_q)\pa \le
\varepsilon.
\end{eqnarray*}
Suppose that $\pa v\pa_{\theta_0}=1$ for some $v\in V$. Say
$v=\sum_{q\in \mathcal{J}}v_q\otimes u_q$. Then $\pa v_q\pa=\pa
v_q\otimes u_q\pa_{\theta_0}\le 1$ for each $q\in \mathcal{J}$.
Thus
\begin{eqnarray*}
\pa v\pa_{\theta}&=&\pa (I\otimes \phi_{\theta})(v)\pa\le \pa (I\otimes \phi_{\theta})
(v)-(I\otimes \phi_{\theta_0})(v)\pa +\pa (I\otimes
\phi_{\theta_0})(v)\pa \\ &\le &
\varepsilon+\pa v\pa_{\theta_0}=
\varepsilon+1.
\end{eqnarray*}
Similarly, $\pa v\pa_{\theta}\ge 1-\varepsilon$. Therefore
$(1-\varepsilon)\pa \cdot \pa_{\theta_0}\le \pa \cdot
\pa_{\theta}\le (1+\varepsilon)\pa \cdot \pa_{\theta_0}$ on $V$ throughout
$\mathcal{U}$. Now suppose that $L_{\theta_0}(w)=1$ for some $w\in
V$. Say $w=\sum_{q\in \mathcal{J}}w_q\otimes u_q$.
Then by (\ref{L to d 2:eq}) $\pa
dw_q\pa=\pa dv_q\otimes u_q\pa_{\theta_0}=
L_{\theta_0}(v_q\otimes u_q) \le 1$. Let $w'=\sum_{q\in \mathcal{J}}(dw_q)\otimes
u_q\in \sum_{q\in \mathcal{J}}C(M, Cl^{\Ce}M)\otimes u_q$. Then
\begin{eqnarray*}
L_{\theta}(w)&\overset{(\ref{L to d 2:eq})}=&
\pa  w'\pa_{\theta}=\pa (I\otimes \phi_{\theta})(w')\pa\\
&\ge &\pa (I\otimes
\phi_{\theta_0})(w')\pa -\pa (I\otimes \phi_{\theta_0})(w')-(I\otimes \phi_{\theta})(w')\pa \\
&\overset{(\ref{L to d 2:eq})}\ge & L_{\theta_0}(v)-\varepsilon= 1-\varepsilon.
\end{eqnarray*}
Therefore
$(1-\varepsilon)L_{\theta_0}\le L_{\theta}$ on $V$ throughout
$\mathcal{U}$.
\end{proof}

Combining \cite[Lemma 4.4]{Li9}
and Lemma~\ref{unif:lemma}
together we see that
the field $(\{(W_{\theta}, L_{\theta}|_{W_{\theta}})\},
(\Gamma^M)_{\sa})$ is a continuous field of compact quantum metric
  spaces.
Let $R=r_M+C\int_{\mathbb{T}^n}\mathnormal{l}(x)\, dx$, where
$r_M$ is the radius of $M$ equipped with the geodesic distance and
the constant $C$ was defined in \cite[Proposition 5.5]{Li9}. At
the end of \cite{Li9}
it was proved that
the radius of $(W_{\theta},
  L_{\theta}|_{W_{\theta}})$ is no bigger than $R$ for each $\theta$.
Let $\varepsilon>0$ be given. Pick a finite subset $\mathcal{J}\subseteq \Ze^n$
for $\varepsilon$
in \cite[Lemma 4.4]{Li9}.
Then $\dist^R_{\oq}(W_{\theta},V_{\theta} )\le
\dist^{W_{\theta}}_{\rH}(\cD_{ R}(W_{\theta}), \cD_{
R}(V_{\theta}))\le \varepsilon$ for all $\theta\in \Theta$. By
Lemmas~\ref{unif:lemma}, \ref{unif to cover:lemma}, and
Theorem~\ref{criterion 2 of conv:thm} $\dist^R_{\oq}(V_{\theta},
V_{\theta_0})\to 0$ as $\theta\to \theta_0$. Thus there is a
neighborhood $\mathcal{U}$ of $\theta_0$ such that
$\dist^R_{\oq}(V_{\theta}, V_{\theta_0})<\varepsilon$ throughout
$\mathcal{U}$. Then clearly
$\dist^R_{\oq}(W_{\theta},W_{\theta_0})\le 3\varepsilon$ throughout
  $\mathcal{U}$.
Therefore $\dist^R_{\oq}(W_{\theta},W_{\theta_0})\to 0$ as
$\theta\to \theta_0$. By Theorem~\ref{dist_q=dist_oq:thm} we get
$\dist_{\oq}(W_{\theta},W_{\theta_0})\to 0$. This finishes the proof
of Theorem~\ref{theta-deform cont:thm}.

\end{document}